\documentclass[11pt,amssymb]{amsart}
\usepackage{amssymb}
\usepackage[mathscr]{eucal}

\newcommand{\Z}{{\mathbb Z}}
\newcommand{\C}{{\mathbb C}}
\newcommand{\Q}{{\mathbb Q}}
\newcommand{\tr}{\mathrm{Tr}\:}

\newcommand{\R}{{\mathbb R}}
\newcommand{\A}{\mathbb{A}}
\newcommand{\rk}{\mathrm{rk}}

\newcommand{\cK}{\mathcal{K}}

\newcommand{\fX}{\mathfrak{X}}
\newcommand{\cO}{\mathcal{O}}
\newcommand{\HH}{\mathbb{H}}
\newcommand{\cS}{\mathcal{S}}

\newcommand{\Ga}{\mathrm{Gal}}
\newtheorem{thm}{Theorem}[section]
\newtheorem{lemma}[thm]{Lemma}
\newtheorem{prop}[thm]{Proposition}
\newtheorem{cor}[thm]{Corollary}

\newcommand{\cG}{\mathcal{G}}

\begin{document}

\font\brus=wncyr10.240pk scaled 1200 \font\rus=wncyr10.240pk

\title[Arithmetic Groups and Locally Symmetric Spaces]{Weakly commensurable
arithmetic groups, lengths of closed geodesics and
isospectral locally symmetric spaces}

\author[Prasad]{Gopal Prasad}

\author[Rapinchuk]{Andrei S. Rapinchuk}

\address{Department of Mathematics, University of Michigan, Ann
Arbor, MI 48109}
\email{gprasad@umich.edu}

\address{Department of Mathematics, University of Virginia,
Charlottesville, VA 22904}
\email{asr3x@virginia.edu}
\maketitle

\centerline{\it Dedicated to the memory of A.R.'s mother Izabella B.
Rapinchuk}

\vskip5mm

\section{Introduction}

The goal of this paper is two-fold. First, we introduce and analyze
a~new relationship between (Zariski-dense) abstract subgroups of the
group of $F$-rational points of a connected semi-simple algebraic
group defined over a field $F$, which we call {\it weak
commensurability.} This relationship is expressed in terms of the
eigenvalues of individual elements, and does not involve any
structural connections between the subgroups. Nevertheless, it turns
out that weakly commensurable $S$-arithmetic subgroups of a
given group always split into finitely many commensurability
classes, and that in certain types of groups, any two weakly
commensurable $S$-arithmetic subgroups are actually commensurable.
Second, we use results and conjectures in transcendental number
theory to relate weak commensurability with interesting differential
geometric problems on length-commensurable, and isospectral, locally
symmetric spaces, and to settle a~series of open questions in this
area by applying our results on weakly commensurable arithmetic (and
more general) subgroups. These applications lead us to believe that
the notion of weak commensurability is likely to become a useful
tool in the theory of Lie groups and related areas.

\vskip3mm

We begin with the definition of weak commensurability. Let $G$ be a
connected semi-simple algebraic group defined over a field $F$.

\vskip2mm

\noindent {\bf Definitions. 1.} Two semi-simple elements $g_1 ,
g_2$ of $G(F)$ are  {\it weakly commensurable} if
there exist maximal $F$-tori $T_1 , T_2$ of $G$ such that $g_i
\in T_i(F),$ and for some characters $\chi_i$ of $T_i$ (defined over an algebraic closure $\overline F$ of $F$), we have
$$
\chi_1(g_1) = \chi_2(g_2) \neq 1. \footnote{In other
words, the subgroup of ${\overline F}^{\times}$ generated by the  eigenvalues
(in a faithful representation of $G$) of $g_1$ intersects the subgroup generated by the eigenvalues of $g_2$ nontrivially.}
$$

{\bf 2.} Two subgroups $\Gamma_1$ and $\Gamma_2$  of $G(F)$ are {\it weakly
commensurable} if given a semi-simple element $\gamma_1 \in \Gamma_1$ of infinite order, there is a semi-simple element $\gamma_2\in\Gamma_2$ of infinite order which is weakly commensurable to $\gamma_1$, and given a semi-simple
element $\gamma_2 \in \Gamma_2$ of infinite order, there is a semi-simple element $\gamma_1\in\Gamma_1$ of infinite order which is weakly
commensurable to $\gamma_2$.
\vskip1mm

Our first basic result is the following.

\vskip1mm

\noindent {\bf Theorem A.} {\it Let $G$ be a connected  absolutely simple algebraic group
defined over a field $F$ of characteristic zero. Let $\Gamma_1$ and $\Gamma_2$ be two
finitely generated Zariski-dense subgroups of $G(F)$, and
$K_{\Gamma_i}$ be the subfield of $F$ generated by the traces ${\tr}
\mathrm{Ad}\:\gamma$, in the adjoint representation,  of  $\gamma \in \Gamma_i$. If
$\Gamma_1$ and $\Gamma_2$ are weakly commensurable, then
$K_{\Gamma_1} = K_{\Gamma_2}.$}

\vskip3mm

Most of the results of this paper are on arithmetic subgroups. In
fact, the central issue for us is what can be said about two forms
over number fields, of a connected absolutely simple $F$-group $G$,
given that these forms contain weakly commensurable Zariski-dense
$S$-arithmetic subgroups. To give the precise statements (see
Theorems B--E), we need to describe our set-up more carefully. Let
$G$ be a connected  absolutely simple algebraic group defined over a
field $F$ of characteristic zero. Suppose we are given a number
field $K$, an embedding $K \hookrightarrow F$, and an algebraic
$K$-group $G_0$ such that the $F$-group $_{F}G_0$ obtained by
extension of scalars $K\hookrightarrow F$, is $F$-isomorphic to $G$
(in other words, $G_0$ is an $F/K$-form of $G$). Then we have an 
embedding $\iota \colon G_0(K) \hookrightarrow G(F),$ which is
well-defined up to an $F$-automorphism of $G.$ Now let $S$ be a
finite set of places of $K$ which contains the set $V^{K}_{\infty}$
of all archimedean places, but does not contain any nonarchimedean
place where $G_0$ is anisotropic. Let $\cO_{K}(S)$ denote the ring
of $S$-integers in $K$ (with $\cO_{K} = \cO_{K}(V^{K}_{\infty})$
denoting the ring of algebraic integers in $K$), and let
$G_0(\cO_{K}(S))$ denote the corresponding $S$-arithmetic subgroup
defined in terms of a fixed $K$-embedding $G_0 \hookrightarrow
\mathrm{GL}_n,$ i.e., $G_0(\cO_{K}(S))= G_0(K)\cap{\rm
GL}_n(\cO_{K}(S))$. We will say that two subgroups $\Gamma' ,
\Gamma''$ of $G(F)$ are commensurable up to an $F$-automorphism of
$G$ if there exists an $F$-automorphism $\sigma$ of $G$ such that
$\sigma(\Gamma')$ and $\Gamma''$ are commensurable in the usual
sense (i.e., their intersection has finite index in both of them).
Then any subgroup $\Gamma$ of $G(F)$ which is commensurable with
$\iota (G_0(\cO_{K}(S)))$ up to an $F$-automorphism of $G$, will be
called a $(G_0 , K , S)$-{\it arithmetic subgroup}\footnote{Notice
that if $G_0$ is anisotropic over $K_v,$ where $v$ is a
nonarchimedean place of $K$, then $G_0 (\cO_{K}(S))$ is
commensurable with $G_0(\cO_{K}(S \cup \{ v\})),$ so the classes of
$S$- and $(S \cup \{ v \} )$-arithmetic subgroups coincide. Thus,
the above assumption on $S$ is necessary if one wants to recover $S$
from a given $S$-arithmetic subgroup.}. As usual, $(G_0 , K,
V^{K}_{\infty})$-arithmetic subgroups will simply be called $(G_0 ,
K)$-arithmetic. Now, let $\Gamma_i \subset G(F),$ where $i = 1 , 2,$
be Zariski-dense $(G_i , K_i ,S_i)$-arithmetic subgroups. The key
question for us is when does the fact that $\Gamma_1$ and $\Gamma_2$
are weakly commensurable imply that they are commensurable up to an
$F$-automorphism of $G,$  i.e., $K_1 = K_2,$ $S_1 = S_2$ and $G_1$
and $G_2$ are $K$-isomorphic (cf.\,Proposition \ref{P:P1}). Theorems
B--E address this question.

\vskip3mm

\noindent {\bf Theorem B.} {\it If Zariski-dense $(G_i , K_i ,
S_i)$-arithmetic subgroups $\Gamma_i$ of $G(F)$ are weakly commensurable for $i = 1 ,
2,$ then $K_1 = K_2$ and $S_1 = S_2.$}

\vskip3mm

Examples 6.5 and 6.6 show that the existence of weakly commensurable
$S$-arithmetic subgroups does not guarantee that $G_1$ and $G_2$ are
always isomorphic over $K := K_1 = K_2$. In the next theorem we list
the cases where it can be asserted that $G_1$ and $G_2$ are
$K$-isomorphic, and then give a general finiteness result for the
number of $K$-isomorphism classes.

\vskip3mm

\noindent {\bf Theorem C.} {\it Suppose $G$ is not of type $A_n$ $(n
> 1),$ $D_{2n+1}$ $(n \geqslant 1)$, $D_4$ or $E_6.$ If $G(F)$ contains
Zariski-dense weakly commensurable $(G_i , K , S)$-arithmetic
subgroups $\Gamma_i$ for $i = 1 , 2,$ then $G_1 \simeq G_2$ over
$K,$ and hence $\Gamma_1$ and $\Gamma_2$ are commensurable up to
an $F$-automorphism of $G.$}

\vskip3mm

\noindent {\bf Theorem D.} {\it Let $\Gamma_1$ be a Zariski-dense
$(G_1 , K , S)$-arithmetic subgroup of $G(F).$ Then the set of
$K$-isomorphism classes of $K$-forms $G_2$ of $G$ such that $G(F)$
contains a Zariski-dense $(G_2 , K , S)$-arithmetic subgroup weakly
commensurable to $\Gamma_1$ is finite. In other words, the set of
all  $(K , S)$-arithmetic subgroups of $G(F)$ which are weakly
commensurable to a given $(K , S)$-arith-metic subgroup is a union
of finitely many commensurability classes.} \vskip3mm

 A noteworthy fact about weak commensurability
is that it has the following implication for the existence of
unipotent elements in arithmetic subgroups (even though it is
formulated entirely in terms of semi-simple ones).

\vskip3mm

\noindent {\bf Theorem E.} {\it Assume that $G(F)$ contains
Zariski-dense $(G_1 ,K , S)$- and $(G_2 , K , S)$-arithmetic
subgroups which are weakly commensurable. Then the Tits indices of
$G_1/K$ and $G_2/K$, and for every place $v$ of $K$, the Tits
indices of $G_1/K_v$ and $G_2/K_v$, are isomorphic. In particular,
$\mathrm{rk}_K\: G_1 = \mathrm{rk}_K\: G_2,$ and consequently if
$G_1$ is $K$-isotropic, then so is $G_2$.} \vskip1mm

(For a description of Tits index of a semi-simple algebraic group, see \S 7.)

\vskip3mm

The following result asserts that a
lattice which is weakly commensurable with an $S$-arithmetic group is arithmetic.

\vskip3mm

\noindent{\bf Theorem F.} {\it Let $G$ be a connected absolutely simple
algebraic group over a nondiscrete locally compact field $F$ of characteristic zero,
and let $\Gamma_1$ and $\Gamma_2$ be two Zariski-dense lattices in
$G(F).$ Assume that $\Gamma_1$ is $(K , S)$-arithmetic. If
$\Gamma_1$ and $\Gamma_2$ are weakly commensurable, then $\Gamma_2$
is also $(K , S)$-arithmetic.}

\vskip3mm

The proofs of these theorems use a variety of algebraic and
number-theoretic techniques. One of the key ingredients is a new
method for constructing elements with special properties in a given
Zariski-dense subgroup of a semi-simple algebraic group developed in
our papers \cite{PR1}-\cite{PR3} to answer questions of Y.~Benoist,
G.A.~Margulis, R.~Spatzier et al arising in geometry. This method is
described in \S 3 below in a considerably modified form required for
the proofs of Theorems A--F. Among other important ingredients of
our proofs are Tits' classification of semi-simple algebraic groups
over nonalgebraically closed fields (cf.\,\cite{Ti}), and results
on Galois cohomology of semi-simple groups over local and global
fields.  As a by-product of our argument, we obtain an almost
complete solution of the old problem whether an absolutely simple
group over a number field is determined by the set of isomorphism
classes of its maximal tori (cf.\,Theorem \ref{T:D2}). It is our
belief that the notion of weak commensurability, and the techniques
involved in its analysis,  in conjunction with the results of
\cite{PR1}-\cite{PR3}, will have numerous applications in the theory
of Lie groups, ergodic theory, and (differential) geometry. In fact,
the results on weak commensurability stated above were motivated by,
and actually enabled us to settle, some problems about the lengths
of closed geodesics in, and isospectrality of, arithmetically
defined locally symmetric spaces. We now proceed to describe these
geometric applications. \vskip1mm

For a Riemannian manifold $M,$ the {\it length spectrum}
$\mathcal{L}(M)$ (resp.,\,the {\it weak length spectrum} $L(M)$) is
defined to be the set of lengths of closed geodesics in $M$ with
multiplicities (resp.,\,without multiplicities), cf.\,\cite{LMNR}.
The following question has received considerable attention: to what
extent do $\mathcal{L}(M), $ $L(M)$, or the spectrum of the Laplace
operator, determine $M$? It turns out that all these sets are
interrelated: for example, two compact hyperbolic 2-manifolds are
isospectral\footnote{Two compact Riemannian manifolds are said to be
{\it isospectral} if their Laplace operators have the same
eigenvalues with the same multiplicities, cf.\,\S \ref{S:Sp}.} if
and only if they have the same length spectrum, cf.\,\cite{Mc}; two
hyperbolic 3-manifolds are isospectral if and only if they have the
same {\it complex-length} spectrum, cf.\,\cite{Ga}.
Furthermore, it is known that isospectral compact locally symmetric
spaces of nonpositive curvature have the same weak length spectrum,
see Theorem \ref{T:S1} below. The first examples of isospectral but
not isometric (although commensurable\footnote{Two manifolds are
called {\it commensurable} if they admit a common finite-sheeted
cover.}) compact hyperbolic 2- and 3-manifolds were given in
\cite{Vi}. Recently, in \cite{LSV}, noncommensurable isospectral
locally symmetric spaces have been constructed. On the other hand,
in 1985 Sunada \cite{Su} described a general method for producing
examples of nonisometric (but commensurable) isospectral manifolds.
A variant of Sunada's construction has been used in \cite{LMNR} to
give examples of hyperbolic manifolds with equal weak length spectra
but different volumes.  Earlier, in \cite{Re}, the same approach was
used to produce nonisometric hyperbolic 3-manifolds with equal weak
length spectra. It should be pointed out that Sunada's construction,
which is the only known general method for constructing manifolds
with the same (weak) length, or Laplace operator, spectra, {\it
always} produces commensurable manifolds (in particular, the
examples in \cite{LMNR} and \cite{Re} are commensurable).  So, the
following question was raised (cf.,\,\,for example, \cite{Re}):

\vskip2mm

\noindent $(1)$\,\,  {\it Let $M_1$ and $M_2$ be two (hyperbolic)
manifolds (of finite volume or even compact). Suppose $L(M_1) =
L(M_2).$ Are $M_1$ and $M_2$ necessarily commensurable?} \hfill

\vskip2mm

\noindent  One may generalize this question by introducing the
notion of length-commen- surability, which in particular allows us
to replace the manifolds under consideration with commensurable
ones: we say that $M_1$ and $M_2$ are {\it length-commensurable} if
$\Q \cdot L(M_1) = \Q \cdot L(M_2).$ Now, $(1)$ can be reformulated
as follows:

\vskip2mm

\noindent $(2)$\,\, {\it Suppose $M_1$ and $M_2$ are
length-commensurable. Are they commensurable?}

\vskip2mm

\noindent In \cite{Re},  an affirmative answer (to $(1)$) was given
for arithmetically defined hyperbolic 2-manifolds, and very recently
in \cite{CHLR} a similar result has been obtained for hyperbolic
3-manifolds. The results of this paper, combined with that of \cite{PR6}, \S 9, 
for groups of type $D_{2n}$, provide an affirmative answer
to $(2)$  for arithmetically defined hyperbolic
manifolds of dimensions $2n$ and $4n+3$, but a negative answer for 
hyperbolic manifolds of dimension $4n + 1$, and for complex hyperbolic manifolds. 
In fact, we analyze the problem in the general context of arithmetically defined
locally symmetric spaces. \vskip2mm

Let $G$ be a connected semi-simple real algebraic subgroup of $\mathrm{SL}_n$, 
$\cG$ be $G(\R)$
considered as a Lie group, and let $\cK$ be a maximal
compact subgroup of $\cG.$ Then $\fX = \cK \backslash \cG$ is the
symmetric space of $\cG.$ Given a discrete torsion-free subgroup
$\Gamma$ of $\cG,$ the quotient $\fX_{\Gamma} = \fX/\Gamma$ is a
locally symmetric space. We say that $\fX_{\Gamma}$ is {\it
arithmetically defined} if $\Gamma$ is an arithmetic subgroup of $\cG$ (cf.\:\cite{M}, Ch.\,IX). According to the following theorem,
length-commensurability of locally symmetric spaces is closely
related to weak commensurability of the corresponding discrete
subgroups.

\vskip3mm \noindent{\bf Theorem \ref{T:G-3}.} {\it Let
$\Gamma_1$, $\Gamma_2$ be discrete torsion-free subgroups of
$\mathcal{G}$. If $\Gamma_1$ and $\Gamma_2$ are not weakly
commensurable, then, possibly after interchanging them, the
following assertions hold:

\vskip2mm

{\rm \ (i)} \parbox[t]{11.5cm}{If $\mathrm{rk}_{\R}\: G = 1$,  then
there exists $\lambda_1 \in L(\mathfrak{X}_{\Gamma_1})$ such that
for {\rm any} $\lambda_2 \in L(\mathfrak{X}_{\Gamma_2}),$ the ratio
$\lambda_1/\lambda_2$ is irrational.}

\vskip2mm

{\rm (ii)} \parbox[t]{11.5cm}{If there exists a number field $K$
such that both $\Gamma_1$ and $\Gamma_2$ can be conjugated into
$\mathrm{SL}_n(K)$, and Schanuel's conjecture holds, then there
exists $\lambda_1 \in L(\mathfrak{X}_{\Gamma_1})$ which is
algebraically independent from {\rm any} $\lambda_2 \in
L(\mathfrak{X}_{\Gamma_2}).$}

\vskip2mm

\noindent In either case, (under the above assumptions)
$\fX_{\Gamma_1}$ and $\fX_{\Gamma_2}$ are not length-commensurable.}

\vskip3mm

We would like to emphasize that while the results below for rank one
locally symmetric spaces (which include hyperbolic spaces of all
types) are unconditional, the results for spaces of
higher rank depend on the validity of the well-known conjecture in
transcendental number theory due to Schanuel (see \S \ref{S:G} for
the statement); needless to say that our results in \S\S\:2-7, 9 on
weak commensurability (in particular, Theorems A-F) {\it do not}
involve {\it any} transcendental number theory.

\vskip1mm

Henceforth,  we will assume that $G$ is a connected  absolutely
simple real algebraic group; the corresponding real Lie group $\cG =
G(\R)$ will then also be called ``absolutely simple''.
Let $\Gamma_1$ and $\Gamma_2$ be torsion-free discrete subgroups of $\cG$ such that the associated locally symmetric spaces $\fX_{\Gamma_1}$ and $\fX_{\Gamma_2}$ are
length-commensurable. Then Theorem \ref{T:G-3} implies that $\Gamma_1$ and $\Gamma_2$ are weakly commensurable. Now observing that $\fX_{\Gamma_1}$ and $\fX_{\Gamma_2}$ are
commensurable as manifolds if and only if $\Gamma_1$ and $\Gamma_2$
are commensurable up to an $\R$-automorphism of $G,$  from Theorems C and D, we obtain

\vskip3mm

\noindent {\bf Theorem \ref{T:LS1}}. {\it Each class of
length-commensurable arithmetically defined locally symmetric spaces
of $\cG = G(\R)$ is a union of finitely many commensurability
classes. It in fact consists of a single commensurabilty class if
$G$ is not of type $A_n$ $(n > 1),$ $D_{2n+1}$ ($n\geqslant 1$), $D_4$ or $E_6.$}

\vskip3mm

Furthermore, Theorem E implies the following rather surprising
result which has so far defied attempts to prove it purely
geometrically.

\vskip3mm

\noindent {\bf Theorem \ref{T:G2010}.} {\it Let $\fX_{\Gamma_1}$ and
$\fX_{\Gamma_2}$ be two arithmetically defined locally symmetric
spaces of the same absolutely simple real Lie group $\cG$. If they
are length-commensurable, then the compactness of one of them
implies the compactness of the other.}

\vskip3mm

Theorem A shows that length-commensurability provides some
information about the fundamental groups even without any
assumptions of arithmeticity.

\vskip3mm

\noindent {\bf Theorem \ref{T:G2011}.} {\it Let $\fX_{\Gamma_1}$ and
$\fX_{\Gamma_2}$ be two  locally symmetric spaces of the same
absolutely simple real Lie group $\cG$, modulo torsion-free lattices $\Gamma_1$ and
$\Gamma_2$. Denote by $K_{\Gamma_i}$ the field generated by
the traces ${\tr}\mathrm{Ad}\: \gamma$ for $\gamma \in \Gamma_i.$ If
$\fX_{\Gamma_1}$ and $\fX_{\Gamma_2}$ are length-commensurable, then
$K_{\Gamma_1} = K_{\Gamma_2}.$}

\vskip3mm

In \S 9, we present a general cohomological construction which, in
particular, enables us to give examples of
length-commensurable, but not commensurable, arithmetically defined
locally symmetric spaces associated to an absolutely simple Lie
group of any of the following types: $A_n,$ $D_{2n+1}$ $(n > 1),$ or
$E_6,$ see 9.14 (thus, the second assertion of Theorem \ref{T:LS1}
definitely cannot be extended to these types). Towards this end, we
establish a new local-global principle for the existence of an
embedding of a given $K$-torus as a maximal torus in an absolutely
simple simply connected $K$-group (for the precise assertion, see
Theorem 9.5). Using this local-global principle, we  show that there
exist nonisomorphic $K$-forms $G_1$ and $G_2$ of an absolutely
simple $K$-group of each of the types $A_n$, $D_{2n+1}$  ($n>1$), or
$E_6$, such that (i) $G_1$ is isomorphic to $G_2$ over $K_v$, for
all places $v$ of $K$ (so $G_1(A)$ is isomorphic to $G_2(A)$ as a topological group, where $A$ is the ad\`ele ring of $K$), and (ii) given a maximal $K$-torus $T_i$ of
$G_i$, there is an isomorphism $G_i\rightarrow G_{3-i}$ whose
restriction to $T_i$ is defined over $K$. Such $K$-forms are likely
to be of interest in Langlands program. Given such nonisomorphic
$K$-forms $G_1$ and $G_2$, any arithmetic subgroup $\Gamma_1$ of
$G_1(K)$ is weakly commensurable, but not commensurable, to any
arithmetic subgroup $\Gamma_2$ of $G_2(K)$, and the associated
locally symmetric spaces $\mathfrak{X}_{\Gamma_1}$ and
$\mathfrak{X}_{\Gamma_2}$ are length-commensurable but not
commensurable (see Proposition 9.13 and 9.14). \vskip3mm

Since isospectral compact locally symmetric spaces of nonpositive
curvature are length-commensurable  (Theorem \ref{T:S1}), the
following theorem, which answers Mark Kac's famous question ``Can
one hear the shape of a drum?'' for arithmetically defined compact
locally symmetric spaces, is a consequence of Theorem \ref{T:LS1}.

\vskip3mm

\noindent{\bf Theorem \ref{T:S3}.} {\it Any two arithmetically
defined compact isospectral locally symmetric spaces of an
absolutely simple real Lie group of type other than $A_n$ $(n > 1),$
$D_{2n+1}$ $(n\geqslant 1)$, $D_4$ and $E_6$, are commensurable to each other.}

\vskip2mm

We finally mention some results dealing with arithmeticity. If
$\Gamma_1$ and $\Gamma_2$ are torsion-free lattices in $\cG$ such
that $\fX_{\Gamma_1}$ and $\fX_{\Gamma_2}$ are length-commensurable,
then Theorem \ref{T:G-3} implies that $\Gamma_1$ and $\Gamma_2$ are
weakly commensurable, and according to Theorem F, if one of them is
arithmetic, then so is the other (Theorem \ref{T:G25}).
Moreover, if $\Gamma_1$ and $\Gamma_2$ are torsion-free cocompact
discrete subgroups of $\cG$ such that $\fX_{\Gamma_1}$ and
$\fX_{\Gamma_2}$ are isospectral, then by Theorem \ref{T:S1} they
are length-commensurable, and consequently, we again see that if one
of the $\Gamma_i$ is arithmetic then so is the other (Theorem
\ref{T:S4}).

\vskip4mm

{\it Notations and conventions.} Unless stated otherwise, all our
fields will be of characteristic zero. For a number field $K,$ we
let $V^K$ (resp.,\,$V^K_{\infty}$ and $V^K_f$) denote the set of all
places (resp., the subsets of archimedean and nonarchimedean
places). For a torus $T,$ we let $X(T)$ denote the character group,
and for a morphism $\pi \colon T_1 \to T_2$ between two tori,  we let $\pi^*
\colon X(T_2) \to X(T_1)$ denote the induced homomorphism of
the character groups. If $T$ is defined over $K$, then $K_T$ will
denote the (minimal) splitting field of $T$ over $K$ and $X(T)$ will
be considered as a module over the Galois group $\Ga(K_T/K).$
\vskip1mm

In the sequel, all number fields are assumed to be contained in the
field $\C$ of complex numbers. For a subfield $K$ (resp.,\,$K_i$) of
$\C$, $\overline K$ (resp.,\,${\overline K}_i$) will denote its
algebraic closure in $\C$. For a place $v$ of a number field $K$
(resp.,\,$K_i$), ${\overline K}_v$ (resp.,\,${{\overline K}_i}_v$)
will denote an algebraic closure of the completion $K_v$
(resp.,\,${K_i}_v$) of $K$ (resp.,\,$K_i$) at $v$. In particular,
${\overline \Q}_p$ will denote an algebraic closure of $\Q_p$.
\vskip1mm

For an algebraic $K$-subgroup $G$ of $\mathrm{GL}_n$, and a subring $R$ of a 
commutative $K$-algebra $C$, $G(R)$ will denote the group $G(C)\cap \mathrm{GL}_n(R)$.

\vskip4mm

{\it Acknowledgements.} We thank Alejandro Uribe and Steve Zelditch
for the proof of Theorem \ref{T:S1}, Skip Garibaldi and Peter Sarnak for their
comments and corrections, and Peter Abramenko and Robert Griess for useful discussions on Weyl groups and root systems.

Both the authors were supported by the Institute for Advanced Study,
the NSF (grants DMS-0400640 and DMS-0502120), the BSF (grant \#
200071), the Humboldt Foundation and SFB 701 (Bielefeld). Part of
this work was done while the first-named author visited Nagoya
University. He would like to thank Shigeyuki Kondo for his kind
hospitality.

\section{Preliminaries}\label{S:P}
We begin with a simple comment on the notion of weak
commensurability of semi-simple elements.
\begin{lemma}\label{L:P0}
Let $\gamma_1 , \gamma_2 \in G(F)$ be semi-simple elements. The
following conditions are equivalent:

\vskip3mm

{\rm (1)} \parbox[t]{11.5cm}{$\gamma_1$ and $\gamma_2$ are weakly
commensurable, i.e., {\rm there exist} maximal $F$-tori $T_i$ of $G$
for $i = 1 , 2$ such that $\gamma_i \in T_i(F)$ and $\chi_1(\gamma_1) =
\chi_2(\gamma_2) \neq 1$ for some characters $\chi_i \in X(T_i);$}

\vskip2mm

{\rm (2)} \parbox[t]{11.5cm}{{\rm for any} maximal $F$-tori $T_i$ of
$G$ with $\gamma_i \in T_i(F),$ there exist characters $\chi_i \in
X(T_i)$ such that $\chi_1(\gamma_1) = \chi_2(\gamma_2) \neq 1.$}
\end{lemma}

While (2) trivially implies (1), the opposite implication follows
from the fact that if $C_i$ is the Zariski-closure of the subgroup
generated by $\gamma_i$, then for {\it any} torus $T_i$ containing
$C_i,$ the restriction map $X(T_i) \to X(C_i)$ is surjective
(cf.\,\cite{Bo}, 8.2).

\begin{cor}\label{C:P1}
For $i = 1 , 2,$ let $K_i$ be a subfield of $F,$ $G_i$ be an
$F/K_i$-form of $G,$ and $\gamma_i \in G_i(K_i) \hookrightarrow
G(F)$ be a semi-simple element. Then $\gamma_1$ and $\gamma_2$ are
weakly commensurable if and only if there exist maximal $K_i$-tori
$T_i$ of $G_i$ such that $\chi_1(\gamma_1) = \chi_2(\gamma_2) \neq
1$ for some $\chi_i \in X(T_i).$
\end{cor}

This follows from the lemma because every semi-simple $\gamma_i \in
G_i(K_i)$ is contained in a maximal $K_i$-torus of $G_i.$

\vskip2mm

We will now prove two elementary lemmas on weak commensurability of
subgroups. The first lemma enables one to replace each of the two
weakly commensurable subgroups with a commensurable subgroup.
\begin{lemma}\label{L:P1}
Let $\Gamma_1$ and $\Gamma_2$ be two weakly commensurable finitely
generated Zariski-dense subgroups of $G(F)$. For $i = 1,\,2$, if
$\Delta_i $ is a subgroup of $G(F)$ commensurable with $\Gamma_i$
for $i = 1 , 2$, then the subgroups $\Delta_1$ and $\Delta_2$ are
weakly commensurable.
\end{lemma}
\begin{proof}
We recall that a subgroup $\Delta$ of ${\rm GL}_n(K)$ is {\it
neat}\, if for every $\delta \in \Delta,$ the subgroup of
${\overline K}^{\times}$ generated by the eigenvalues of $\Delta$ is
torsion-free. According to a result proved by Borel (cf.\,\cite{Ra},
Theorem 6.11) every finitely generated subgroup of ${\rm GL}_n(K)$
contains a neat subgroup of finite index. We fix a neat subgroup 
$\Theta$  of   $\Gamma_1 \cap \Delta_1$ of finite index,  then
$[\Delta_1 : \Theta] < \infty.$ Given a semi-simple element
$\delta_1 \in \Delta_1$ of infinite order, we can pick $n_1 \geqslant
1$ so that $\gamma_1 := \delta_1^{n_1} \in \Theta.$ Since $\Gamma_1$
and $\Gamma_2$ are weakly commensurable, there exists a semi-simple element 
$\gamma_2 \in\Gamma_2$ of infinite order so that
$$
\chi_1(\gamma_1) = \chi_2(\gamma_2) \neq 1,
$$
where for $i=1,\,2$, $\chi_i$ is a character of a maximal $F$-torus $T_i$ of $G$ such that
$\delta_1\in T_1(F)$ and $\gamma_2\in T_2(F)$. Now, pick $n_2
\geqslant 1$ so that $\delta_2 := \gamma_2^{n_2} \in \Gamma_2 \cap
\Delta_2.$ Then
\begin{equation}\label{E:P1}
(n_2\chi_1)(\gamma_1) = ((n_1n_2)\chi_1)(\delta_1) =
\chi_2(\delta_2).
\end{equation}
It remains to observe that since $\chi_1(\gamma_1) \neq 1$ belongs
to the subgroup generated by the eigenvalues of $\gamma_1,$ which 
is torsion-free, it is not a root of unity. This
implies that the common value in (\ref{E:P1}) is $\neq 1,$ and
therefore $\delta_1$ and $\delta_2$ are weakly commensurable. Thus,
every semi-simple $\delta_1 \in \Delta_1$ of infinite order is
weakly commensurable to some semi-simple $\delta_2 \in \Delta_2$ of infinite order, 
and by symmetry, every semi-simple $\delta_2 \in \Delta_2$ of
infinite order is weakly commensurable to some semi-simple $\delta_1
\in \Delta_1$ of infinite order, which makes $\Delta_1$ and $\Delta_2$ weakly
commensurable.
\end{proof}

\vskip2mm

The next lemma shows that in the analysis of weak commensurability
of subgroups, one can replace the ambient algebraic group with an
isogenous group.
\begin{lemma}\label{L:P2}
Let $\pi \colon G \to G'$ be an $F$-isogeny of connected semi-simple algebraic
$F$-groups, and let $\Gamma_1$, $\Gamma_2 $ be two finitely
generated Zariski-dense subgroups of $G(F)$. Then $\Gamma_1$ and
$\Gamma_2$ are weakly commensurable if and only if their images
$\Gamma'_1 = \pi(\Gamma_1)$ and $\Gamma'_2 = \pi(\Gamma_2) \subset
G'(F)$ are weakly commensurable.
\end{lemma}
\begin{proof}
One direction is almost immediate. Namely, suppose $\Gamma'_1$ and
$\Gamma'_2$ are weakly commensurable. Then for a given semi-simple
element $\gamma_1$ of $\Gamma_1$ of infinite order, there exists a
semi-simple element $\gamma_2 \in \Gamma_2$ of infinite order 
so that for $i=1,\,2$, there exist a maximal
$F$-torus $T'_i$ of $G'$, and a character
$\chi'_i$ of  $T'_i$, such that $\pi(\gamma_i)\in T'_i(F)$, and
$$
\chi'_1(\pi(\gamma_1)) = \chi'_2(\pi(\gamma_2)) \neq 1.
$$
Then, for $i=1,\,2$,  $T_i := \pi^{-1}(T'_i)$ is  a  maximal $F$-torus of $G$,
$\gamma_i\in T_i(F)$, and for their characters $\chi_i =
\pi^*(\chi'_i)$ we have
$$
\chi_1(\gamma_1) = \chi_2(\gamma_2) \neq 1.
$$
This, combined with a ``symmetric"  argument, implies that
$\Gamma_1$ and $\Gamma_2$ are weakly commensurable.

Conversely, suppose that $\Gamma_1$ and $\Gamma_2$ are weakly
commensurable, and for $i=1$, $2$, pick neat subgroups $\Delta_i$ of
$\Gamma_i$ of finite index. By Lemma \ref{L:P1}, it is enough to
show that $\pi(\Delta_1)$ and $\pi(\Delta_2)$ are weakly
commensurable. Let $\delta_1$ be a nontrivial semi-simple element of
$\Delta_1$. Then there exists a semi-simple element $\delta_2 \in \Delta_2$ such 
that for $i=1,\,2$, there exist a
maximal $F$-torus $T_i$ of $G$, with $\delta_i\in T_i(F)$,  and a character
$\chi_i$ of  $T_i$, so that
\begin{equation}\label{E:P50}
\chi_1(\delta_1) = \chi_2(\delta_2) \neq 1.
\end{equation}
Set $T'_i = \pi(T_i)$. Then  $\pi (\delta_i)\in T'_i(F)$. If $m = \vert \ker \pi
\vert$, then there exist characters $\chi'_i \in X(T'_i)$ such that
$m\chi_i = \pi^*(\chi'_i).$ Since $\Delta_1$ is neat, the common
value in (\ref{E:P50}) is not an $m$-th root of unity, and then
$$
\chi'_1(\pi(\delta_1)) = \chi'_2(\pi(\delta_2)) = \chi_1(\delta_1)^m
\neq 1.
$$
This, together with a ``symmetric" argument, implies that
$\pi(\Delta_1)$ and $\pi(\Delta_2)$ are weakly commensurable.
\end{proof}

\vskip1mm

Next, we prove the following (known) proposition which characterizes
commensurable $S$-arithmetic subgroups. Since we have not been able
to find a reference for its proof, we give a complete argument.
\begin{prop}\label{P:P1}
Let $G$ be a connected absolutely simple algebraic group over a field $F$ of
characteristic zero, and let $\Gamma_i \subset G(F),$ for $i = 1 ,
2,$  be a Zariski-dense $(G_i , K_i , S_i)$-arithmetic subgroup.
Then $\Gamma_1$ and $\Gamma_2$ are commensurable up to an
$F$-automorphism of $G$ if and only if $K_1 = K_2 =: K,$ $S_1 =
S_2,$ and $G_1$ and $G_2$ are $K$-isomorphic.
\end{prop}
\begin{proof}
Let $\iota_i \colon G_i \to G$ be an $F$-isomorphism used to define
$(G_i , K_i , S_i)$-arithmetic subgroups, where $G_i$ is defined
over a subfield $K_i$ of $F$ (technically, $\iota_i$ should have
been written as $\iota_i \colon _{F}G_i \to G,$ but our simplified
notation will not lead to a confusion). One implication is obvious.
Namely, suppose $K_1 = K_2 =: K,$ $S_1 = S_2 =: S,$ and let $\tau
\colon G_1 \to G_2$ be a $K$-isomorphism. Then $\tau(G_1(\cO_K(S)))$
is commensurable with $G_2(\cO_K(S)),$ and $\sigma := \iota_2 \circ
\tau \circ \iota_1^{-1}$ is an $F$-automorphism of $G.$ Clearly,
$\sigma(\iota_1(G_1(\cO_K(S))))$ is commensurable with
$\iota_2(G_2(\cO_K(S))),$ implying that $\sigma(\Gamma_1)$ is
commensurable with $\Gamma_2,$ as required.

Conversely, suppose $\sigma$ is an $F$-automorphism of $G$ such that
$\sigma(\Gamma_1)$ and $\Gamma_2$ are commensurable. Then the assertion 
that $K_1 = K_2$ is an immediate consequence of the following lemma applied to
$\sigma(\Gamma_1) \cap \Gamma_2$ which is both $(G_1 , K_1 ,
S_1)$- and $(G_2 , K_2 , S_2)$-arithmetic.
\begin{lemma}\label{L:P2007}
Let $G$ be a connected absolutely simple algebraic group over a field $F$ of
characteristic zero, and $\Gamma \subset G(F)$ be a Zariski-dense
$(G_0 , K , S)$-arithmetic subgroup. Then the subfield $K_{\Gamma}$
of $F$ generated by $\tr \mathrm{Ad}_G(\gamma)$ for $\gamma \in
\Gamma$ coincides with $K.$
\end{lemma}
\noindent {\it Proof.} We will assume (as we may) that the group $G$
is adjoint and is realized as a linear group by means of the adjoint
representation on the Lie algebra $\mathfrak{g}$ of $G.$ By
definition, there exists an $F$-isomorphism $\iota \colon G_0 \simeq
G,$ and we will use its differential $d\iota$ to identify the Lie
algebra $\mathfrak{g}_0$ of $G_0$ with $\mathfrak{g}.$ Set
$$
\Gamma_0 := \iota^{-1}(\Gamma) \subset G_0(F).
$$
Then $\Gamma_0$ is a $(K , S)$-arithmetic subgroup. As $G_0$ is of adjoint type, $\Gamma_0$ is contained in $G_0(K)$ (see,
for example, Proposition 1.2 of \cite{BP}). This implies that $\tr\mathrm{Ad}_{G_0}(\Gamma_0) =
\tr \mathrm{Ad}_G(\Gamma) \subset K$, hence the inclusion $K_{\Gamma} \subset K.$ To prove the reverse inclusion, we observe that according to
Theorem 1 of Vinberg\:\cite{Vn}, there exists a basis of
$\mathfrak{g}_0$ in which $\Gamma_0$ is represented by matrices
with entries in $K_{\Gamma},$ and then $G_0$ is actually defined
over $K_{\Gamma}.$ Let $A \subset \mathrm{End}\ \mathfrak{g}$ be the
linear span of $\Gamma_0.$ Then $A$ is invariant under conjugation
by $\Gamma_0,$ hence by $G_0,$ so we can consider the corresponding
(faithful) representation $\rho \colon G_0 \to \mathrm{GL}(A).$ Let $A_0$ be the 
$K_{\Gamma}$-span of $\Gamma_0.$ Any subgroup $\Gamma'_0$ of  $\Gamma_0$ 
of finite index has the same Zariski-closure as $\Gamma_0$ (viz., $G_0$), and hence the same
$K_{\Gamma}$-span (viz., $A_0$). Since for any $g \in G_0(K)$, the
intersection $\Gamma_0 \cap g^{-1}\Gamma_0g$ is of finite index
in $\Gamma_0,$ we see that $A_0$ is invariant under conjugation
by $G_0(K),$ and therefore, in terms of a basis of $A$ contained in $A_0$, $\rho(G_0(K))$ is represented by matrices
with entries in $K_{\Gamma}.$ Thus, $G_0(K_{\Gamma}) = G_0(K),$ so
the lemma is implied by the following.

\vskip3mm

\begin{lemma}\label{L:P3}
Let $G$ be a connected reductive algebraic group of positive dimension defined
over an infinite field $K.$ Then for any nontrivial extension $F/K,$
$G(K) \neq G(F).$
\end{lemma}
\begin{proof}
We may assume that $G$ is connected. It is known that $G$ is
unirational over $K$ (cf.\,\cite{Bo}, Theorem 18.2), i.e., there
exists a dominant $K$-rational map $f \colon \A^n \to G.$ We pick a
line $\ell$ in  $\A^n$, defined over $K$, such that $f$ restricts to
a nonconstant map on $\ell$. Let $C$ be the Zariski-closure of
${f(\ell (K))}$. Then $C$ is a curve defined over $K$; furthermore,
by L\"uroth's theorem, $C$ is rational over $K,$ i.e., it is
$K$-isomorphic to an open subvariety of $\A^1.$ This immediately
implies that $C(K) \neq C(F),$ and our claim follows.
\end{proof}

\vskip1mm

Now to complete the proof of Proposition \ref{P:P1}, consider the 
$F$-isomorphism $\tau = \iota_2^{-1} \circ \sigma
\circ \iota_1$ between $G_1$ and $G_2.$ We can obviously choose
subgroups $\Delta_i$ of $G_i(\cO_{K_i}(S_i))$ of finite index so
that $\sigma(\iota_1(\Delta_1)) = \iota_2(\Delta_2),$ and then
$\tau(\Delta_1) = \Delta_2.$ Since $\Delta_i$ is a Zariski-dense
subgroup of $G_i(K),$ where $K := K_1 = K_2,$ we see that $\tau$ is
in fact defined over $K.$ Next, take any  $v \notin S_1.$ Since the
closure of $\Delta_1$ in $G_1(K_v)$ is compact, we obtain that the
closure of $\Delta_2 = \tau(\Delta_1)$ in $G_2(K_v)$ is also
compact. If we assume that $v \in S_2$, then the facts that
$G_2(K_v)$ is noncompact by our construction, and $\Delta_2$ is a
lattice in $\prod_{w \in S_2} G_2(K_w)$, yield a contradiction.
Thus, $v \notin S_2,$ proving the inclusion $S_2 \subset S_1.$ The
opposite inclusion is proved similarly, so $S_1 = S_2.$
\end{proof}

\vskip3mm

\noindent {\bf Remark 2.8.} The assertion of Lemma \ref{L:P3}
remains true also over a finite field $K$ for any connected
reductive group $G$ which is not a torus. Indeed, in this case $G$
is quasi-split over $K$ (cf.\,\cite{Bo}, Proposition 16.6), and
therefore it contains a 1-dimensional split torus $C.$ Clearly,
$C(K) \neq C(F),$ implying that $G(K) \neq G(F).$

Let now $G = T$ be a torus over $K = \mathbb{F}_q,$ and let $F =
\mathbb{F}_{q^m}$ with $m > 1.$ It follows from (\cite{Vos}, 9.1)
that
$$
\vert T(K) \vert = \prod_{i = 1}^d (q - \lambda_i) \ \ \text{and} \
\ \vert T(F) \vert = \prod_{i = 1}^d (q^m - \lambda_i^m),
$$
where $\lambda_i$ are certain complex roots of unity and $d = \dim
T.$ We have
$$
\vert q - \lambda_i\vert \leqslant q + 1 \ \ \text{and} \ \ \vert
q^m - \lambda_i^m \vert \geqslant q^m - 1,
$$
so if $q^m - q > 2,$ which is always the case unless $q = 2=m ,$
then $\vert T(F) \vert > \vert T(K) \vert.$ Suppose now that $q =2=
m.$ Clearly, $\vert T(K) \vert = \vert T(F) \vert$ is possible only
if $\vert q - \lambda_i \vert = q + 1,$ i.e., $\lambda_i = -1,$ for
all $i.$ This means that $T \simeq
\left(R_{F/K}^{(1)}(\mathrm{GL}_1) \right)^d,$ where
$R_{F/K}^{(1)}(\mathrm{GL}_1)$ is the norm one torus associated with
the extension $F/K = \mathbb{F}_4/\mathbb{F}_2.$ For these tori we
have $T(K) = T(F),$ and our argument shows that these are the only
exceptions to Lemma \ref{L:P3} over finite fields.

\section{Results on irreducible tori}\label{S:IT}

A pivotal role in the proof of Theorems A-F is played by a
reformulation of Theorem 3 of \cite{PR2}. To explain this
reformulation, we need to introduce some additional notation. 

\vskip1mm

Let $K$ be an infinite field and $G$ be a
connected absolutely simple algebraic $K$-subgroup of
${\rm GL}_n$. Let $T$ be a maximal $K$-torus of
$G.$ As usual, $\Phi = \Phi(G , T)$
will denote the root system of $G$ with respect to
$T$, and $W(\Phi)$, or $W(G,T)$, the
Weyl group of $\Phi$. We shall denote by $K_{T}$
the (minimal) splitting field of $T$ in a fixed algebraic
closure $\overline{K}$ of $K$. Then there exists a natural
injective homomorphism $\theta_{T} \colon
\Ga(K_{T}/K) \to \mathrm{Aut}(\Phi).$
The following result is a strengthening of Theorem 3(i) of
\cite{PR2}, which does not require any significant changes in the
proof.
\begin{thm}\label{T:IR1}
Let $G$ be a connected absolutely simple algebraic group defined
over a~finitely generated field $K$ of characteristic
zero, and $L$ be a finitely generated field containing
$K.$ Let $r$ be the number of nontrivial conjugacy
classes of the Weyl group of $G,$ and suppose that we are
given $r$  inequivalent nontrivial discrete valuations $v_1, \ldots ,
v_r$ of $K$ such that the completion $K_{v_i}$
is locally compact and contains $L,$ and $G$
splits over $K_{v_i},$ for  $i = 1, \ldots , r.$ There
exist
 maximal $K_{v_i}$-tori $T(v_i)$ of
$G,$ one for each $i \in \{1, \ldots , r\},$ with the
property that for any maximal $K$-torus $T$ of
$G$ which is conjugate to $T(v_i)$ by an element
of $G(K_{v_i})$ for all $i = 1, \ldots , r,$ we
have
\begin{equation}\label{E:P100}
\theta_{T}(\Ga(L_{T}/L))
\supset W(G, T),
\end{equation}
where $L_{T} =
K_{T}L$ is the splitting field of
$T$ over $L$  so that
$\Ga(L_{T}/L)$ can be identified with
a subgroup of $\Ga(K_{T}/K)$.
\end{thm}
We will now derive a series of corollaries that will be used in the
subsequent sections.
\begin{cor}\label{C:IR1}
Let $G,$ $K$ and $L$ be as in Theorem
\ref{T:IR1}, and let $V$ be a finite set of nontrivial
valuations of $K$ such that for every $v \in V,$
the completion $K_v$ is locally compact. Suppose that for
each $v \in V$ we are given a maximal
$K_v$-torus $T(v)$ of $G.$ Then there
exists a maximal $K$-torus $T$ of $G$
for which {\rm (\ref{E:P100})} holds and which is conjugate to
$T(v)$ by an element of $G(K_v),$ for
all $v \in V.$
\end{cor}
\begin{proof}
Let $r$ denote the number of nontrivial conjugacy classes in the Weyl group of
$G.$ Enlarging $L$ if necessary, we
assume that $G$ splits over $L.$ By Proposition
1 of \cite{PR2}, there exists an infinite set $\Pi$ of rational
primes such that for each $p \in \Pi$ there exists an embedding
$\iota_p \colon L \to \Q_p.$ It follows that one can pick
$r$ distinct primes $p_1, \ldots , p_r \in \Pi$ so that for the
valuations $v_i$ of $K$ obtained as pullbacks of the
$p_i$-adic valuations $v_{p_i}$ on $\Q_{p_i},$ the set $R
= \{ v_1, \ldots , v_r\}$ is disjoint from $V.$ Now, let
$T(v_i),$ for $i = 1, \ldots , r,$ be the tori as in  
Theorem \ref{T:IR1}. Since the completions
$K_v$ for $v \in R \cup V$ are locally
compact, it follows from the Implicit Function Theorem that the tori
in the $G(K_v)$-conjugacy class of
$T(v)$ correspond to points of an open subset of
$\mathscr{T}(K_v),$ where $\mathscr{T}$ is the variety of maximal
tori of $G.$ Since $\mathscr{T}$ has the weak
approximation property (cf.\,\cite{PlR}, Corollary 3 in \S 7.2), there
exists a maximal $K$-torus $T$ of $G$ which is
conjugate to $T(v)$ by an element of
$G(K_v)$ for all $v \in R \cup
V.$ It follows from our construction that this torus has
the desired properties.
\end{proof}

\vskip2mm

To reformulate the above results for individual elements instead of
tori, we need the following lemma. We will call a subset of a topological
group {\it solid} if it intersects every open subgroup of that
group.

\begin{lemma}\label{L:P5}
Let $v$ be a nontrivial valuation of $K$ with locally
compact completion $K_v,$ and let $T$ be a
maximal $K_v$-torus of $G.$ Consider the map
$$
\varphi \colon G \times T \longrightarrow
G, \ \ (g , t) \mapsto gtg^{-1}.
$$
Then $$\mathscr{U}(T , v) :=
\varphi(G(K_v) ,
T_{\mathrm{reg}}(K_v)),$$ where
$T_{\mathrm{reg}}$ is the Zariski-open subvariety of
$T$ of regular elements, is a solid open subset of
$G(K_v).$
\end{lemma}

\begin{proof}
Indeed, one easily verifies that the differential $d_{(g ,
t)}\varphi$ is surjective for any $(g , t) \in
G(K_v) \times
T_{\mathrm{reg}}(K_v),$ so the openness of
$\mathscr{U}(T , v)$ follows from the Implicit Function
Theorem. Furthermore, for any open subgroup $\Omega$ of
$G(K_v),$ the set $T(K_v)
\cap \Omega$ is Zariski-dense in $T$ (cf.\,\cite{PlR},
Lemma 3.2),  and therefore it contains an element of
$T_{\mathrm{reg}}(K_v).$ So,
$\mathscr{U}(T , v) \cap \Omega \neq \emptyset$.
\end{proof}

\vskip2mm

\begin{cor}\label{C:IR2}
Let $G,$ $K,$ $L$ and $r$ be as in
Theorem \ref{T:IR1}. Furthermore, let $v_1, \ldots , v_r$ be $r$
valuations of $K$ with the properties specified in Theorem
\ref{T:IR1}, and let
$$
\delta \colon G(K) \hookrightarrow \prod_{i =
1}^r G(K_{v_i}) =: \mathscr{G}
$$
be the diagonal embedding. Then there exists a solid open subset
$\mathscr{U} \subset \mathscr{G}$ such that any $\gamma \in
G(K)$ satisfying $\delta(\gamma) \in
\mathscr{U}$ is regular semi-simple, and for the torus $T
= Z_{G}(\gamma)^{\circ},$ condition {\rm (\ref{E:P100})}
holds.
\end{cor}
Indeed, let $T(v_i),$ where $i = 1, \ldots , r,$ be the
tori given by Theorem \ref{T:IR1}. Then it is easy to see that
the set
$$
\mathscr{U} = \prod_{i=1}^r \mathscr{U}(T(v_i) , v_i)
$$
(notations as in Lemma \ref{L:P5}), satisfies all our requirements.

\vskip2mm

In \cite{PR1}, a $K$-torus $T$ was called
$K$-{\it irreducible} if it has no proper
$K$-subtori, which  is equivalent to the condition that
the absolute Galois group $\Ga(\overline {K}/K)$
acts irreducibly on the $\Q$-vector space $X(T)
\otimes_{\Z} \Q.$ It follows that, in our previous notation, a
maximal $K$-torus $T$ of $G$ such that
$\theta_{T}(\Ga(K_{T}/K))$
$\supset W(G , T)$ is $K$-irreducible
(cf.\:\cite{Bou}, Ch.\,VI, \S 1, n$^{\circ}$ 2). If $T$ is a $K$-irreducible torus, then the cyclic group generated by any element of $T(K)$ of infinite order is Zariski-dense in $T$.

We will use the
following general fact about $K$-irreducible tori.
\begin{lemma}\label{L:P6}
Let $T$ be a $K$-irreducible torus, and
$K_{T}$ be its splitting field over
$K$. Let $t \in T(K)$ be an element of
infinite order, and $\chi \in X(T)$ be a nontrivial
character. Then for $\lambda := \chi(t),$ the Galois conjugates
$\sigma(\lambda),$ with $\sigma \in \Ga(K_{T}
/K),$ generate $K_{T}$ over
$K.$
\end{lemma}
\begin{proof}
We need to show that if $\tau \in \Ga(K_{T}
/K)$  is such that
$$
\tau(\sigma(\lambda)) = \sigma(\lambda) \ \  \text{for all}  \ \
\sigma \in \Ga(K_{T} /K),
$$
then $\tau = \mathrm{id}.$ For such a $\tau$ we have
$$
(\sigma^{-1}\tau\sigma)(\chi(t)) =
((\sigma^{-1}\tau\sigma)(\chi))(t) = \chi(t).
$$
Hence, the character $(\sigma^{-1}\tau\sigma)(\chi) - \chi$ takes the
value $1$ at $t,$ and therefore, $\tau(\sigma(\chi)) = \sigma(\chi)$
because $t$ generates a Zariski-dense subgroup of $T$. But the fact that
$T$ is $K$-irreducible implies that the
characters $\sigma(\chi)$, for $\sigma \in
\Ga(K_{T} /K)$, span $X(T)
\otimes_{\Z} \Q,$ so $\tau = \mathrm{id}.$
\end{proof}

\section{The isogeny theorem}\label{S:Is}

In this
section, $K$ will be a field of arbitrary characteristic and
$K^{\mathfrak s}$ a fixed separable closure of $K$.  Let $G$ be a connected
absolutely simple algebraic $K$-group. Let $T$ and $T'$ be two
maximal tori of $G$, and $L$ any field extension of $K$ such that
both the tori are defined and split over it. Given systems $\Delta
\subset \Phi(G,T)$ and $\Delta' \subset \Phi(G , T')$ of simple
roots, there exists $g \in G(L)$ such that the corresponding inner
automorphism $i_g$ of $G$ maps $T$ onto $T'$, and the induced
homomorphism $i_g^* \colon X(T') \to X(T)$ of the character groups
maps $\Delta'$ onto $\Delta$. Such a $g$ is determined uniquely up
to an element of $T(L),$ which implies that the identification
$\Delta \simeq \Delta'$ induced by $i_g^*$ is independent of the choice of $g$. We will always
employ this identification of $\Delta$ with $\Delta'$ in
the sequel.

\vskip1mm

Now let $T$ be a maximal $K$-torus of $G.$ Fix a system
$\Delta \subset \Phi(G , T)$ of simple roots. Then for any $\sigma \in
\Ga({K^{\mathfrak s}}/K),$ there exists a unique $w_{\sigma} \in
W(G, T)$ such that $w_{\sigma}(\sigma(\Delta)) = \Delta.$ The
correspondence $\alpha \mapsto w_{\sigma}(\sigma(\alpha))$ defines
an action of $\Ga(K^{\mathfrak{s}}/K)$ on $\Delta,$ which is called
the $*$-{\it action} (cf.\,\cite{Ti0}) .
\vskip1mm

The following lemma describes some properties of the $*$-action, and
of the aforementioned identification of $\Delta$ with $\Delta'$, which will be used later in
the paper.

\begin{lemma}\label{L:P700}
{\rm({\it a})} Let $T$ and $T'$ be two maximal $K$-tori of $G,$ and
let $\Delta \subset \Phi(G , T)$ and $\Delta' \subset \Phi(G , T')$
be two systems of simple roots. Pick $g \in G(K^{\mathfrak s})$ so
that $i_g(T) = T'$ and $i_g^*(\Delta') = \Delta.$ Then $i_g^*$
commutes with the $*$-action of $\Ga(K^{\mathfrak s}/K)$ on
$\Delta'$ and $\Delta$ respectively. In particular, it carries the
orbits of the $*$-action on $\Delta'$ to the orbits of the
$*$-action on $\Delta.$

\vskip3mm

\noindent {\rm({\it b})}  The following conditions are equivalent:

\vskip2mm \ \ {\rm (i)} $G$ is an inner form (i.e., an
inner twist of a split group) over $K;$

\vskip1mm

\ {\rm (ii)} \parbox[t]{11.5cm}{$*$-action is trivial for some
(equivalently, any) maximal $K$-torus $T$ and a
system of simple roots $\Delta \subset \Phi(G ,T);$}

\vskip1mm

{\rm (iii)}
\parbox[t]{11.5cm}{$\theta_{T}(\Ga(K_{T}/K))
\subset W(G , T)$ for some (equivalently, any)
maximal $K$-torus $T$ of $G.$}

\vskip3mm

\noindent {\rm({\it c})} \parbox[t]{11.8cm}{The minimal Galois extension
$L$ of $K$ over which $G$ becomes an
inner form admits the following (equivalent) characterizations:}

\vskip2mm

\ {\rm (i)} \parbox[t]{11.5cm}{$L = (K^{\mathfrak s})^{\mathcal
{H}}$, where $\mathcal{H}$ is the kernel of the $*$-action;}

\vskip1mm

{\rm (ii)} \parbox[t]{11.5cm}{$L =
\left(K_{T}\right)^{\mathcal{H}_{T}}$, where $\mathcal{H}_T=
\theta_{T}^{-1}(\theta_{T}(\Ga(K_{T}/K)) \cap W(G , T)).$}
\end{lemma}
\begin{proof}
({\it a}): Let $\sigma \in \Ga(K^{\mathfrak s}/K) $, and
pick $w_{\sigma} \in W(G , T)$ and $w'_{\sigma}
\in W(G , T')$ so that
$w_{\sigma}(\sigma(\Delta)) = \Delta$ and
$w'_{\sigma}(\sigma(\Delta')) = \Delta'.$ We need to show that
\begin{equation}\label{E:P701}
i_g^*(w'_{\sigma}(\sigma(\alpha'))) =
w_{\sigma}(\sigma(i_g^*(\alpha'))) \ \ \text{for all} \ \ \alpha'
\in \Delta'.
\end{equation}
Since both $T$ and $T'$ are
defined over $K$, we have $g^{-1}\sigma(g) \in
N_{G}(T),$ and we let $u_{\sigma}$ denote the
corresponding element of $W(G , T).$ Then
$$
\sigma(i_g^*(\alpha')) = u_{\sigma}(i_g^*(\sigma(\alpha'))).
$$
Now, we observe that both $i_g^* \circ w'_{\sigma}
\circ \sigma$ and $w_{\sigma} \circ \sigma \circ i_g^* = w_{\sigma}
\circ u_{\sigma} \circ i_g^* \circ \sigma$ take $\Delta'$ to
$\Delta.$ This means that $$\widetilde{w} := (i_g^*)^{-1} \circ
u_{\sigma}^{-1} \circ w_{\sigma}^{-1} \circ i_g^* \circ w'_{\sigma}$$
leaves the system of simple roots $\sigma(\Delta')$ invariant. On
the other hand, $\widetilde{w} \in W(G , T').$ So,
$\widetilde{w} = 1,$ and (\ref{E:P701}) follows.

\vskip2mm

({\it b}): It follows from ({\it a}) that if the $*$-action is trivial on some
$\Delta \subset \Phi(G , T)$ for some maximal
$K$-torus $T$, then it is trivial on any $\Delta'
\subset \Phi(G , T')$ for  any maximal
$K$-torus $T'.$ On the other hand, it follows
from the description of the $*$-action on
$\Delta \subset \Phi(G , T)$ that its triviality is equivalent to
the following:
\begin{equation}\label{E:P705}
\theta_{T}(\Ga(K_{T}/K))
\subset W(G , T).
\end{equation}
This shows that (ii) and (iii) are equivalent. It remains to show
that (i) is equivalent to the inclusion (\ref{E:P705}). For this, we
assume, as we clearly may, that $G$ is adjoint. Let
$G_0$ be the $K$-split adjoint group of the same
type as $G$, and $T_0$ be a $K$-split
maximal torus of $G_0.$ Pick an isomorphism $\varphi
\colon G_0 \to G$ such that
$\varphi(T_0) = T$. Then
$$\alpha_{\sigma} = \varphi^{-1} \circ \sigma(\varphi) \ \
\text{for} \ \ \sigma \in \Ga(K^{\mathfrak s}/K),$$
defines a 1-cocycle $\alpha \in Z^1(K , \mathrm{Aut}\:
G_0)$ associated to $G.$ For any $\chi \in
X(T)$ we have  $\chi \circ \varphi \in X(T_0),$
and therefore, $\sigma(\chi \circ \varphi) = \chi \circ \varphi$ as
$T_0$ is $K$-split. An easy computation then
shows that
\begin{equation}\label{E:P201}
\sigma(\chi) = \chi \circ (\varphi \circ \alpha_{\sigma}^{-1} \circ
\varphi^{-1}).
\end{equation}
Next, (i) amounts to the assertion that $\alpha$ is cohomologous  to a $\mathrm{Int}\:G_0(K^{\mathfrak{s}})$-valued Galois cocycle $\beta: \sigma\mapsto\beta_{\sigma}$,  $\sigma \in \Ga(K^{\mathfrak{s}}/K)$, i.e., there exists $\gamma \in \mathrm{Aut}\:
G_0$ such that $\alpha_{\sigma} = \gamma^{-1} \circ
\beta_{\sigma} \circ \sigma(\gamma),$ for all $\sigma \in
\Ga(K^{\mathfrak s}/K).$ Let us show that then in
fact
\begin{equation}\label{E:P706}
\alpha_{\sigma} \in \mathrm{Int}\: G_0 \ \ \text{for all}
\ \ \sigma \in \Ga(K^{\mathfrak s}/K).
\end{equation}
Indeed, it is well-known that
$$
\mathrm{Aut}\: G_0 = \mathrm{Int}\: G_0 \rtimes
\Psi(T_0 , B_0),
$$
where $\Psi(T_0 , B_0)$ is a subgroup of the
group of all $K$-rational automorphisms of $G_0$
that leave invariant $T_0$ and a Borel
$K$-subgroup $B_0$ containing $T_0.$
Since all the elements of $\Psi(T_0 , B_0)$ are
$K$-rational, by writing $\gamma$ in the form $\gamma =
\delta \circ \psi$ with $\delta \in \mathrm{Int}\: G_0$
and $\psi \in  \Psi(T_0 , B_0),$ we obtain that
$$
\alpha_{\sigma} = \psi^{-1} \circ (\delta^{-1} \circ \beta_{\sigma}
\circ \sigma(\delta)) \circ \psi.
$$
So, since $\mathrm{Int}\: G_0 \vartriangleleft
\mathrm{Aut}\: G_0,$ we obtain (\ref{E:P706}). In addition,
since both $T_0$ and $T$ are defined over
$K$, we have $\alpha_{\sigma}(T_0) = T_0,$ and therefore for $\varkappa_{\sigma} := \varphi
\circ \alpha_{\sigma}^{-1} \circ \varphi^{-1} \,(\in \mathrm{Aut}\: G)$, $\varkappa_{\sigma}(T) = T.$
Thus, if
$G$ is an inner form, then $\varkappa_{\sigma}$ is an inner
automorphism of $G$ which leaves  $T$ invariant. Then
its restriction $\varkappa_{\sigma} \vert T$ is given by
an element of the Weyl group $W(G , T),$ so
(\ref{E:P201}) yields the inclusion (\ref{E:P705}). Conversely,
(\ref{E:P705}) in conjunction with (\ref{E:P201}) implies that
$\varkappa_{\sigma} \vert T$ is induced by an element of
$W(G , T).$ But then $\varkappa_{\sigma}$ itself
is inner, which implies that $G$ is an inner form.

\vskip2mm

({\it c}): Characterization (i) immediately follows from part ({\it b}). For
(ii), let $F =
(K_{T})^{\mathcal{H}_{T}}.$ Since
$G$ is an inner form over $L$ and splits over
$K_{T},$ by ({\it b}), for $L_{T}
= K_{T}$ we have
$$
\theta_{T}(\Ga(L_{T}/L))
\subset W(G , T),
$$
implying that $F \subset L.$ On the other hand,
using the definition of $F$ we see that
$$
\theta_{T}(\Ga(F_{T}/F))
\subset W(G , T).
$$
Then, again by ({\it b}),  $G$ is an inner form over
$F,$ and therefore $L\subset F.$
Thus, $L = F,$ as claimed.
\end{proof}
\vskip2mm

\begin{thm}\label{T:Is1} {\rm (Isogeny theorem.)}
Let $G$ be an absolutely simple algebraic group over an infinite
field $F$. For $i=1,\,2$, let $G_i$ be a form of $G$ over an
infinite subfield $K$ of $F,$ and let $L_i$ be the minimal Galois
extension of $K$ over which $G_i$ is an inner form of a split group.
Suppose that for $i = 1 , 2$, we are given a semi-simple element
$\gamma_i \in G_i(K)$ contained in a maximal $K$-torus $T_i$ of
$G_i.$ Assume that $\gamma_1$ has infinite order and that
$\theta_{T_1}(\Ga(K_{T_1}/K)) \supset W(G_1 , T_1).$ If
$\gamma_1$ and $\gamma_2$ are weakly commensurable, then there
exists a $K$-isogeny $\pi \colon T_2 \to T_1$ which carries
$\gamma_2^{m_2}$ to $\gamma_1^{m_1}$ for some integers $m_1 , m_2
\geqslant 1.$ Moreover, if $L_1 = L_2$\footnote{cf.\,Theorem
\ref{T:BC2}(2)}, then $\pi^* \colon X(T_1) \otimes_{\Z} \Q \to
X(T_2) \otimes_{\Z} \Q$ has the property $\pi^*(\Q \cdot \Phi(G_1 ,
T_1)) = \Q \cdot \Phi(G_2 , T_2),$ and in fact, if $G$ is of type
different from $B_2 = C_2,$ $F_4$ or $G_2,$ a suitable rational
multiple of $\pi^*$ maps $\Phi(G_1 , T_1)$ onto $\Phi(G_2 , T_2).$
\end{thm}
\begin{proof}
By Lemma~\ref{L:P0},  there exist characters $\chi_i \in X(T_i)$
such that
$$
\chi_1(\gamma_1) = \chi_2(\gamma_2) =: \lambda \neq 1.
$$
We will proceed by showing first that $T_2$ is $K$-irreducible and that
the splitting fields $K_{T_1}$ and $K_{T_2}$ coincide. The first
assertion requires the following.
\begin{lemma}\label{L:Is1}
Let $\Phi$ be an irreducible root system, and let $H$ be a
subgroup of $\mathrm{Aut}(\Phi)$ which contains the Weyl group
$W(\Phi).$ Then any subgroup $H'$ of  $\mathrm{Aut}(\Phi)$ which
admits a surjective homomorphism  $H'
\stackrel{\delta}{\twoheadrightarrow} H,$ acts irreducibly on the
$\Q$-vector space $\Q[\Phi]$ spanned by $\Phi.$
\end{lemma}
\begin{proof}
If $\mathrm{Aut}(\Phi) = W(\Phi),$ then our assumption implies that
$H'= W(\Phi),$ and there is nothing to prove. Next, we consider
the cases $\Phi = A_n$ $(n > 1)$ or $E_6$ where $\mathrm{Aut}(\Phi)
= W(\Phi) \times S,$ with $S = \{ \pm I \}$ (cf.\,\cite{Bou}, Tables
I and V). It is enough to show that $H'S = \mathrm{Aut}(\Phi)$ as
then any $H'$-invariant subspace would be
$\mathrm{Aut}(\Phi)$-invariant. But $H'S \neq \mathrm{Aut}(\Phi)$
can occur only when $\delta$ is an isomorphism of $H'$ onto $H =
W(\Phi)$ and $S \subset H'.$ Then $\delta(S) \subset W(\Phi)$ would
be a central subgroup of $\mathrm{Aut}(\Phi)$ of order two, hence
$\delta(S) = S$ by Schur's Lemma, which is impossible. It remains to
consider the case $\Phi = D_n$ $(n \geqslant 4).$ First, suppose
that $n \geqslant 5.$ Then $\mathrm{Aut}(\Phi)/W(\Phi)$ has order
two and $W(\Phi) = D \rtimes S_n$, where in terms of a suitable
basis $e_1,\, \ldots\, , e_n$ of $V = \Q[\Phi],$ the group $D$
consists of $\mathrm{diag}(\varepsilon_1,\, \ldots\, ,
\varepsilon_n)$ with $\varepsilon_i = \pm 1$ and $\varepsilon_1
\cdots \varepsilon_n = 1,$ and $S_n$ permutes the basic vectors
(cf.\,\,\cite{Bou}, Table IV). Then $H' \cap D$ has order $2^{n-1}$
or $2^{n-2},$ and therefore has at least $n - 2$ distinct weight
subspaces. At least one of those subspaces is 1-dimensional,
hence is spanned by a basic vector $e_i.$ So, if the action of $H'$
is not irreducible, there is a proper invariant subspace $W \subset
V$ containing $e_i.$ But $H' \cap S_n$ has index $\leqslant 2$ in
$S_n,$ hence it contains $A_n.$ Since $A_n$ acts on the basic
vectors transitively, we obtain $W = V$ -- a contradiction. It
remains to consider the case $\Phi = D_4.$ In this case also the above
description of $W(\Phi)$ remains valid, but
$\mathrm{Aut}(\Phi)/W(\Phi) \simeq S_3.$ It is well-known that a
Sylow 2-subgroup $P$ of $S_4$ acts on $\{1,2,3, 4 \}$ transitively,
which easily implies that the Sylow 2-subgroup $W(\Phi)_2 = D
\rtimes P$ acts on $V$ irreducibly. Pick a Sylow 2-subgroup $A
\subset H'$ and let $B$ be a Sylow 2-subgroup of
$\mathrm{Aut}(\Phi)$ that contains $A.$ Clearly, $[B : A] \leqslant
2,$ so it follows from the irreducibility of $B$ and Clifford's
Lemma that if $V$ is not $H'$-irreducible, then $V = V_1 \oplus
V_2$, where the $V_i$'s  are 2-dimensional $H'$-invariant
subspaces. But then the image of $H'$ in each of ${\mathrm{GL}}(V_i)$s would
be conjugate to a subgroup of ${\mathrm{O}}_2(\R),$ hence it is cyclic or
dihedral, implying that $H'$ has derived length $\leqslant 2,$
which is not the case. Thus, the action of $H'$ is irreducible.
\end{proof}

Clearly, $T_1$ is $K$-irreducible, so according to Lemma \ref{L:P6},
the conjugates $\sigma(\lambda)$ with $\sigma \in \Ga(K^{\mathfrak
s}/K),$ generate the splitting field $K_{T_1}.$ At the same time,
since $K_{T_2}/K$ is a Galois extension, all these conjugates belong
to $K_{T_2},$ yielding the inclusion $K_{T_1} \subset K_{T_2}$ and
hence a surjective homomorphism
$$
\Ga(K_{T_2}/K) \longrightarrow \Ga(K_{T_1}/K).
$$
It now follows from Lemma \ref{L:Is1} that
$\theta_{T_2}(\Ga(K_{T_2}/K))$ acts irreducibly on $\Q[\Phi(G_2 ,
T_2)] \simeq X(T_2) \otimes_{\Z} \Q,$ implying that $T_2$ is
$K$-irreducible. Now, $\gamma_2$ has infinite order as $\gamma_2^m =
1$ would imply $(m\chi_1)(\gamma_1) = 1$, which is impossible since
$\gamma_1$ generates a Zariski-dense subgroup of $T_1.$ It follows
that $\gamma_2$ generates a Zariski-dense subgroup of $T_2,$ and
therefore the conjugates $\sigma(\lambda),$ where $\sigma \in
\Ga(K^{\mathfrak s}/K),$ generate $K_{T_2}$ as well, yielding
$K_{T_1} = K_{T_2} =: \cK.$

Let $\mathscr{G} = \Ga({\cK}/K).$ We next show that there is an isomorphism of
$\Q[\mathscr{G}]$-modules
$$
\rho \colon \Q \otimes_{\Z} X(T_1) \longrightarrow \Q
\otimes_{\Z} X(T_2)
$$
that takes $\chi_1$ to $\chi_2.$ For this, we consider
$$\nu_i \colon \Q[\mathscr{G}] \to \Q \otimes_{\Z} X(T_i),  \ \ \sum a_{\sigma} \sigma
\mapsto \sum a_{\sigma} \sigma(\chi_i);$$ clearly, $\nu_i(\Z[\mathscr{G}])
\subset X(T_i).$  The irreducibility of $T_i$ implies that $\nu_i$
is surjective for $i = 1 , 2,$ so it is enough to show that
\begin{equation}\label{E:Is5}
\mathrm{Ker}\: \nu_1 = \mathrm{Ker}\: \nu_2.
\end{equation}
For this we observe that given $a = \sum a_{\sigma} \sigma \in
\Z[\mathscr{G}],$ we have
$$
\nu_1(a)(\gamma_1) = \prod \sigma(\chi_1)(\gamma_1) ^{a_{\sigma}} =
\prod \sigma(\lambda)^{a_\sigma} = \prod \sigma(\chi_2)(\gamma_2)^{a_{\sigma}} = \nu_2(a)(\gamma_2).
$$
Since for $i = 1 , 2$, $\gamma_i$ generates a Zariski-dense subgroup
of $T_i,$ the above computation shows that $\nu_1(a) = 0$ is
equivalent to $\nu_2(a) = 0,$ and (\ref{E:Is5}) follows.

The subgroup $\Theta := \nu_1(\Z[\mathscr{G}])$ has finite index, say $d,$
in $X(T_1).$ Then the multiplication by $d$ followed by $\rho$
defines the homomorphism
$$\pi^* \colon X(T_1) \to \nu_2(\Z[\mathscr{G}]) \subset X(T_2)$$
of $\mathscr{G}$-modules such that $\pi^*(\chi_1) = d\chi_2.$ Let $\pi \colon T_2 \to T_1$ be
the $K$-isogeny corresponding to $\pi^*$. Then $\chi_1(\pi(\gamma))
= \chi_2(\gamma)^d$ for every $\gamma \in T_2,$ and in particular,
$$
\chi_1(\pi(\gamma_2)) = \chi_2(\gamma_2)^d = \chi_1(\gamma_1^d).
$$
Applying the elements of $\mathscr{G},$ we see that $\chi(\pi(\gamma_2)) =
\chi(\gamma_1^d)$ for all $\chi \in \Theta,$ and therefore,
$$
\chi(\pi(\gamma_2)^d) = \chi(\gamma_1^{d^2}) \ \ \text{for every} \ \
\chi \in X(T_1).
$$
Thus, $\pi(\gamma_2)^d = \gamma_1^{d^2},$ so the first assertion of
the theorem holds with $m_1 = d^2,$ $m_2 = d.$ The second
assertion of Theorem \ref{T:Is1} will be deduced from:
\begin{lemma}\label{L:Is500}
For $i=1,\,2$, let $\Phi_i$ be an irreducible reduced root system
contained in, and spanning, the $\Q$-vector space $V_i$. We assume
that $\Phi_1$ is isomorphic to $\Phi_2$, and there is an isomorphism
$\mu \colon W(\Phi_1) \to W(\Phi_2)$ of the corresponding Weyl
groups, and a linear isomorphism $\lambda\colon V_1\to V_2$
compatible with $\mu$ (i.e., $\lambda(w(v)) = \mu(w)(\lambda(v))$
for all $v \in V_1$ and $w \in W(\Phi_1)$). Then $\lambda(\Q \cdot
\Phi_1) = \Q \cdot \Phi_2,$ and in fact, if $\Phi_1$ and $\Phi_2$
are not of type $B_2 = C_2,$ $F_4$ or $G_2,$ a suitable rational
multiple of $\lambda$ maps $\Phi_1$ onto $\Phi_2.$
\end{lemma}
\begin{proof}
We equip $V_i$ with a positive definite $W(\Phi_i)$-invariant inner product
scaled so that the short (long) roots in
$\Phi_1$ and $\Phi_2$ have the same length in the respective spaces.
We note that as $V_i$ is an absolutely irreducible $W(\Phi_i)$-module, any two $W(\Phi_i)$-invariant inner products
on $V_i$ are multiples of each other, see \cite{Bou}, Ch.\,VI, \S 1,
Proposition 7. This implies, in particular, that $\lambda$ is a multiple of an isometry.  For a root $\alpha \in \Phi_1,$ let $w_{\alpha} \in
W(\Phi_1)$ be the corresponding reflection. Then $\mu(w_{\alpha})$
is the reflection of $V_2$ with respect to $\lambda(\alpha).$ On the
other hand, $\mu(w_{\alpha})\in W(\Phi_2),$ so it follows from
(\cite{Bou}, Ch.\,V, \S 3, Cor.\,\,in n$^{\circ}$ 2) that
$\mu(w_{\alpha}) = w_{\bar{\alpha}}$  for some $\bar{\alpha} \in
\Phi_2.$ So,  $\lambda(\alpha) = t_{\alpha} \bar{\alpha}$ for some
$t_{\alpha} \in \Q,$ and our first assertion follows.

Now fix an arbitrary (resp., an arbitrary short) root $\alpha_0 \in
\Phi_1$ if all  roots have the same length (resp., if $\Phi_1$
contains roots of unequal lengths). Replacing $\lambda$ with
$t_{\alpha_0}^{-1}\lambda,$ we assume that $\lambda(\alpha_0) =
\bar{\alpha}_0.$ If all roots have the same length, then $W(\Phi_1)
\cdot \alpha_0 = \Phi_1$ and $W(\Phi_2) \cdot \bar{\alpha}_0 =
\Phi_2$ (\cite{Bou}, Ch.\,VI, \S 1, Proposition 11), yielding
$\lambda(\Phi_1) = \Phi_2.$ It remains now to deal with the root
systems of types $B_n$ and $C_n$ with $n > 2.$ Then $W(\Phi_1) \cdot
\alpha_0$ is the subset $\Phi_1^{\mathrm{short}}$ of all short
roots, and $W(\Phi_2) \cdot \bar{\alpha}_0$ is either
$\Phi_2^{\mathrm{short}}$ or $\Phi_2^{\mathrm{long}}$ depending on
whether $\bar{\alpha}_0$ is short or long (cf.\:{\it loc.\,cit.}).
But for the types under consideration, $\vert
\Phi_1^{\mathrm{short}} \vert \neq \vert
\Phi_2^{\mathrm{long}}\vert,$ and therefore,
$\lambda(\Phi_1^{\mathrm{short}}) = \Phi_2^{\mathrm{short}}.$ Since
$\alpha_0$ and $\lambda(\alpha_0)$ have the same length, $\lambda$
is an isometry. So, if $\beta_0 \in \Phi_1^{\mathrm{long}}$,  then
writing $\lambda(\beta_0) = t_{\beta_0} \bar{\beta}_0$ and observing
that the squared-length of $\beta_0$ is twice the squared-length of
any short root, we conclude that $\bar{\beta}_0$ cannot be a short
root. Therefore, $\bar{\beta}_0$ is long, $t_{\beta_0}= \pm 1$, and
it follows that $\lambda(\Phi_1) = \Phi_2$.
\end{proof}

\vskip3mm

Set $L := L_1 = L_2.$ Then it follows from Lemma \ref{L:P700} that
$$
\theta_{T_1}(\Ga(L_{T_1}/L)) = W(G_1 , T_1) \ \  \text{and} \
\ \theta_{T_2}(\Ga(L_{T_2}/L)) \subset W(G_2 , T_2).
$$
Since $L_{T_1} = L_{T_2},$ we see that the composite map
$$
\mu \colon W(G_1 , T_1)
\stackrel{\theta_{T_1}^{-1}}{\longrightarrow} \Ga(L_{T_1}/L) =
\Ga(L_{T_2}/L) \stackrel{\theta_{T_2}}{\longrightarrow} W(G_2 ,T_2)
$$
is an isomorphism of the Weyl groups compatible with $\pi^* \colon
\Q \otimes_{\Z} X(T_1)  \to \Q \otimes_{\Z} X(T_2).$ Now, the second
assertion of Theorem \ref{T:Is1} follows from Lemma \ref{L:Is500}.
\end{proof}

\vskip1mm

\noindent {\bf Remark 4.5.} The second assertion of Theorem
\ref{T:Is1} has the following consequence. We assume $\pi^*$ scaled
so that $\pi^*(\Phi(G_1 , T_1)) = \Phi(G_2 , T_2)$.  Then it induces
a $K$-isomorphism $\bar{\pi} \colon \overline{T}_2 \to
\overline{T}_1$ of the corresponding tori in the adjoint groups
$\overline{G}_i,$ which still has the property
$\bar{\pi}(\bar{\gamma}_2^{\bar{m}_2}) = \bar{\gamma}_1^{\bar{m}_1}$
for some integers $\bar{m}_1 , \bar{m}_2 \geqslant 1,$ where
$\bar{\gamma}_i$ is the image of $\gamma_i$ in $\bar{T}_i.$
Furthermore, if $Y_i$ is the dual in $V_i$ (where $V_i$ is as in
Lemma \ref{L:Is500}) of the lattice $X_i$ spanned by $\Phi_i$, then
$Y_i$ is the character group of the maximal $K$-torus
${\widetilde{T}}_i$, corresponding to the maximal torus $T_i$, of
the simply connected cover $\widetilde{G}_i$  of $G_i$, and $\pi^*$
induces an isomorphism $Y_1\to Y_2$, which in turn induces a
$K$-isomorphism $\widetilde\pi \colon {\widetilde{T}}_2\to
{\widetilde{T}}_1$. Both $\widetilde{\pi}$ and $\bar{\pi}$ extend to
$K^{\mathfrak s}$-isomorphisms $\widetilde{G}_2\to \widetilde{G}_1$
and $\overline{G}_2 \to \overline{G}_1.$ Also, if $\Delta_1$ is a
system of simple roots in $\Phi(G_1, T_1)$, and $\Delta_2 =
\pi^*(\Delta_1)$, then $\pi^*$ commutes with the $*$-action of
$\Ga(K^{\mathfrak s}/K)$ on $\Delta_1$ and $\Delta_2$ respectively.

\section{Proof of Theorems A, B and F}\label{S:AE}

We begin this section with the following two
auxiliary propositions, the first of which is a variant of Proposition
1 of  \cite{PR2}.
\begin{prop}\label{P:F201}
Let $\mathscr{F}_1 \varsubsetneqq \mathscr{F}_2 \subset \mathscr{E}$ be a tower of finitely
generated fields of characteristic zero, and let $\mathscr{R} \subset \mathscr{E}$
be a finitely generated subring. Then there exists an infinite set
of rational primes $\Pi$ such that for each $p \in \Pi,$ there are
embeddings $\iota' , \iota'' \colon \mathscr{E} \to \Q_p$ with the following
properties:

\vskip2.5mm

{\rm (1)} both $\iota'(\mathscr{R})$ and $\iota''(\mathscr{R})$ are contained in
$\Z_p;$

\vskip1mm

{\rm (2)} $\iota' \vert \mathscr{F}_1 = \iota'' \vert \mathscr{F}_1,$ but $\iota'
\vert \mathscr{F}_2 \neq \iota'' \vert \mathscr{F}_2.$
\end{prop}
\begin{proof}
First, we observe that there exists a transcendence basis $t_1,
\ldots , t_n$ of $\mathscr{E}$ over $\Q$ such that for $\cK := \Q(t_1,
\ldots , t_n)$ we have $\cK\mathscr{F}_1 \neq \cK\mathscr{F}_2.$ Indeed, let $t_1,
\ldots , t_{n_1}$ be an arbitrary transcendence basis of $\mathscr{F}_1$
over $\Q,$ and $t_{n_1 + 1}, \ldots , t_{n_2}$ be a transcendence
basis of $\mathscr{F}_2$ over $\mathscr{F}_1$ such that $$\mathscr{F}_2 \neq \mathscr{F}_1(t_{n_1 +
1}, \ldots , t_{n_2}).$$ Then, for $\cK_0 := \Q(t_1, \ldots ,
t_{n_2}),$ we have
$$
\cK_0\mathscr{F}_1 = \mathscr{F}_1(t_{n_1 + 1}, \ldots , t_{n_2}) \neq \mathscr{F}_2.
$$
Now, let $t_{n_2 + 1}, \ldots , t_n$ be a transcendence basis of
$\mathscr{E}$ over $\mathscr{F}_2.$ Then, of course, $(\cK_0\mathscr{F}_1)(t_{n_2 + 1},
\ldots , t_n) \neq \mathscr{F}_2(t_{n_2 + 1}, \ldots , t_n),$ and therefore,
$$
\cK\mathscr{F}_1 = (\cK_0\mathscr{F}_1)(t_{n_2 + 1}, \ldots , t_n) \neq \mathscr{F}_2(t_{n_2
+ 1}, \ldots , t_n) = \cK\mathscr{F}_2,
$$
as required.

Obviously, $\mathscr{E}$ is a finite extension of $\cK\mathscr{F}_1$. Let
$\mathscr{M}$ denote the Galois closure of $\mathscr{E}$ over $\cK\mathscr{F}_1.$ Then there
exists $\sigma \in \Ga(\mathscr{M}/\cK\mathscr{F}_1)$ which acts nontrivially on
$\cK\mathscr{F}_2,$ and hence on $\mathscr{F}_2.$ Let $\mathscr{R}_0$ be the subring
generated by $\mathscr{R}$ and $\sigma(\mathscr{R}).$ Since $\mathscr{M}$ is a finitely
generated field and $\mathscr{R}_0$ is a finitely generated ring, by
Proposition~1 of \cite{PR2}, one can find an infinite set of
rational primes $\Pi$ such that, for  every $p \in \Pi$, there exists an
embedding $\iota_p \colon \mathscr{M} \to \Q_p$ with the property
$\iota_p(\mathscr{R}_0) \subset \Z_p.$ Then, for $p\in\Pi$,  the embeddings
$$
\iota' = \iota_p \vert \mathscr{E} \ \ \text{and} \ \ \iota'' = (\iota_p
\circ \sigma) \vert \mathscr{E}
$$
satisfy both of our requirements.
\end{proof}

\vskip3mm

\begin{prop}\label{P:F202}
Let $G$ be a connected absolutely simple adjoint algebraic group defined over
a field $\mathscr{F}$ of characteristic zero. Let $\mathscr{E}$ be an extension of
$\mathscr{F},$ $\Gamma \subset G(\mathscr{E})$ be a Zariski-dense subgroup, and
$\cK_{\Gamma}$ be the subfield generated by the traces $\tr
\mathrm{Ad}\: \gamma$, in the adjoint representation, of  $\gamma \in \Gamma.$ 
Given two embeddings
$\iota^{(1)} , \iota^{(2)} \colon \mathscr{E} \to \Q_p$ such that
$\iota^{(1)} \vert \mathscr{F} = \iota^{(2)} \vert \mathscr{F} =: \iota,$ we
consider $G$ as a $\Q_p$-group via extension of scalars $\iota : \mathscr{F}
\to \Q_p$, and let $\rho^{(1)} , \rho^{(2)} \colon G(\mathscr{E}) \to
G(\Q_p)$ denote the homomorphisms induced by $\iota^{(1)}$ and
$\iota^{(2)},$ respectively. If

\vskip2mm

{\rm (a)} $\rho^{(i)}(\Gamma)$ is relatively compact for $i = 1 ,
2;$

\vskip1mm

{\rm (b)} $\iota^{(1)} \vert \cK_{\Gamma} \neq \iota^{(2)} \vert
\cK_{\Gamma};$

\vskip2mm

\noindent then the closure of the image of the diagonal homomorphism
$$
\rho \colon \Gamma \to G(\Q_p) \times G(\Q_p), \ \   \gamma \mapsto
(\rho^{(1)}(\gamma) , \rho^{(2)}(\gamma)),
$$
in the $p$-adic topology, is open.
\end{prop}
\begin{proof}
We begin by showing that the image of $\rho$ is Zariski-dense in $G
\times G.$
\begin{lemma}\label{L:F201}
Let $G$ be a connected simple adjoint algebraic group, and let $\rho_i \colon
\Gamma \to G,$ where $i = 1 , 2,$ be two homomorphisms of a group
$\Gamma$ with Zariski-dense images. Then either
\begin{equation}\label{E:F120}
\tr\mathrm{Ad}\: \rho_1(\gamma)   = \tr\mathrm{Ad}\:
\rho_2(\gamma)   \ \ \text{for all} \ \ \gamma \in \Gamma,
\end{equation}
or the image of the homomorphism
$$
\rho \colon \Gamma \to G \times G, \ \ \ \gamma \mapsto
(\rho_1(\gamma) , \rho_2(\gamma)),
$$
is Zariski-dense in $G \times G.$
\end{lemma}
\begin{proof}
Let $H$ be the Zariski-closure of $\rho(\Gamma)$ in $G \times G,$
and assume that $H \neq G \times G.$ Since both $\rho_1$ and
$\rho_2$ have Zariski-dense images, for the corresponding
projections we have
$$
\mathrm{pr}_i(H) = G, \ \ i = 1 , 2.
$$
Set $H_i = H\: \cap\:\: \ker\: \mathrm{pr}_i.$ Then
$\mathrm{pr}_2(H_1)$ is a normal subgroup of $G,$ and therefore it
is either $G$ or is trivial.  Furthermore, if it equals $G,$ then as
$\mathrm{pr}_1(H) = G$, we easily see that $H = G \times G.$
Similarly, $\mathrm{pr}_1(H_2)$ is ether $G$ or is trivial, and in
the former case $H = G \times G.$ Thus, since $H \neq G \times G,$
we see that  $H_i$ is trivial for $i = 1 , 2.$ This means that
$\mathrm{pr}_i$ induces an isomorphism $\epsilon_i \colon H \to G$
for $i = 1 , 2.$ Then $\sigma := \epsilon_2 \circ \epsilon_1^{-1}$
is an automorphism of $G,$ and
$$
H = \{\: (g , \sigma(g)) \ \vert \ g  \in G \}.
$$
It follows that $\rho_2 = \sigma \circ \rho_1,$ which implies
(\ref{E:F120}).
\end{proof}

We now return to the proof of
Proposition \ref{P:F202}, and denote by $\mathscr{H}$ the closure of
$\rho(\Gamma)$ in $G(\Q_p) \times G(\Q_p)$ in the $p$-adic topology.
Then $\mathscr{H}$ is a $p$-adic Lie group (cf.\:\cite{Bou}, Ch.\,III, \S 8,
Th\'eor\`eme 2), and we let $\mathfrak{h}$ denote its Lie algebra. It
follows from condition (b) that (\ref{E:F120}) does not hold, and
hence by Lemma \ref{L:F201}, $\rho(\Gamma)$ is Zariski-dense in $G
\times G.$ This immediately implies (cf.\,\cite{PlR}, Proposition 3.4)
that $\mathfrak{h}$ is an ideal of $\mathfrak{g} \times
\mathfrak{g},$ where $\mathfrak{g}$ is the Lie algebra of $G(\Q_p)$
as a $p$-adic Lie group. If the projection of $\mathfrak{h}$ 
to, say, the first component, is zero, then the image of $\rho^{(1)}$ would be
discrete, hence finite (in view of condition (a)), which is
impossible. Thus, $\mathfrak{h}$ has nonzero projections to both
components, and therefore, being an ideal of $\mathfrak{g} \times
\mathfrak{g},$ must coincide with $\mathfrak{g} \times \mathfrak{g}$
since $\mathfrak{g}$ is simple. But this means that $\mathscr{H}$ is open in
$G(\Q_p) \times G(\Q_p)$.
\end{proof}

We will now prove Theorem A. Without any loss of generality, we
assume (as we may) that the group $G$ is adjoint and use its linear 
realization given by the adjoint representation. Since $\Gamma_i$ is
finitely generated, it is contained in $\mathrm{GL}_n(F_i)$ for some
finitely generated field $F_i.$ Then the field $K_i := K_{\Gamma_i}$
is a subfield of $F_i,$ and therefore it is finitely generated, for
$i = 1 , 2.$ By symmetry, it is enough to establish the inclusion
$K_1 \subset K_2.$ Assume the contrary, and set $K = K_1K_2.$ By
Theorem 1 of Vinberg\:\cite{Vn}, one can choose a basis (which we fix
for the rest of the proof) of the Lie algebra of $G$ so that the elements of
$\Gamma_2$ are represented by matrices with entries in $K_2$ with
respect to the basis. Then
$G$ is defined over $K_2$ (hence also over $K$) and $\Gamma_2
\subset G(K_2).$ Now, pick a finitely generated extension $L$ of $K$
over which $G$ splits and
has the property that $\Gamma_1 \subset G(L).$ Furthermore, pick a
finitely generated subring $R$ of $L$ such that $\Gamma_1 \subset
G(R).$ Let $r$ be the number of nontrivial conjugacy classes of
the Weyl group of $G.$ By Proposition 1 of \cite{PR2}, there exist
rational primes $p_1, \ldots , p_r$ and embeddings $\iota_j \colon L
\to \Q_{p_j}$ such that $\iota_j(R) \subset \Z_{p_j}.$ Let $\rho_j
\colon \Gamma_1 \to G(\Z_{p_j})$ be the corresponding homomorphisms.
Then according to Lemma 2 of \cite{PR2}, the closure of the image of
the homomorphism
$$
\delta \colon \Gamma_1 \to G(\Z_{p_1}) \times \cdots \times
G(\Z_{p_r}), \ \ \ \gamma \mapsto (\rho_1(\gamma), \ldots ,
\rho_r(\gamma)),
$$
is open. Furthermore, by Corollary \ref{C:IR2}, there exists a solid
open subset $U \subset G(\Z_{p_1}) \times \cdots \times G(\Z_{p_r})$
such that any $\gamma \in \Gamma_1$ $(\subset G(L))$, with
$\delta(\gamma) \in U,$ is regular semi-simple and for the $L$-torus $T =
Z_G(\gamma)^{\circ},$ we have
\begin{equation}\label{E:F501}
\theta_T(\Ga(L_T/L)) \supset W(G , T),
\end{equation}
where $L_T/L$ is the splitting field of $T.$

Next, applying Proposition \ref{P:F201} to the tower
$$
K_2 \varsubsetneqq K \subset L
$$
we find a prime $p \notin \{ p_1, \ldots , p_r \}$ such that there
exists a pair of embeddings $\iota^{(1)} , \iota^{(2)} \colon L \to
\Q_p$ that have the same restriction to $K_2,$ but different
restrictions to $K,$  hence $K_1,$ and also satisfy
$\iota^{(i)}(R) \subset \Z_p$ for $i = 1 , 2.$ Consider the
resulting homomorphisms $\rho^{(1)} , \rho^{(2)} \colon \Gamma_1 \to
G(\Z_p)$ (as in Proposition \ref{P:F202}, $G$ is considered to be a
$\Q_p$-group by extension of scalars $K_2 \to \Q_p$ in terms of the
embedding $\iota^{(1)} \vert K_2 = \iota^{(2)} \vert K_2$). Since
$\iota^{(1)}$ and $\iota^{(2)}$ have different restrictions to $K_1
= K_{\Gamma_1},$ by Proposition \ref{P:F202}, the closure of the
image of the homomorphism
$$
\Gamma_1 \to G(\Z_p) \times G(\Z_p), \ \ \gamma \mapsto
(\rho^{(1)}(\gamma) , \rho^{(2)}(\gamma)),
$$
is open in $G(\Z_p) \times G(\Z_p).$ Since $p \notin \{ p_1, \ldots
, p_r\},$ it follows that the closure of the image of
$$
\rho \colon \Gamma_1 \to G(\Z_{p_1}) \times \cdots \times
G(\Z_{p_r}) \times G(\Z_p) \times G(\Z_p),
$$
$$
\gamma \mapsto\rho(\gamma):= (\rho_1(\gamma), \ldots , \rho_r(\gamma),
\rho^{(1)}(\gamma), \rho^{(2)}(\gamma))= (\delta(\gamma),\rho^{(1)}(\gamma),\rho^{(2)}(\gamma)),
$$
is open as well. Since $L\subset \Q_p$, $G$ splits over $\Q_p$. We fix
a $\Q_p$-split maximal $\Q_p$-torus $\mathscr{T}_1$ of $G.$ According to
\cite{PlR}, Theorem 6.21 (for a different proof,
see \cite{B}, \S 2.4), $G$ contains a $\Q_p$-anisotropic
maximal $\Q_p$-torus $\mathscr{T}_2.$ For $i=1,\, 2$, let $U_i =
\mathscr{U}(\mathscr{T}_i , v_p)$ in the notation of Lemma
\ref{L:P5}, where $v_p$ is the $p$-adic valuation on $\Q_p.$
Since the sets $U,$ $U_1$ and $U_2$ are solid in the
corresponding groups, it follows from our preceding observation
about the openness of the closure of $\mathrm{Im}\, \rho$ that there exists
$\gamma_1 \in \Gamma_1$ such that
$$
\rho(\gamma_1) \in U \times U_1\times U_2.
$$
Let $T_1 = Z_G(\gamma_1)^{\circ}.$ Since $\Gamma_1$ and $\Gamma_2$
are weakly commensurable, there exist a maximal $K_2$-torus $T_2$ of
$G$, and $\gamma_2 \in \Gamma_2\cap T_2(K_2)$ such that
$$
\chi_1(\gamma_1) = \chi_2(\gamma_2) =: \lambda \neq 1
$$
for some characters $\chi_i \in X(T_i).$ Since $\gamma_2 \in
T_2(K_2),$ $\lambda$ is algebraic over $K_2.$ Furthermore, even
though $\gamma_1$ may not have entries in $K_1,$ by Vinberg's
theorem, it is conjugate to a matrix with entries in $K_1.$ It
follows that the torus $T_1$ is definable over $K_1$ and $\gamma_1
\in T_1(K_1),$ hence $\lambda$ is algebraic over $K_1$ as well. For $i
= 1,\, 2,$ let $\mathscr{K}_i$ be the field generated over $K_i$ by
the conjugates $\sigma(\lambda)$ with $\sigma \in
\Ga(\overline{K}_i/K_i),$ and let $\mathscr{L}$ be the field
generated over $L$ by the conjugates $\sigma(\lambda)$ with $\sigma
\in \Ga(\overline{L}/L).$ We claim that
\begin{equation}\label{E:AE2}
{\mathscr K}_1L = {\mathscr L} = {\mathscr K}_2L.
\end{equation}
By looking at the minimal polynomials of $\lambda$ over $K_i$ and
$L,$ we immediately see that $\mathscr{L} \subset \mathscr{K}_iL$
for $i = 1 , 2.$ For the opposite inclusion, we first observe that
as $\delta(\gamma_1)\in U$, it follows from (\ref{E:F501})
that $T_1$ is $L$-irreducible, and therefore by Lemma \ref{L:P6},
$\mathscr{L}$ coincides with $L_{T_1},$ the splitting field of $T_1$
over $L.$ Thus, again from (\ref{E:F501}),
\begin{equation}\label{E:AE3}
\vert \Ga({\mathscr L}/L) \vert \geqslant \vert W(G, T_1) \vert.
\end{equation}
On the other hand, for both $i = 1, \, 2,$ the field ${\mathscr
K}_iL$ is contained in the splitting field $L_{T_i}$  of $T_i$ over
$L$,  and since $G$ is of inner type over $L$ (as it splits over $L$), we obtain from
Lemma \ref{L:P700}(b) that $\theta_{T_i}(\Ga(L_{T_i}/L)) \subset
W(G, T_i).$ Thus,
$$
\vert \Ga({\mathscr K}_iL/L) \vert \leqslant \vert W(G, T_i)
\vert ;
$$
combining this with (\ref{E:AE3}), we obtain (\ref{E:AE2}).

To complete the argument, we let $v_1$ and $v_2$ denote the
valuations of $K_1$ obtained by pulling back the $p$-adic valuation
on $\Q_p$ under the embeddings $\iota^{(1)} \vert K_1$ and
$\iota^{(2)} \vert K_1$, of $K_1$ into $\Q_p,$ respectively. Then, of
course, the completion ${K_{1}}_{v_i}$ can be identified with
$\Q_p$ for $i = 1 , 2.$ It follows from the description of the open
sets $U_i$ that as $\mathscr{T}_1$ splits over
${K_{1}}_{v_1}$, $T_1$ also splits over ${K_{1}}_{v_1}$, and
as $\mathscr{T}_2$ is anisotropic over ${K_{1}}_{v_2}$, so is
$T_1$. Therefore, given a nontrivial character $\chi \in X(T_1),$
there exists $\sigma \in \Ga
({\overline{K}_{1}}_{v_2}/{K_{1}}_{v_2})$ such that
$\sigma(\chi) \neq \chi.$ Then, in view of the Zariski-density of
the subgroup generated by $\gamma_1,$ we have
$$
\sigma(\chi)(\gamma_1) = \sigma(\chi(\gamma_1)) \neq \chi(\gamma_1),
$$
and consequently,
\begin{equation}\label{E:AE1}
\chi(\gamma_1) \notin {K_{1}}_{v_2} \ \ \text{for {\it any}
nontrivial} \ \ \chi\in X(T_1).
\end{equation}
Now, we extend our original embeddings $\iota^{(1)} , \iota^{(2)}
\colon L \to \Q_p$ to embeddings $\tilde{\iota}^{(1)} ,
\tilde{\iota}^{(2)} \colon \mathscr{L} \to \overline{\Q}_p.$ As
$T_1$ splits over ${K_{1}}_{v_1},$
$$
\sigma(\lambda) = \sigma(\chi_1)(\gamma_1) \in {K_{1}}_{v_1} \ \
\text{for all} \ \ \sigma \in \Ga(\overline{K}_1/K_1),
$$
and therefore, $\tilde{\iota}^{(1)}({\mathscr K}_1) \subset \Q_p.$
Then $\tilde{\iota}^{(1)}({\mathscr L}) \subset \Q_p,$ which, in
view of (\ref{E:AE2}),  implies that $\tilde{\iota}^{(1)}({\mathscr
K}_2) \subset \Q_p.$ On the other hand, it follows from
(\ref{E:AE1}) that $\tilde{\iota}^{(2)}({\mathscr K}_1) \not\subset
\Q_p,$ so $\tilde{\iota}^{(2)}({\mathscr K}_2) \not\subset \Q_p.$
But $\tilde{\iota}^{(1)}$ and $\tilde{\iota}^{(2)}$ have the same
restriction to $K_2,$ and since ${\mathscr K}_2/K_2$ is a Galois
extension, the restrictions $\tilde{\iota}^{(1)} \vert {\mathscr
K}_2$ and $\tilde{\iota}^{(2)} \vert {\mathscr K}_2$ differ by an
element of $\Ga({\mathscr K}_2/K_2),$ which shows that the
assertions 
$$
\tilde{\iota}^{(1)}({\mathscr K}_2) \subset \Q_p \ \ \text{and} \ \
\tilde{\iota}^{(2)}({\mathscr K}_2) \not\subset \Q_p
$$
are incompatible. A contradiction, which shows that our assumption
that $K_1 \not\subset K_2$ is false, and therefore, $K_1 \subset
K_2.$ This proves Theorem A.

\vskip2mm

\noindent {\bf Remark 5.4.} As we will prove soon, weakly
commensurable Zariski-dense $S$-arithmetic subgroups share not only
the field of definition, but also many other important
characteristics (cf.\,Theorems B, C and E). For arbitrary finitely
generated Zariski-dense subgroups, however, we cannot say much
beyond Theorem A. One of the reasons is that at this point,
classification results for semi-simple groups over general fields
are quite scarce. Here is one intriguing basic question in this
direction: {\it let $D_1$ and $D_2$ be two quaternion algebras over
a field $K.$ Assume that $D_1$ and $D_2$ are {\rm weakly
isomorphic,} i.e., have the same maximal subfields. Are they
isomorphic?} The answer is easily seen to be positive when
$K$ is a global field. On the other hand, M.~Rost has informed us
that over large fields (like those used in the proof of the
Merkurjev-Suslin theorem), the answer can be negative. However, for
finitely generated fields (and the fields that arise in the context
of the present paper are finitely generated), the question remains
open (apparently, even for such fields as $K = \Q(x)$). Furthermore,
if the answer turns out to be negative, one would like to know if
every class of weakly isomorphic quaternion algebras splits into
finitely many isomorphism classes (for a finitely generated field
$K$). Of course, one can ask similar questions for other types of
algebraic groups (defining two $K$-forms of the same group to be
{\it weakly isomorphic} if they have the same maximal $K$-tori).

\vskip3mm

{\bf Proof of Theorem B.}  For $i = 1 , 2,$ let $\Gamma_i$ be a
Zariski-dense $(G_i , K_i , S_i)$-arithmetic subgroup of $G(F),$ and
assume that $\Gamma_1$ and $\Gamma_2$ are weakly commensurable. By
Lemma \ref{L:P2007}, the field $K_{\Gamma_i}$ generated by $\tr
\mathrm{Ad} \gamma$ for $\gamma \in \Gamma_i,$ coincides with $K_i.$
Since the $\Gamma_i$'s are finitely generated (cf.\,\cite{PlR},
Theorem 6.1), we can now use Theorem A to conclude that
$$
K_1 = K_{\Gamma_1} = K_{\Gamma_2} = K_2 =: K.
$$

In view of the obvious symmetry, to prove that $S_1 = S_2,$ it is
enough to prove the inclusion $S_1 \subset S_2.$ Suppose there
exists $v_0 \in S_1 \setminus S_2.$ Our restrictions on $S_i$ imply
that the group $G_1$ is $K_{v_0}$-isotropic, so there exists a
maximal $K_{v_0}$-torus $T(v_0)$ of $G_1$ which is
$K_{v_0}$-isotropic. Then by Corollary~\ref{C:IR1}, there exists a
maximal $K$-torus $T_1$ of $G_1$ for which
\begin{equation}\label{E:AE10}
\theta_{T_1}(\Ga(K_{T_1}/K)) \supset W(G_1 , T_1).
\end{equation}
and which is conjugate to $T(v_0)$ under an element of
$G_1(K_{v_0}),$ hence is $K_{v_0}$-isotropic.

Clearly, $T_1$ is $K$-anisotropic, so the quotient
${T_{1}}_{S_1}/T_1(\cO_K(S_1))$ is compact, where
${T_{1}}_{S_1}=\prod_{v\in S_1}T_1(K_v)$ (cf.\,\cite{PlR}, Theorem
5.7), which implies that the quotient of $T_1(K_{v_0})$ by the
closure $C$ of  $T_1(\cO_K(S_1))$ in $T_1(K_{v_0})$ is also
compact. But as $T_1$ is $K_{v_0}$-isotropic, the group
$T_1(K_{v_0})$ is noncompact, and we conclude $C$ is noncompact as
well. Since $T_1(\cO_K(S_1))$ is a finitely generated abelian group
(cf.\:\cite{PlR}, Theorem 5.12), this implies that there exists
$\gamma_1 \in T_1(\cO_K(S_1))$ such that the closure of the cyclic
group $\langle \gamma_1 \rangle$ in $T_1(K_{v_0})$ is noncompact. We
can in fact assume that $\gamma_1 \in \Gamma_1 \cap
T_1(\cO_K(S_1)).$ By our assumption, $\gamma_1$ is weakly
commensurable to a semi-simple element $\gamma_2$ of $\Gamma_2$. Let
$T_2$ be a maximal $K$-torus containing $\gamma_2.$ Then according
to Theorem \ref{T:Is1}, there exists a $K$-isogeny $\pi \colon T_2
\to T_1$ such that $\pi(\gamma_2^{m_2}) = \gamma_1^{m_1}$ for some
integers $m_1 , m_2 \geqslant 1.$ $\pi$ induces a continuous
homomorphism $\pi_{v_0} \colon T_2(K_{v_0}) \to T_1(K_{v_0}).$ But
since $v_0 \notin S_2$ and $\Gamma_2$ is $S_2$-arithmetic, the
subgroup $\langle \gamma_2 \rangle$ has compact closure in
$T_2(K_{v_0}),$ and we obtain that $\langle \gamma_1^{m_1} \rangle,$
and hence $\langle \gamma_1 \rangle,$ has compact closure in
$T_1(K_{v_0})$; a contradiction. \hfill $\Box$

\vskip3mm

{\bf Proof of Theorem F.}\footnote{Of course, if
$\mathrm{rk}_F\: G \geqslant 2$, then $\Gamma_2$ is automatically
arithmetic by Margulis' Arithmeticity Theorem (cf.\:\cite{M}, Ch.\,IX), 
so we only need to consider the case $\mathrm{rk}_F\: G = 1.$
Our argument, however, does not depend on $\mathrm{rk}_F\: G.$} We
will assume (as we may) that $G$ is adjoint and is realized as a
linear group via the adjoint representation on its Lie algebra
$\mathfrak{g}.$ Suppose that $\Gamma_1$ is $(G_1 , K,
S)$-arithmetic; then, in particular, $\Gamma_1 \subset G_1(K)$ as
$G_1$ is adjoint (see, for example, \cite{BP}, Proposition 1.2). Let $v_0$ be the
valuation of $K$ obtained as the pullback of the normalized valuation on $F$
using the embedding $K\hookrightarrow F$. Then of course $K_{v_0}
\subset F.$  Furthermore, $v_0 \in S.$ Indeed, if $v_0 \notin S$, then
$v_0$ is nonarchimedean and the group $G_1(\cO_{K}(S))$ is
relatively compact in $G_1(K_{v_0}).$ Since $\Gamma_1$ is
commensurable with $G_1(\cO_{K}(S)),$ it would then be relatively
compact in $G_1(K_{v_0})$, and so in $G(F).$ However, as $\Gamma_1$ 
is discrete, it would be finite,
which would contradict its
Zariski-density. Moreover, being commensurable with $\Gamma_1$,
$G_1(\cO_{K}(S))$ is discrete in
$G_1(K_{v_0}).$ Combining this with the fact that $G_1(\cO_{K}(S))$
is a lattice in ${G_1}_{S} := \prod_{v \in S} G_1({K_1}_{ v}),$ we obtain
that the group $G_1(K_{v})$ is compact for all $v \in S \setminus \{
v_0 \}$ (so, in particular, ${K_1}_{ v} = \R$ for all archimedean $v \in
S \setminus \{ v_0 \}$). Because of our convention regarding $S,$ we
see that there are in fact only two possibilities: (1) $S =
V^K_{\infty},$ or (2) $v_0$ is nonarchimedean, and $S =
V^K_{\infty} \cup \{ v_0 \}.$ Furthermore, as we have seen 
above, $\Gamma_1$ is relatively compact in $G_1(K_{v})$ for all $v
\notin S.$  Thus, for any $\gamma_1 \in \Gamma_1,$ the cyclic
subgroup $\langle \gamma_1 \rangle$ is relatively compact in
$G_1(K_v)$ for all $v \in V^{K} \setminus \{ v_0 \}.$

Let $K_{\Gamma_i}$ denote the field generated by
the traces of all elements $\gamma \in \Gamma_i.$ Being lattices,
$\Gamma_1$ and $\Gamma_2$ are finitely generated, and therefore
Theorem A applies. Combining the latter with Lemma \ref{L:P2007}, we
conclude that $$K_{\Gamma_1} = K=K_{\Gamma_2} .$$ By Vinberg's theorem \cite{Vi},
there exists a basis of $\mathfrak{g}$ in which $\Gamma_2$ is
represented by matrices with entries in $K$; we fix such a basis
for the rest of the proof. Then $G$ has a $K$-form $G_2$ such that
$\Gamma_2 \subset G_2(K).$ In the sequel, the groups of points of
$G_2$ over subrings of $K$ will be understood in terms of the
realization of $G_2(K)$ as a matrix group using the basis of $\mathfrak{g}$ fixed above.
We claim that $\Gamma_2$ is commensurable with $G_2(\cO_K(S)),$
which will prove our claim. For this it is enough to establish the
following two assertions:

\vskip3mm

({\it a}) $G_2(K_v)$ is compact for all $v \in V^K_{\infty}
\setminus \{ v_0 \}.$

\vskip2mm

({\it b}) $\Gamma_2$ is bounded in $G_2(K_v)$ for all $v \in V^K_f
\setminus \{ v_0 \}.$

\vskip3mm

\noindent Indeed, since $\Gamma_2$ is finitely generated, and
therefore it is contained in $G_2(\cO_v)$ for all but finitely many
$v \in V^K_f,$ we derive from ({\it b}), in either possibility for
$S,$ that $$[\Gamma_2 : \Gamma_2 \cap G_2(\cO_K(S))] < \infty;$$
in particular, $\Gamma_2 \cap G_2(\cO_{K}(S))$ is a lattice in
$G(F),$ and hence in $G_2(K_{v_0}).$ On the other hand, it follows
from ({\it a}) that, for either of the two possibilities for $S,$ the
subgroup  $G_2(\cO_K(S))$ is a lattice in $G_2(K_{v_0}),$ implying
that $[G_2(\cO_K(S)) : \Gamma_2 \cap G_2(\cO_{K}(S))]< \infty.$

\vskip2mm

Let $v \in V^K \setminus \{ v_0 \}$ be such
that at least one of the assertions ({\it a}) and ({\it b}) fails. We will then find a regular
semi-simple element $\gamma_2 \in \Gamma_2$ of infinite order such
that the closure of $\langle \gamma_2 \rangle$ in $G_2(K_v)$ is
noncompact and for the unique maximal $K$-torus $T_2$ of $G_2$
containing $\gamma_2$ we have
\begin{equation}\label{E:Ar1}
\theta_{T_2}(\Ga(K_{T_2}/K)) \supset W(G_2 , T_2).
\end{equation}
Let us assume for a moment that such a $\gamma_2$ exists. Then since 
$\Gamma_1$ and $\Gamma_2$ are weakly commensurable, there
exists  a semi-simple element $\gamma_1 \in \Gamma_1$ which is
weakly commensurable to $\gamma_2.$ Then, if $T_1$  is a maximal
$K$-torus of $G_1$ containing $\gamma_1,$ by Theorem \ref{T:Is1}
there exists a $K$-isogeny $\pi \colon T_1 \to T_2$ which carries
$\gamma_1^{m_1}$ to $\gamma_2^{m_2}$ for some integers $m_1 , m_2
\geqslant 1.$ The isogeny $\pi$ induces a continuous group
homomorphism $\overline{\langle \gamma_1^{m_1}
\rangle} \to \overline{\langle \gamma_2^{m_2} \rangle}$ 
of the closures of the
cyclic subgroups generated by $\gamma_1^{m_1}$ and $\gamma_2^{m_2}$
in $G_1(K_v)$ and $G_2(K_v)$ respectively. As we observed above,
$\overline{\langle \gamma_1^{m_1}\rangle}$ is compact, so
$\overline{\langle \gamma_2^{m_2} \rangle}$ must also be compact, a
contradiction.

To find a $\gamma_2\in \Gamma_2$ with the desired properties we will
use the results of \cite{PR2}. Let us consider first the case where 
$v \in V^K_{\infty}\setminus \{ v_0 \}$ and $G_2(K_v)$ is noncompact
(or, equivalently, $\text{rk}_{K_v}\,G_2>0$). It has been shown in
\cite{PR2} (cf.\,the proof of Theorem 2) that there exists a regular
$\R$-regular\footnote{Given a connected semi-simple algebraic group
$G$ defined over a local field $L$, an element $x \in G(L)$ is
called $L$-regular if the number of eigenvalues, counted with
multiplicity, of modulus 1 of $\mathrm{Ad}\: x$ is minimum
possible.} semi-simple element $\gamma_2 \in \Gamma_2$ for which the
corresponding torus $T_2$ satisfies (\ref{E:Ar1}). Since the fact
that $\gamma_2$ is $\R$-regular clearly implies that the closure of
$\langle \gamma_2 \rangle$ is noncompact, we see that $\gamma_2$ has
the required properties.

Assuming now that $v\in V^K_f \setminus \{v_0\}$ and $\Gamma_2$ is 
unbounded in $G_2(K_v)$, we will prove the existence of a $\gamma_2\in \Gamma_2$ 
with the desired properties. For this, we will use the
results of \cite{PR2} in conjunction with the following result of
Weisfeiler (\cite{W}, Theorem 10.5): there exists a finite subset
$\cS$ of  $V^K$ containing $V^K_{\infty}$ such that (i) the subgroup
$\widetilde{\Gamma}_2 := \Gamma_2 \cap G_2(\cO(\cS))$ is
Zariski-dense in $G_2$, (ii) for every $v \in V^K \setminus \cS,$
the closure of $\widetilde{\Gamma}_2$ in $G_2(K_v)$ is open, and
(iii) for any $v \in \cS \setminus V^K_{\infty},$ the subgroup
$\widetilde{\Gamma}_2$ is discrete in $G_2(K_v).$ Pick such a set
$\cS,$ and first consider the case where $v \in \cS \setminus
V^K_{\infty}.$ Since $\widetilde{\Gamma}_2$ is Zariski-dense, by
\cite{PR2}, there exists a regular semi-simple element $\gamma_2 \in
\widetilde{\Gamma}_2$ of infinite order such that the corresponding
torus $T_2$ satisfies (\ref{E:Ar1}). But since
$\widetilde{\Gamma}_2$ is discrete in $G_2(K_v),$ the subgroup
$\langle \gamma_2 \rangle$ is automatically unbounded. Now, let $v
\in V^K \setminus \cS,$ and suppose that $\Gamma_2$ is unbounded in
$G_2(K_v).$ Then $G_2$ is $K_v$-isotropic and the closure of
$\Gamma_2$ in $G_2(K_v)$ is unbounded and open, so it contains the
normal subgroup $G_2(K_v)^+$ of $G_2(K_v)$ generated by the
unipotent elements (cf.\,\cite{Pr}), which is known to be an open
subgroup of $G_2(K_v)$ of finite index (cf.\,\cite{PlR}, Theorem 3.3
and Proposition 3.17). Now we fix a maximal $K_v$-torus $T_2^{v}$ of $G_2$
which contains a maximal $K_v$-split torus of the latter. Consider
the solid open subset $\mathscr{U} = \mathscr{U}(T_2^v , v)$ of
$G_2(K_v)$ given in Lemma \ref{L:P5} for $G= G_2$ and $T = T_2$. Then $\Omega_2^{v} :=
\mathscr{U} \cap G_2(K_v)^+$ is a nonempty open subset of $G_2(K_v)^+.$ Hence, 
$\Gamma_2 \cap \Omega_2^{v}$ is dense in $\Omega_2^{v}.$
Pick a $y\in \Gamma_2\cap \Omega_2^{v}$. Then an
argument similar to the one used to prove Theorem 2 in \cite{PR2}
(where instead of using Lemma 3.5 of \cite{PRa}, we use Proposition
2.6 of \cite{Pr0}) shows that there exists $x\in \Gamma_2$ such
that, for a suitable large positive integer $n$, $\gamma_2 := xy^n$
is regular $K_v$-regular, and for the unique maximal $K$-torus $T_2$
of $G_2$ containing $\gamma_2$, (\ref{E:Ar1}) holds. At the same
time, since $\gamma_2$ is $K_v$-regular, the subgroup $\langle
\gamma_2 \rangle$ is unbounded in $G_2(K_v).$ Thus, $\gamma_2$ has the 
required properties. This completes the proof of Theorem F.

\vskip3mm

\noindent {\bf Remark 5.5.} Let $\Gamma$ be a torsion-free
Zariski-dense subgroup of $G(F)$. For any positive integer $m$, the
normal subgroup $\Gamma^{(m)}$ of $\Gamma$, generated by the $m$-th
powers of the elements in $\Gamma$, is weakly commensurable with
$\Gamma$. On the other hand, it is known, see \cite{Ol}, that if
$\Gamma$ is a cocompact lattice in a real semi-simple Lie group of
real rank $1$, then there
 exists an integer $m$ such that $\Gamma^{(m)}$ is of infinite
 index in $\Gamma$. This shows that the requirement that $\Gamma_2$ be a lattice
in Theorem F cannot be omitted in case $G$ is of $F$-rank $1$. The
question whether or not a (discrete) subgroup weakly commensurable
to an irreducible  lattice (which is, of course, automatically
arithmetic) in a real semi-simple Lie group of real rank $>1$, is
itself a lattice, remains open. We would like to point out, however,
that no variant of the above method for constructing
counter-examples is likely to work in the higher rank case.

\vskip1mm

More precisely, let again $\Gamma$ be a torsion-free Zariski-dense
subgroup of $G(F).$ Given a map $\varphi \colon \Gamma\rightarrow
{\mathbb N},$ we let $\Gamma_{\varphi}$ denote the subgroup of
$\Gamma$ generated by $\gamma^{\varphi (\gamma)}$ for all
$\gamma\in\Gamma.$ This subgroup is obviously weakly commensurable
to $\Gamma$ for {\it any} choice of $\varphi.$ However, in contrast
to the case of cocompact lattices in rank one groups
discussed in the previous paragraph, or even finite index subgroups
of $\mathrm{SL}_2(\Z),$ where the subgroup $\Gamma^{(m)}$ (which corresponds
to $\varphi \equiv m$) has infinite index in $\Gamma$ for a suitable
$m,$ the subgroup $\Gamma_{\varphi}$ always has finite index in
$\Gamma$ if $\Gamma$ is ``boundedly generated" (this fact was
pointed out  to us by Thomas Delzant). Several
non-cocompact arithmetic lattices in the higher rank case are known
to be boundedly generated (see \cite{ER} for the definition of, and
most recent results on, ``bounded generation"), and for them
considering subgroups of the form $\Gamma_{\varphi}$ will never lead
to a~weakly commensurable subgroup of infinite index.

\vskip3mm

\noindent {\bf 5.6.\,A question.} Given two weakly
commensurable Zariski-dense subgroups of $G(F)$ (where $F$ is a nondiscrete
locally compact field), is it true that discreteness of  one of them
implies that of the other?

\vskip0.5cm

\section{The invariance of rank and the proof of Theorems C and
D}\label{S:BC}

In view of Theorem B, weakly commensurable Zariski-dense
$S$-arithmetic subgroups necessarily have the same field of
definition $K$ and correspond to the same set of places of $K.$ So now
the focus of our study of such subgroups shifts to identifying
common characteristics of the $K$-forms $G_i$ used to construct
them.

\begin{prop}\label{P:BC0}
Let $V_0$ be a finite set of places of $K$. Let $\Gamma_i$ be a Zariski-dense $(G_i , K , S)$-arithmetic
subgroup of $G(F)$ for $i = 1,\, 2.$ Let $L_i$ be the smallest Galois extension of $K$ over which $G_i$ is inner. If $\Gamma_1$ and $\Gamma_2$ are
weakly commensurable, then there exists a maximal $K$-torus $T_1$ of $G_1$ which contains a maximal $K_{v_0}$-split torus of $G_1$ for all $v_0\in V_0$, a maximal $K$-torus $T_2$ of $G_2$, and a $K$-isogeny $\pi:\,T_2\to T_1$. Moreover, if $L_1 =L_2$, and $G$ is either simply connected or adjoint, and it is not of type $B_2=C_2$, $F_4$ or $G_2$, then we can assume that $\pi$ is an isomorphism, and $\pi^*(\Phi(G_1,T_1))=\Phi(G_2,T_2)$.
\end{prop}
\begin{proof}
Using Corollary \ref{C:IR1}, we can find a maximal $K$-torus $T_1$
of $G_1$ which contains a maximal $K_v$-split torus of $G_1$ for
every $v \in S \cup V_0$, and for which
$$\theta_{T_1}(\Ga(K_{T_1}/K)) \supset W(G_1 , T_1).$$
Then the group ${T_{1}}_{S} :=\prod_{v\in S}T_1(K_v)$ is
noncompact, and since the quotient ${T_{1}}_{S}/T_1(\cO_K(S))$ is
compact as $T_1$ is $K$-anisotropic, we infer that $T_1(\cO_K(S))$
is infinite. Therefore, $\Gamma_1\cap T_1(K)$ contains an element
$\gamma_1$ of infinite order. By our assumption, $\gamma_1$ is
weakly commensurable to some semi-simple $\gamma_2 \in \Gamma_2 \cap
G_2(K)$. Let $T_2$ be a maximal $K$-torus of $G_2$ that contains
$\gamma_2.$ According to Theorem \ref{T:Is1}, there exists a
$K$-isogeny $\pi \colon T_2 \to T_1$. The second assertion of the proposition follows
from Theorem \ref{T:Is1} and Remark 4.5. \end{proof}

\begin{thm}\label{T:BC1}
Let $\Gamma_i$ be a Zariski-dense $(G_i , K , S)$-arithmetic
subgroup of $G(F)$ for $i = 1,\, 2.$ If $\Gamma_1$ and $\Gamma_2$ are
weakly commensurable, then
$$\mathrm{rk}_{K_v}\,G_1 =
\mathrm{rk}_{K_v}\,G_2 \ \ \text{for all} \ \  v \in V^K.$$
\end{thm}
\begin{proof}
Fix $v_0 \in V^K.$ By symmetry, it is enough to show that
$$
\mathrm{rk}_{K_{v_0}}\,G_1 \leqslant \mathrm{rk}_{K_{v_0}}\,G_2.
$$
Applying the preceding proposition to  $V_0 =\{ v_0\}$, for $i=1,\,2$, we can find a maximal $K$-torus $T_i$ of $G_i$ such that $T_1$ contains a maximal $K_{v_0}$-split torus of $G_1$, and there is a $K$-isogeny  $\pi: T_2\to T_1$. From this we see that 
\vskip1mm

\centerline{$\text{rk}_{K_{v_0}}\,G_1=\text{rk}_{K_{v_0}}\,T_1=\text{rk}_{K_{v_0}}\,T_2\leqslant \text{rk}_{K_{v_0}}\,G_2.$}
\end{proof}

\vskip1mm

For a connected absolutely simple algebraic group $G_0$ defined over a number
field $K$, we let $\Sigma(G_0 , K)$ (resp., $\Sigma^q(G_0 , K)$) be
the set of places $v$ of  $K$ such that $G_0$ is split (resp., is
quasi-split but not split) over $K_v$ (of course, $\Sigma^q(G_0 ,
K)$ is empty if $G_0$ is an inner form of a split group over $K$).
\begin{thm}\label{T:BC2}
Let $\Gamma_i$ be a Zariski-dense $(G_i , K , S)$-arithmetic
subgroup of $G(F)$ for $i = 1 , 2.$ If $\Gamma_1$ and $\Gamma_2$ are
weakly commensurable, then

\vskip3mm

{\rm (1)} $\Sigma(G_1 , K) = \Sigma(G_2 , K);$

\vskip2mm

{\rm (2)} \parbox[t]{11.5cm}{if $L_i$ is the minimal Galois
extension of $K$ over which $G_i$ becomes an~inner form (of a split
group), then $L_1 = L_2;$}

\vskip2mm

{\rm (3)} $\Sigma^q(G_1 , K) = \Sigma^q(G_2 , K).$
\end{thm}
\begin{proof}
Assertion (1) immediately follows from the preceding theorem. To prove
(2), by symmetry it is enough to show that $L_1 \subset L_2.$
Assume, if possible, that $L_1$ is not contained in $L_2$. Then
$L_1L_2$ is a Galois extension of $K$ that properly contains $L_2.$
It follows from Chebotarev's Density Theorem that there are
infinitely many $v \in V^K_f$ which split completely in $L_2$ but not
in $L_1.$ Also, $G_2$ is quasi-split over $K_v$ for all but finitely
many $v \in V^K_f,$ cf.\,\cite{PlR}, Theorem 6.7. So there exists a
$v\in V^K_f$ which splits completely in $L_2$ but not in $L_1$, and $G_2$
is quasi-split over $K_v$. Then $G_2$
actually splits over $K_v,$ i.e., $v \in \Sigma(G_2, K),$ but since
$v$ does not split in $L_1,$ we have $v \notin \Sigma(G_1 , K),$
this contradicts assertion (1). Now assertion (3) follows at once from
Theorem \ref{T:BC1}.
\end{proof}

\vskip1mm

\noindent {\bf Remark 6.4.} Technically, Theorem 6.2 and Theorem
6.3, parts (1) and (3), are consequences of the assertion in Theorem
E (to be proved in the next section) that in the situation at hand,
the Tits indices of $G_1$ and $G_2$ over $K_v$ are identical, for
all $v \in V^K.$ We decided to include the above straightforward
proofs for the following two reasons: first, the assertions of
Theorems 6.2 and 6.3 are actually used in the proof of Theorem E,
and second, we would like to show the reader that all theorems {\it
except} Theorem E can be obtained without using the Tits index.

\vskip1mm

Before we proceed to the proofs of Theorems C and D, we briefly
recall the classification of absolutely simple algebraic groups of a
given type over a  field $K$ (cf.\,\cite{S}, \cite{Ti}). Any such
group is an inner twist of a $K$-quasi-split group of the
given type. So, fix a $K$-quasi-split group $G_0.$ Notice that $G_0$
is completely determined by specifying (in addition to its type)
the minimal Galois extension $L/K$ over which it splits; this
extension necessarily has degree 1 (which means that $G_0$ splits
over $K$) if the type is different from $A_n$ $(n > 1),$ $D_n$ $(n
\geqslant 4)$, or $E_6,$ can have degree 1 or 2 for the types $A_n$, $D_n$ and $E_6$,
and can also be either a cyclic extension of degree 3 or a Galois
extension with the Galois group $S_3$ for type $D_4.$ Furthermore,
the $K$-isomorphism classes of inner twists of $G_0$
correspond bijectively to the elements lying in the image of the natural map
$$H^1(K , \overline{G}_0) \longrightarrow
H^1(K,{\mathrm{Aut}}\,{\overline{G}}_0),$$ where $\overline{G}_0$ is
the adjoint group of $G_0$ identified with its group of inner
automorphisms. When $K$ is a number field, one
considers the natural ``global-to-local" map
$$
H^1(K , \overline{G}_0) \stackrel{\omega}{\longrightarrow} \bigoplus_{v
\in V^K} H^1(K_v , \overline{G}_0),
$$
and also the truncated maps
$$
H^1(K , \overline{G}_0) \stackrel{\omega_S}{\longrightarrow} \bigoplus_{v
\notin S} H^1(K_v , \overline{G}_0),
$$
for every finite subset $S$ of $V^K.$ It is known that $\omega$ is
injective (cf.\,\cite{PlR}, Theorem 6.22) and $\mathrm{Ker}\:
\omega_S$ is finite (cf.\:\cite{S}, Theorem 7 in Ch.\,III, \S 4.6).

\vskip3mm

{\bf Proof of Theorem C.} For $G$ of type $D_{2n}$, $n>2$, Theorem C 
is Theorem 9.1 of \cite{PR6}, therefore to prove the theorem we 
assume that $G$ is not of type $A$, $D$ or $E_6$. 

Let $G_0$ be the $K$-split form of $G$. For the groups of the types under 
consideration we have
$\mathrm{Aut}\: G_0 = \overline{G}_0,$ so the group $G_i$, for $i = 1
, 2$, is obtained from $G_0$ by twisting with a Galois cocycle
representing an appropriate element $c_i$ of $H^1(K ,
\overline{G}_0).$ We need to show that $c_1 = c_2.$ For this we
notice that according to Theorem \ref{T:BC1}, we have
$\mathrm{rk}_{K_v}\: G_1 = \mathrm{rk}_{K_v}\: G_2,$ for all $v \in
V^K.$ But for the types currently being considered this implies that
\begin{equation}\label{E:BC30}
G_1 \simeq G_2 \ \ \text{over} \ \  K_v.
\end{equation}
Indeed, for $v$ real, this follows from the classification of real
forms of absolutely simple Lie algebras/real algebraic groups
(cf.\:\cite{H}, Ch.\:X, \S 6, or \cite{Ti}). Now let $v$ be nonarchimedean. For the groups under consideration, the center $Z$ of the
corresponding simply connected group is a subgroup of $\mu_2$, the kernel
of the endomorphism $x\mapsto x^2$ of ${\mathrm{GL}}_1$. In
view of the bijection between $H^1(K_v , \overline{G}_0)$ and $H^2(K_v ,
Z)$ (cf.\,\cite{PlR}, Corollary to Theorem 6.20), we see that $H^1(K_v, \overline{G}_0)\simeq H^2(K_v, Z) ={\mathrm {Br}}(K_v)_2$, and hence, $\vert H^1(K_v , \overline{G}_0) \vert \leqslant 2,$ which means that there
exists at most one nonsplit form, and therefore the equality of
ranks implies the isomorphism between the forms. If we now let
$$
\omega_v \colon H^1(K , \overline{G}_0) \longrightarrow H^1(K_v ,
\overline{G}_0)
$$
denote the restriction map, then the isomorphism (\ref{E:BC30})
implies that $\omega_v(c_1) = \omega_v(c_2),$ for all $v \in
V^K.$ Then, $\omega(c_1) = \omega(c_2),$ and therefore, $c_1 =
c_2,$ as required.\hfill $\Box$

\vskip3mm

{\bf Proof of Theorem D.} By Theorem \ref{T:BC2}\,(2), the groups
$G_1$ and $G_2$ become inner forms over the same Galois extension of $K$. 
Let $G_0$ be the unique quasi-split inner twist of $G_1$ over $K$. Next, let
$$
V_i = V^K \setminus \left(\Sigma(G_i , K) \cup \Sigma^q(G_i ,K)
\right)
$$
be the set of places $v$ of $K$ where $G_i$ is not quasi-split. It
is well-known that $V_i$ is finite (cf.\:\cite{PlR}, Theorem 6.7).
Furthermore, it follows from Theorem \ref{T:BC2},\:(1) and (3), that
$V_1 = V_2 =: V.$ Thus, by fixing $G_1$ we automatically fix a
finite set of places $V$ such that any $G_2$ as in the statement of
the theorem is quasi-split outside $V.$ Now, consider $\xi_2 \in
H^1(K , \overline{G}_0)$ which twists $G_0$ into $G_2.$ Then for all
$v \notin V,$ the group $G_2$ is quasi-split over $K_v,$ hence it is
$K_v$-isomorphic to $G_0,$ which means that $\omega_v(\xi_2)$ is
trivial. (Here we use the fact that for a quasi-split group $G_0$
over any field $F,$ the map $H^1(F , \overline{G}_0) \to H^1(F ,
\mathrm{Aut}\: G_0)$ has trivial kernel, which follows from the
observation that $\mathrm{Aut}\: G_0$ is a semi-direct product over
$F,$ of $\overline{G}_0$ and a finite $F$-group of
automorphisms corresponding to the symmetries of the Dynkin
diagram.) Thus, $\xi_2 \in \mathrm{Ker}\, \omega_V,$ so the
finiteness of this kernel yields the finiteness of the number of
$K$-isomorphism classes of possible $K$-groups $G_2$ with the
properties described in the theorem.\hfill $\Box$

\vskip3mm

We conclude this section with two explicit examples demonstrating that in
groups of type $A_n,$ $n > 1,$ the collection
of weakly commensurable arithmetic subgroups may consist of more
than one commensurability classes. Later, in \S 9, the idea underlying
these examples will be developed into a new general technique for
constructing nonisomorphic $K$-groups of type $A_n$, $D_{2n+1}$ $(n
> 1)$ and $E_6$ which contain weakly commensurable arithmetic subgroups.

\vskip2mm

\noindent {\bf Example 6.5.} Take $G = \mathrm{SL}_d,$ where $d >
2,$ over $F = \R$ (so that $G$ is of type  $A_n$ with $n = d-1
> 1$), and fix a real number field $K$. Pick four arbitrary nonarchimedean
places $v_1, v_2, v_3, v_4 \in V^K_f.$ Let $D_1$ and $D_2$ be the
central division algebras of degree $d$ over $K$ whose local
invariants ($\in \Q/\Z$) are respectively
$$
n_v^{(1)} = \left\{\begin{array}{rcl} 0\ , & v \neq v_i,\  i\leqslant 4\\
1/d \ ,&  v = v_1 \ \ \text{or}\  \  v_2 \\ -1/d \ ,&  v = v_3 \ \ \text{or}
\ \ v_4
\end{array} \right. \ \ \
\text{and} \ \ \ \  n_v^{(2)} = \left\{\begin{array}{rcl} 0\ , & v \neq v_i, \ i\leqslant 4\\
1/d \ , & v = v_1\ \ \text{or} \ \ v_3 \\ -1/d\ , & v = v_2\  \ \text{or}
\ \ v_4.
\end{array} \right.
$$
Then as $d > 2,$  the algebras $D_1$ and $D_2$ are neither
isomorphic nor anti-isomorphic. So the algebraic $K$-groups $G_1 =
{\rm{SL}}_{1,D_1}$ and $G_2 = {\rm{SL}}_{1,D_2},$ which are inner
$K$-forms of $G,$ are not $K$-isomorphic. Thus, for any finite $S
\subset V^K,$ containing $V^K_{\infty},$ the corresponding $(G_i , K
, S)$-arithmetic subgroups $\Gamma_i \subset G(F)$ are not
commensurable (cf.\,Proposition \ref{P:P1}).  On the other hand, if
$D$ is a central division algebra of degree $d$ over $K$, then an
extension $L/K$ of degree $d$ is isomorphic to a maximal subfield of
$D$ if and only if for every $v \in V^K$, and any extension $w \vert
v,$ the local degree $[L_w : K_v]$ annihilates the corresponding
local invariant $n_v$ of $D$ (cf.\,\cite{Pie}, Corollary b in \S 18.4). It
follows that the maximal subfields of either $D_1$ or  $D_2$ are
characterized as those extensions $L/K$ of degree $d$ for which
$[L_{w_i} : K_{v_i}] = d$ for $i = 1,\,2,\,3,\,4.$ Thus, $D_1$ and
$D_2$ have the {\it same} maximal subfields, which easily implies
that $\Gamma_1$ and $\Gamma_2$ are weakly commensurable. Indeed, let
$\gamma_1 \in \Gamma_1$ be a semi-simple element of infinite order,
and let $T_1$ be a maximal $K$-torus of $G_1$ that contains
$\gamma_1.$ Since $D_1$ and $D_2$ have the same maximal subfields,
there exists a $K$-isomorphism $T_1 \stackrel{\varphi}{\simeq} T_2$
with a maximal $K$-torus $T_2$ of $G_2.$ Then the subgroup
$\varphi(T_1(K) \cap \Gamma_1)$ is an $S$-arithmetic subgroup of $T_2(K),$
so there exists $n > 0$ such that $\gamma_2 := \varphi(\gamma_1)^n
\in \Gamma_2.$ Let $\chi_1 \in X(T_1)$ be a character such that
$\chi_1(\gamma_1)$ is not a root of unity. Then for $\chi_2 =
(\varphi^*)^{-1}(\chi_1) \in X(T_2)$ we have
$$
(n\chi_1)(\gamma_1) = \chi_1(\gamma_1)^n = \chi_2(\gamma_2) \neq 1,
$$
which implies that $\Gamma_1$ and $\Gamma_2$ are weakly
commensurable.

\vskip1mm

This example can be refined in two ways. First, by picking a
sufficiently large number of nonarchimedean places and modifying the
above construction accordingly, one can construct an arbitrarily
large number of noncommensurable weakly commensurable $S$-arithmetic
subgroups of the group $G(F) = \mathrm{SL}_d(\R).$ Second, suppose
$d>2$ is even, and consider the real algebraic group $G =
\mathrm{SL}_{d/2 , \HH},$ where $\HH$ is the division algebra of
Hamiltonian quaternions. Assume that $K$ is a number field that
admits a real embedding $K \hookrightarrow \R =: F,$ and we let
$v_{\infty}$ denote the real place corresponding to this embedding.
In addition to the four places $v_1, v_2, v_3, v_4 \in V^K_f$ fixed
in the above example, we pick a fifth place $v_5 \in V^K_f \setminus
\{v_1, v_2, v_3, v_4\},$ and consider the central division algebras
$D_1$ and $D_2$ of degree $d$ over $K$ with the same local
invariants at $v_1, v_2, v_3, v_4$ as above, and having the
invariant $1/2$ at $v_{\infty}$ and $v_5,$ and $0$ everywhere else.
Then for any finite $S \subset V^K$ containing $V^K_{\infty}$ (in
particular, for $S = V^K_{\infty}$ itself), the corresponding $(G_i
, K , S)$-arithmetic subgroups are weakly commensurable, but not
commensurable, and in addition are contained in $G(F) =
\mathrm{SL}_{d/2}(\HH).$ Furthemore, by increasing the number of
places picked, we can construct an arbitrarily large number of
noncommensurable weakly commensurable $S$-arithmetic subgroups of
$\mathrm{SL}_{d/2}(\HH).$

The above  construction implemented for $K = \Q$
and $d = 4$ has the following geometric significance. Over $\R,$ the
group $G$ is isomorphic to the spinor group of a real quadratic form
with signature $(5 , 1),$ and therefore the associated symmetric
space is the real hyperbolic 5-space. So, the noncommensurable
arithmetic subgroups constructed above give rise to noncommensurable
length-commensurable compact hyperbolic 5-manifolds (cf.\,Remark
8.11). We will elaborate on this observation in \S 9, where, in
particular, noncommensurable length-commensurable compact hyperbolic
manifolds of dimension $4n+1$ ($n\geqslant 1$) will be constructed.

\vskip2mm

\noindent {\bf Example 6.6.} Let $K$ be a number field and $L$ be a
quadratic extension of $K$. For $i=1,\,2$, let $v_i$ be a
nonarchimedean place of $K$ which splits in $L$, and $v_i',\,v_i''$
be the places of $L$ lying over $v_i$. Let $d>1$ be an odd integer.
Let $D_1$ and $D_2$ be the division algebra over $L$ of degree $d$
whose local invariants are respectively
$$
n_v^{(1)} = \left\{\begin{array}{rcl}  1/d \ ,& v = v_1' \ \
\text{or}\  \  v_2' \\ -1/d \ , & v = v_1'' \ \ \text{or} \ \ v_2'',
\end{array} \right.  \ \ \
\text{and} \ \ \ \  n_v^{(2)} = \left\{\begin{array}{rcl} 1/d \ , & v = v_1' \ \  \text{or} \  \
v_2''  \\  -1/d  \ , &  v = v_1''   \ \  \text{or} \ \ v_2',
\end{array} \right.
$$
and whose local invariant at every other place of $L$ is zero. Then
for $i = 1 , 2,$ the algebra $D_i$ admits an involution $\sigma_i$
of the second kind such that the fixed field $L^{\sigma_i}$
coincides with $K.$ Let $G_i$ be the absolutely simple $K$-group with
$$G_i(K) = \{ x\in D_i^{\times}\,|\
x\sigma_i(x) =1,\ {\rm{Nrd}}\,x = 1 \}.$$ Then $G_i$ is an outer 
form of type $A_n$ with $n = d - 1 > 1.$ For
simplicity, let us assume that the involutions are chosen so that
$G_1$ and $G_2$ are quasi-split at every real place of $K$ which
does not split in $L$ (as $d$ is odd,  $G_1$ and $G_2$ split
at all the other real places of $K$). Furthermore, since $d$ is odd,
$G_1$ and $G_2$ are automatically quasi-split at every
nonarchimedean place of $K$ which does not split in $L.$ Thus, it
follows from Proposition A.2 of Appendix A in \cite{PR0} and the
subsequent discussion that for an extension $P/L$ of degree $d$
provided with an automorphism $\tau$ of order two which induces the
nontrivial automorphism of $L/K,$ an embedding $(P , \tau) \to (D_i
, \sigma_i)$ as algebras with involution exists if and only if $[P_w
: K_{v_j}] = d$ for $j = 1 , 2$ and $w \vert v_j.$ This easily
implies that the maximal $\sigma_1$-invariant subfields in $D_1$ are
the same as the maximal $\sigma_2$-invariant subfields in $D_2,$ and
therefore $G_1$ and $G_2$ have the same maximal $K$-tori. Then as in
the previous example, we conclude that for any $S,$ the
$S$-arithmetic subgroups of $G_1$ and $G_2$ are weakly
commensurable. On the other hand, it follows from our choice of
local invariants that $G_1$ and $G_2$ are not isomorphic even over
$L,$ so the constructed $S$-arithmetic subgroups are not
commensurable. A suitable variant of this construction (applied to
$K = \Q,$ $L = \Q(i)$) enables one to construct
length-commensurable, but not commensurable, compact complex
hyperbolic $(d-1)$-manifolds, providing thereby a negative answer to
Question (2) of the introduction for complex hyperbolic manifolds of
any even dimension. We will not give the details here as the general
construction described in \S 9 yields counter-examples in {\it all}
dimensions.


\section{Proof of Theorem E}\label{S:D}

\noindent{\bf Tits index of a semi-simple algebraic group} (cf.\,\cite{Ti0},
or \cite{Spr}, \S 15.5). Let  $G$ be a connected semi-simple
algebraic $K$-group. To describe the Tits index of $G/K$, we pick a
maximal $K$-split torus $T_s$ of $G$ and a maximal $K$-torus $T$ of
$G$ containing $T_s.$ Furthermore, we choose an ordering on the
vector space $X(T_s) \otimes_{\Z} \R,$ lift it to an ordering on
$X(T) \otimes_{\Z} \R$ (cf.\,\cite{Spr}, \S 15.5--we will call such
orderings on these vector spaces {\it coherent}), and let $\Delta
\subset \Phi(G ,T)$ denote the system of simple roots associated
with this ordering. Then the {\it Tits index} of $G/K$ is the data
consisting of $\Delta$ (or the corresponding Dynkin diagram), the
subset of {\it distinguished} roots, and the
$*$-action (see \S 4). We recall that a root $\alpha \in \Delta$ (or the
corresponding vertex in the Dynkin diagram) is {\it distinguished}
if its restriction to $T_s$ is nontrivial. If $\alpha \in \Delta$ is
distinguished, then every root in the orbit $\Omega$ of $\alpha$, under 
the $*$-action,
is distinguished; this is indicated by circling together all the
vertices corresponding to the roots in $\Omega,$ and the latter is
referred to as a {\it distinguished orbit.} We note that
$\text{rk}_{K}\,G$ equals the number of distinguished orbits, and
$G$ is quasi-split over $K$ if and only if every root in $\Delta$ is
distinguished.
\vskip1mm

For a subset $\Theta$ of  $\Delta,$ we let $P_{\Theta}$ denote the
corresponding standard parabolic subgroup which contains the
centralizer of $\left(\bigcap_{\beta \in \Theta} \ker \beta
\right)^{\circ}$ as a Levi subgroup. Then for a subset $\Omega$ of 
$\Delta,$ the subgroup $P_{\Delta \setminus \Omega}$ is
defined over $K$ if and only if $\Omega$ is $*$-invariant and
consists entirely of distinguished roots (in other words, it is a
union of distinguished orbits). In particular, a root $\alpha \in
\Delta$ is distinguished if and only if for its $*$-orbit $\Omega$
the subgroup $P_{\Delta \setminus \Omega}$ is defined over $K$.

\vskip1mm

In the proof of Theorem E, we will need to work with the Tits
indices of a given connected absolutely simple algebraic $K$-group
$G$ over various completions of $K.$ For this purpose, we fix a
maximal $K$-torus $T$ of $G$ and a system of simple roots $\Delta
\subset \Phi(G , T).$ Given a field extension $L/K,$ we choose a
maximal $L$-torus $T'$ containing a maximal $L$-split torus $T'_s$
of $G$, and a system of simple roots $\Delta'\subset \Phi(G , T')$
determined by some coherent orderings on $X(T'_s) \otimes_{\Z} \R$
and $X(T') \otimes_{\Z} \R.$ We say that $\alpha \in \Delta$ {\it
corresponds to a distinguished vertex in the Tits index of} $G/L$ if
the root $\alpha' \in \Delta'$ corresponding to $\alpha$, under
the identification of $\Delta$ with $\Delta'$ described at the beginning of \S \ref{S:Is}, 
is distinguished. The
set of all $\alpha\in\Delta$ which correspond to distinguished
vertices in the Tits index of $G/L$ will be denoted
$\Delta^{(d)}(L).$ It follows from Lemma \ref{L:P700}(a), and the
above discussion, that $\alpha \in \Delta^{(d)}(L)$ if and only if
for the orbit $\Omega$ of $\alpha$ under the $*$-action of $\Ga(\overline{L}/L)$, 
a suitable conjugate of $P_{\Delta \setminus \Omega}$ is defined over $L.$ More generally,
for an arbitrary subset $\Omega$ of  $\Delta,$ a suitable
conjugate of $P_{\Delta \setminus \Omega}$ is defined over $L$ if
and only if $\Omega$ is invariant under the $*$-action of
$\Ga(\overline{L}/L)$  and contained in $\Delta^{(d)}(L).$ Thus,
$\mathrm{rk}_L\: G$ equals the number of orbits of the $*$-action of
$\Ga(\overline{L}/L)$ on $\Delta^{(d)}(L),$ and  $G$ is quasi-split
over $L$ if and only if $\Delta^{(d)}(L)= \Delta$.

\vskip2mm

Let $G$ be a connected absolutely simple algebraic group over a~number field
$K.$ Fix a maximal $K$-torus $T$ of $G$, and a system of simple roots $\Delta \subset \Phi(G , T)$. 
We will say that an orbit in $\Delta$, under the $*$-action of $\Ga(\overline{K}/K)$, 
is {\it distinguished everywhere} if it is contained in  $\Delta^{(d)}(K_v)$ for all $v\in V^K$. 
The following proposition, which is proved using some results of
\cite{PR5}, will not only play a crucial role in the proof of
Theorem E, it is also of independent interest.
\begin{prop}\label{P:DD1}
An orbit under the $*$-action of $\Ga(\overline{K}/K)$ on $\Delta$ is contained in $\Delta^{(d)}(K)$, i.e., it is a distinguished orbit in the Tits index of $G/K$,  if and only if it is distinguished everywhere. Therefore, $\mathrm{rk}_K\: G = r$, where $r$ is the number of orbits which are distinguished everywhere. 
\end{prop}
\begin{proof}
Without any loss of generality, we may (and do) assume that $G$ is
adjoint and $T$ contains a maximal $K$-split torus of $G$. Clearly, the distinguished
orbits in the Tits index of $G/K$ are distinguished everywhere, yielding
the inequality $\mathrm{rk}_K\: G \leqslant r.$ To prove the
opposite inequality, we can assume that $r\geqslant 1$. Let $\Omega_{i_1}, \ldots
, \Omega_{i_r}$ be the orbits in $\Delta$ which are distinguished everywhere. We will prove that 
these are precisely the distinguished orbits in the Tits index of $G/K$. For this, we set
$$
\Omega = \Omega_{i_1} \cup \cdots \cup \Omega_{i_r},
$$
and let $P_{\Delta\setminus\Omega}$ be the corresponding parabolic
subgroup. It will suffice to prove that the conjugacy class of
$P_{\Delta\setminus\Omega}$ contains a subgroup defined over $K$.
The group $G$ is an inner twist of a unique quasi-split $K$-group
$G_0$. Let $T_0$ be the centralizer of a maximal $K$-split torus
$T_0^s$ of $G_0$. Furthermore, let $\Delta_0 \subset \Phi(G_0 ,
T_0)$ be the system of simple roots with respect to some coherent
orderings on $X(T_0^s) \otimes_{\Z} \R$ and $X(T_0) \otimes_{\Z} \R$
(then, in particular, all the roots in $\Delta_0$ are
distinguished). Since $G$ is an inner twist of $G_0,$ we can pick a
$\overline{K}$-isomorphism $f \colon G_0 \to G$ so that the
associated Galois cocycle $$\sigma\mapsto \xi_{\sigma} := f^{-1}
\circ \sigma(f), \ \ \sigma \in \Ga(\overline{K}/K),$$ is of the
form $$\xi_{\sigma} = i_{g_{\sigma}},$$ where $i_z$ denotes the inner
automorphism of $G_0$ corresponding to $z \in G_0(\overline{K}),$ and 
$g:\,\sigma\mapsto g_{\sigma} $ is a Galois cocycle with values in
$G_0(\overline{K}).$ After modifying $f$ by a suitable inner automorphism, we
assume that $f(T_0) = T$ and $f^*(\Delta) =\Delta_0.$ We set 
$\Omega_0 =f^*(\Omega).$ Then for the
parabolic $K$-subgroup $P_{\Delta_0\setminus\Omega_0}$
of $G_0$, we have $f(P_{\Delta_0\setminus\Omega_0}) =
P_{\Delta\setminus\Omega}.$ Let $H_0$ be a Levi $K$-subgroup of
$P_{\Delta_0\setminus\Omega_0}$, and $\omega: H^1(K,H_0)\rightarrow H^1(K,G_0)$ be the Galois-cohomology map induced by the inclusion $H_0\hookrightarrow G_0$.

Take an arbitrary $v \in V^K.$ Then as $\Omega$ is a union of orbits in $\Delta^{(d)}(K_v)$, there
exists $a_v \in G(\overline{K}_v)$ such that
$P_{\Delta\setminus\Omega}^{(v)} := a_v P_{\Delta\setminus\Omega}
a_v^{-1}$ is defined over $K_v$. Set $b_v = f^{-1}(a_v)$ and $f_v =
f \circ i_{b_v}.$ Then $f_v(P_{\Delta_0\setminus\Omega_0}) =
P_{\Delta\setminus\Omega}^{(v)},$ and since both
$P_{\Delta_0\setminus\Omega_0}$ and
$P^{(v)}_{\Delta\setminus\Omega}$ are defined over $K_v,$ for any
$\sigma \in \Ga(\overline{K}_v/K_v),$ the automorphism $$
\xi^{(v)}_{\sigma} := f_v^{-1} \circ \sigma(f_v) = i_{b_v}^{-1}
\circ \xi_{\sigma} \circ i_{\sigma(b_v)} = i_{b_v^{-1} g_{\sigma}
\sigma(b_v)}$$ leaves $P_{\Delta_0\setminus\Omega_0}$ invariant. As
$P_{\Delta_0\setminus\Omega_0}$ coincides with its normalizer in
$G_0$ (cf.\,\cite{Bo}, Theorem 11.16), we conclude that
$b_v^{-1}g_{\sigma}\sigma(b_v)$ lies in
$P_{\Delta_0\setminus\Omega_0}(\overline{K}_v)$. Furthermore, since
the unipotent radical of any parabolic $K_v$-subgroup of a reductive $K_v$-group has trivial
Galois cohomology, we conclude that the cocycle $\sigma\mapsto
b_v^{-1}g_{\sigma}\sigma(b_v)$ is cohomologous  to a
$H_0({\overline{K}}_v)$-valued Galois cocycle ${h}^{(v)}.$ Thus, the
image of the cohomology class $x$ corresponding to the cocycle $g$,
under the restriction map $\rho_v\colon H^1(K , G_0) \to H^1(K_v ,
G_0)$, is equal to the image of  the cohomology class in
$H^1(K_v,H_0)$, corresponding to ${h}^{(v)}$, under the map $H^1(K_v
, H_0) \to H^1(K_v, G_0).$

Now, let $L$ be the minimal Galois extension of $K$ over which $G_0$
splits, and set $P = L$ if $[L : K] \neq 6,$ and let $P$ be any
cubic extension of $K$ contained in $L$ otherwise. Pick $v_0 \in V^K_f$ which does
not split in $P$ (i.e., $P \otimes_{K} K_{v_0}$ is a field). We will assume for the moment that 
$\Omega_0 \ne \Delta_0$ (or, equivalently, $\Omega \ne \Delta$). 
Then using Theorem 2 of \cite{PR5}, we easily conclude that there exists $y \in H^1(K , H_0)$
which maps to $(\rho_v(x))$ under the composite of the following two maps
$$
H^1(K , H_0) \stackrel{\omega}{\longrightarrow} H^1(K , G_0)
\stackrel{\rho=(\rho_v)}{\longrightarrow} \bigoplus_{v \neq v_0} H^1(K_v ,
G_0).$$ 
But according to Theorem 3 in \cite{PR5}, $\rho$ is injective, so
$x = \omega(y).$ This means that there exists $c
\in G_0(\overline{K})$ such that
\begin{equation}\label{E:DD5}
c^{-1} g_{\sigma} \sigma(c) \in H_0(\overline{K}) \ \ \text{for all} \ \ \sigma \in
\Ga(\overline{K}/K).
\end{equation}
We claim that the subgroup $f(c) P_{\Delta\setminus\Omega} f(c)^{-1} =
f(cP_{\Delta_0\setminus\Omega_0}c^{-1})$ is defined over $K.$ Indeed, for $\sigma
\in \Ga(\overline{K}/K)$ we have
$$
\sigma(f(cP_{\Delta_0\setminus\Omega_0}c^{-1})) = \sigma(f)(\sigma(c)
P_{\Delta_0\setminus\Omega_0} \sigma(c)^{-1})$$
$$\ \ \ \ \ \  = f(g_{\sigma}\sigma(c)
P_{\Delta_0\setminus\Omega_0} \sigma(c)^{-1} g_{\sigma}^{-1})
=f(cP_{\Delta_0\setminus\Omega_0}c^{-1})
$$
in view of (\ref{E:DD5}), proving our claim. This proves 
the proposition if $\Omega\ne \Delta$. If $\Omega = \Delta$, then, for all $v\in V^K$, $G$ is quasi-split over $K_v$, and hence is isomorphic to $G_0$ over $K_v$, which implies that $\rho_v(x)$ is trivial for all $v$. From the injectivity of $\rho$ (Theorem 3 of \cite{PR5}) we infer that $x$ is trivial, and so $G$ is isomorphic to $G_0$, and so every $*$-orbit in $\Delta$ is distinguished.    
\end{proof}

\begin{cor}\label{C:D1}
Let $G$ be an absolutely simple $K$-group of one of the following
types: $B_n$ $(n \geqslant 2),$ $C_n$ $(n \geqslant 2),$ $E_7,$
$E_8,$ $F_4$ or $G_2.$ If $G$ is isotropic over $K_v$ for all real
$v \in V^K_{\infty}$, then $G$ is isotropic over $K.$ Additionally,
if $G$ is as above, but not of type $E_7,$ then
\begin{equation}\label{E:DD100}
\rk_K\: G = \min_{v \in V^K}\,\rk_{K_v}\: G.
\end{equation}
\end{cor}
\begin{proof} The groups of these types do not have outer automorphisms, so
given any two maximal $K$-tori $T$ and $T'$ of $G,$ and systems of
simple roots $\Delta \subset \Phi(G , T)$ and $\Delta' \subset
\Phi(G , T'),$ there is a unique isomorphism between $\Phi(G , T)$
and $\Phi(G , T')$ that carries $\Delta$ to $\Delta'$. It
necessarily coincides with the canonical identification as defined
at the beginning of \S 4. Using this remark and inspecting Table II
in \cite{Ti}, we see that for the types listed in the statement, if
for every real place $v$ of $K$, $G$ is isotropic over $K_v$, then
there is a vertex in the Tits index of $G/K$ which corresponds to a
distinguished vertex in the Tits index of $G/K_v,$ for all $v \in
V^K.$ Then it follows from the proposition that this vertex is
distinguished in the Tits index of $G/K,$ and therefore $G$ is
$K$-isotropic. Moreover, if $G$ is not of type $E_7$, then it
follows from the tables in \cite{Ti} that the total number of vertices which are
distinguished in the Tits index of $G/K_v$ for all $v \in V^K$ is
$\min_{v \in V^K}\,\rk_{K_v}\: G,$ so (\ref{E:DD100}) follows from
the proposition.
\end{proof}

{\bf Proof of Theorem E.} If $G$ is of type $B_2=C_2$, $F_4$, or
$G_2$, then its Tits index over any extension $L/K$ is uniquely
determined by its $L$-rank. Therefore, since $\mathrm{rk}_{K_v}\: G_1 =
\mathrm{rk}_{K_v}\: G_2$ according to Theorem \ref{T:BC1}, and
consequently $\mathrm{rk}_K\: G_1 = \mathrm{rk}_K\: G_2$ by
Corollary \ref{C:D1}, all our assertions follow. So, we assume that
$G$ is not of any of the above three types, and in addition is
adjoint.

\vskip1mm

We pick a finite set $V_0$ of places of $K$ such that for every
$v\notin V_0$, both $G_1$ and $G_2$ are quasi-split over $K_v$. By
Theorem \ref{T:BC2}(2), we have $L_1 = L_2,$ so we can use
Proposition \ref{P:BC0} to find  maximal $K$-tori $T_i$ of $G_i$
such that $T_1$ contains a maximal $K_v$-split torus $T_{1 s}^{v}$
of $G_1$ for all $v \in V_0,$ and a $K$-isogeny (actually, a
$K$-isomorphism) $\pi \colon T_2 \to T_1$ such that
$\pi^*(\Phi(G_1,T_1)) =\Phi(G_2,T_2)$. Since $\mathrm{rk}_{K_v}\:
G_1 = \mathrm{rk}_{K_v}\: G_2$ for all $v,$ we see that $T_2$ also
contains a maximal $K_v$-split torus $T_{2s}^v$ of $G_2,$ for all
$v \in V_0.$ Notice that if we choose any system of simple roots
$\Delta_1$ in $\Phi(G_1 , T_1)$ and set $\Delta_2 = \pi^*(\Delta_1)$, 
then because $\pi^*$ commutes with the action of
$\Ga(\overline{K}/K)$ and the corresponding Weyl groups, it also
commutes with the $*$-action of $\Ga(\overline{F}/F)$ for any
extension $F/K.$ Now, let $v \in V_0,$ and let $\Delta_1^v$ be a
system of simple roots in $\Phi(G_1 , T_1)$ that corresponds to a
coherent choice of orderings on $X(T_{1 s}^v) \otimes_{\Z} \R$ and
$X(T_1) \otimes_{\Z} \R.$ Then $\Delta_2^v = \pi^*(\Delta_1^v)$
corresponds to the coherent orderings on $X(T_{2s}^v) \otimes_{\Z}
\R$ and $X(T) \otimes_{\Z} \R.$ Furthermore, since $\pi$ induces an
isomorphism between $T_{2 s}^v$ and $T_{1s}^v,$ we see that $\alpha
\in \Delta_1^v$ has nontrivial restriction to $T_{1s}^v,$ i.e.,\:it is
distinguished in the Tits index of $G_1/K_v$ if and only if
$\pi^*(\alpha)$ has nontrivial restriction to $T_{2s}^v,$ i.e.,\:it is
distinguished in the Tits index of $G_2/K_v.$ This shows that the
Tits indices of $G_2/K_v$ and $G_2/K_v$ are isomorphic for all $v
\in V_0.$ They are also isomorphic for any $v \in V^K \setminus V_0$
because then $G_1$ and $G_2$ are quasi-split, which completes the
proof of the ``local" part of Theorem E.

It remains to prove that the Tits indices of $G_1/K$ and $G_2/K$ are
isomorphic. For this, we fix a system of simple roots $\Delta_1$ of
$\Phi(G_1 , T_1)$ and set $\Delta_2 = \pi^*(\Delta_1).$ If
$\Delta'_1 \subset \Phi(G_1 , T_1)$ is another system of simple
roots and $\Delta'_2 = \pi^*(\Delta'_1)$, then the fact that $\pi^*$
commutes with the action of the corresponding Weyl groups implies
that $\pi^*$ transports the canonical identification $\Delta_1
\simeq \Delta'_1$ to the canonical identification $\Delta_2 \simeq
\Delta'_2$ (another way to see this is to observe that according to
Remark 4.5, $\pi$ extends to a $\overline{K}$-isomorphism $f \colon
G_2 \to G_1$). So, by symmetry, it is enough to prove that if
$\Omega \subset \Delta_1$ is an orbit of the $*$-action of
$\Ga(\overline{K}/K)$ which corresponds to a distinguished orbit in
the Tits index of $G_1/K$, then $\pi^*(\Omega)$ (which is also a
$*$-orbit) corresponds to a distinguished orbit in the Tits index of
$G_2/K.$ According to Proposition \ref{P:DD1}, it is enough to show
that
\begin{equation}\label{E:DD1000}
\pi^*(\Omega) \subset \Delta_2^{(d)}(K_v)
\end{equation}
for all $v \in V^K.$ As $\Delta_2^{(d)}(K_v) =
\Delta_2$ for all $v \in V^K \setminus V_0,$ we only need to
establish (\ref{E:DD1000}) for $v \in V_0.$ But since $\pi^*$
induces a bijection between distinguished vertices in $\Delta_1^v$
and $\Delta_2^v$ in the above notations, we see that
$$
\Delta_2^{(d)}(K_v) = \pi^*(\Delta_1^{(d)}(K_v)),
$$
and (\ref{E:DD1000}) follows, completing the proof of Theorem E.
\hfill $\Box$

\vskip2mm

The following interesting result is an immediate consequence of
Theorems \ref{T:BC2}(2), C, D, and E.
\begin{thm}\label{T:D2}
Let $K$ be a number field and $G$ be a connected absolutely simple
$K$-group. Let $L$ be the smallest Galois extension of $K$ over
which $G$ is an inner form of a split group. Let $\mathfrak{F}$ be
a collection of $K$-forms $G'$ of $G$ such that the set of
$K$-isomorphism classes of maximal $K$-tori of $G'$ equals the set
of $K$-isomorphism classes of  maximal $K$-tori of $G$. Then

\vskip2mm

{\rm (1)} \parbox[t]{11.5cm}{For any group belonging to
$\mathfrak{F}$, $L$ is the smallest Galois extension of $K$ over
which it is an inner form of a split group.}

\vskip2mm

{\rm(2)} \parbox[t]{11.5cm}{For any $G'\in \mathfrak{F}$, the Tits
indices of $G/K$ and $G'/K$, and for every place $v$ of $K$, the
Tits indices of $G/K_v$ and $G'/K_v$,  are isomorphic.}

\vskip2mm

{\rm (3)} \parbox[t]{11.5cm}{If $G$ is not of type $A_n$, $D_{2n+1}$, $D_4$ or
$E_6$, then every $G' \in \mathfrak{F}$ is $K$-isomorphic to $G$.}

\vskip2mm

{\rm (4)} $\mathfrak{F}$ consists of finitely many $K$-isomorphism
classes.
\end{thm}
\begin{proof} Fix $G' \in \mathfrak{F}$ and pick a finite set $S$ of places of $K$ containing all
the archimedean ones so that $\prod_{v\in S}G(K_v)$ and $\prod_{v\in
S} G'(K_v)$ are noncompact. Let $\Gamma$ and $\Gamma'$  be some
$S$-arithmetic subgroup of $G(K)$ and $G'(K),$ respectively. As $G$
and $G'$ have the same $K$-tori, it immediately follows from the
definition of weak commensurability that $\Gamma$ and $\Gamma'$ are
weakly commensurable. Now all the four assertions of the present
theorem follow from Theorems \ref{T:BC2}(2), C, D and E.
\end{proof}

\noindent {\bf Remark 7.4.} In section 9 we will show that assertion
(3) of the preceding theorem is false in general if $G$ is of type
$A_n$, $D_{2n+1}$ $(n > 1)$, or $E_6$.

\section{Lengths of closed geodesics, length-commensurable locally symmetric spaces
and Schanuel's conjecture}\label{S:G}

Let $G$ be a connected semi-simple  real algebraic group, $\mathcal{G} =
G(\R),$ and  let $\mathcal{K}$ be a maximal compact subgroup of
$\mathcal{G}.$ We let $\mathfrak{g}$ and $\mathfrak{k}$ denote the
Lie algebras of $\mathcal{G}$ and $\mathcal{K}$ respectively, and
let $\mathfrak{p}$ denote the orthogonal complement of
$\mathfrak{k}$ in $\mathfrak{g}$ relative to the Killing form
$\langle \ , \ \rangle,$ so that $\mathfrak{g} = \mathfrak{k} \oplus
\mathfrak{p}$ is a Cartan decomposition of $\mathfrak{g}.$ The
corresponding symmetric space $\mathfrak{X} = \mathcal{K} \backslash
\mathcal{G}$ is a Riemannian manifold with the metric induced by the
restriction of the Killing form to $\mathfrak{p}$ (see \cite{H} for
the details).
\vskip2mm

\noindent{\bf Positive characters.} A character $\chi$ of an $\R$-torus
$T$ is said to be {\it positive} if for every $x \in T(\R),$ the value
$\chi(x)$ is a positive real number. Any positive character of $T$
is defined over $\R.$ Given an arbitrary character $\chi \in X(T),$
the character $\chi + \bar{\chi},$ where $\bar{\chi}$ is the
character obtained by applying the complex conjugation to $\chi,$
satisfies
$$
(\chi + \bar{\chi})(x) = \chi(x)\overline{\chi(x)} = \vert \chi(x)
\vert^2
$$
for all $x \in T(\R).$ Thus, for any character $\chi$ and any $x \in
T(\R),$ the square of the absolute value of $\chi(x)$ is the value
assumed by the positive character $\chi + \bar{\chi}$ of $T$ at $x.$
\vskip1mm

Let $S$ be an $\mathbb{R}$-split torus and $T$ be a $\mathbb{R}$-torus containing $S$.
Then every character of $S$ is defined over $\mathbb{R}$. Given a character
$\alpha$ of $S$, let $\chi$ be a complex character of $T$ whose restriction to
$S$ equals $\alpha$. Then the restriction of the positive character
$\chi +\overline{\chi}$ to $S$ is $2\alpha$. Thus every character lying in the
subgroup $2X(S)$ of the character group $X(S)$ of $S$ extends to a positive
character of any $\mathbb{R}$-torus containing $S$.
\vskip2mm

Let $\mathfrak{a}$ be a Cartan subspace contained in $\mathfrak{p}$,
and $\mathcal{A}=\text{exp}\,{\mathfrak{a}}$ be the connected
abelian subgroup of $\mathcal{G}$ with Lie algebra $\mathfrak{a}$.
Let $S$ be the Zariski-closure of $\mathcal{A}$. Then $S$ is a
maximal $\mathbb{R}$-split torus of $G$ and $\mathcal{A} =
S({\mathbb{R}})^\circ$.  We fix a closed Weyl chamber
$\mathfrak{a}^+$ in $\mathfrak{a}$. Let $\{ \alpha_1, \ldots ,
\alpha_r \}$, where $r =\text{rk}_{\mathbb{R}}\,G =\text{dim}\,S$,
be the basis of the root system of $G$, with respect to $S$,
determined the Weyl chamber $\mathfrak{a}^+$, and let $\beta_i =
2\alpha_i$. Then $\beta_1, \ldots, \beta_r$ are linearly independent
positive characters.  In the sequel, we will identify $\mathfrak{a}$
with $\mathbb{R}^r$ by identifying $X\in \mathfrak{a}$ with
$(d\beta_1(X),\ldots , d\beta_r(X))$, where, for $i \in  \{1, \ldots
, r\},$ $d\beta_i$ denotes the differential of $\beta_i$ at the
identity.
\vskip1mm

We will now make some brief comments on the Lyapunov map and its
relations with weak commensurability, and will then proceed to the
core issue of the lengths of closed geodesics and
length-commensurable locally symmetric spaces.

\vskip2mm

\noindent{\bf Lyapunov map.} For an element $g \in \cG,$ we let $g =
g_sg_u$ be its Jordan decomposition. For simplicity, we denote the semi-simple component
$g_s$ by $s$. Let $T$ be a maximal $\mathbb{R}$-torus of $G$ containing
$s$. Let $\mathcal{C}$ be the maximal compact subgroup of
$T(\mathbb{R})$ and $T_s$ be the maximal $\mathbb{R}$-split subtorus
of $T$. Then $T({\mathbb{R}})$ is a direct product of $\mathcal{C}$
and $T_s({\mathbb{R}})^\circ$, so we can write $s = s_e\cdot s_h$,
with $s_e\in \mathcal{C}$, and $s_h \in T_s ({\mathbb{R}})^\circ$.
The elements $s_e$ and $s_h$ are called the {\it elliptic} and the
{\it hyperbolic}  components of $s$ (or of $g$). There is an element
$z\in\mathcal{G}$ which conjugates $\mathcal{C}$ into $\mathcal{K}$
and $T_s ({\mathbb{R}})^\circ$ into $\mathcal{A}$ such that
$zs_hz^{-1} = \text{exp}\,X$, with $X\in \mathfrak{a}^+$. The
element $X$ is the unique element of $\mathfrak{a}^+$ such that the
hyperbolic component $s_h$ of $g$ is a conjugate of $\text{exp}\,X$,
and we will denote it by $\ell (g)$. Thus we get a map (the Lyapunov
map) $\ell : \mathcal{G}\rightarrow {\mathfrak{a}}^+$. Clearly, for
any $g \in \cG$ we have $\ell(g) = \ell(g_s),$ and moreover,
for any positive integer $n$, $\ell (g^n) = n\ell (g)$. \vskip1mm

Continuing with the above notations, we
let $\chi_i,$ for $i \in \{1, \ldots , r\},$ be the unique positive
character of $T$ extending the character $\text{Int}\,z^{-1}\cdot
\beta_i |_{T_s}$, and let $d\chi_i$ denote its differential at the
identity. Since $\chi_i(s)=\chi_i(s_h),$ we have
$$\ell(s) = (d\chi_1(\text{Ad}\,z^{-1}(X)),\ldots, d\chi_r(\text{Ad}\,z^{-1}(X)))
=(\text{log}\,\chi_1(s),\ldots , \text{log}\,\chi_r(s)).$$

\vskip1mm

For a subgroup $\Gamma$ of $\mathcal{G}$, let $\Gamma^{\mathrm{ss}}$ denote the
set of semi-simple elements of $\Gamma$.  From the above description of the
Lyapunov map,  the following proposition is obvious.

\begin{prop}\label{P:G-0} If $\Gamma_1$ and $\Gamma_2$ are two subgroups
of $\mathcal{G}$ such that $\mathbb{Q}\cdot\ell (\Gamma_1^{\mathrm{ss}})=\mathbb{Q}\cdot
\ell(\Gamma_2^{\mathrm{ss}})$, then $\Gamma_1$ and $\Gamma_2$ are weakly commensurable.
\end{prop}

If $\Gamma$ is an arithmetic subgroup of $\mathcal{G}$ and $g
\in \Gamma$, then there exists an integer $n =
n(g)$ such that $g_u^n\in \Gamma$. Then $g_s^n$ lies in
$\Gamma$. On the other hand, if $\Gamma$ is an irreducible nonarithmetic lattice of $\mathcal{G}$ (then $\mathcal{G}$ is of $\R$-rank 1), then it can be shown that there exists a positive integer $n=n(\Gamma)$ such that for every non-semi-simple element  $g$ of $\Gamma$, $g^n$ is unipotent. We conclude that if $\Gamma$ is an arithmetic or nonarithmetic lattice of $\mathcal{G}$, then $\Q\cdot\ell(\Gamma)
=\Q\cdot\ell(\Gamma^{\mathrm{ss}})$. \vskip2mm

\noindent {\bf Lengths of closed geodesics on locally symmetric
spaces.} Given a discrete torsion-free subgroup $\Gamma$ of
$\mathcal{G},$ the quotient $\mathfrak{X}_{\Gamma} :=
\mathfrak{X}/\Gamma$ is a Riemannian locally symmetric space. We
first need to recall some facts about closed geodesics in
$\mathfrak{X}_{\Gamma},$ and in particular the formula for their
length, given in \cite{PR3}. Closed geodesics in $\fX_{\Gamma}$
correspond to semi-simple elements in $\Gamma,$ and are obtained by
a construction similar to the one used to define the Lyapunov map.
More precisely, let $\gamma$ be a fixed semi-simple element of
$\Gamma,$ and let $T$ be a maximal $\R$-torus of $G$ containing
$\gamma.$ As we mentioned above,
$T(\R)$ is a direct product of $\mathcal{C}$ and $T_s(\R)^{\circ},$
where $\mathcal{C}$ is the maximal compact subgroup of $T(\R)$ and
$T_s$ is the maximal $\R$-split subtorus of $T.$ Take {\it any} $z
\in \mathcal{G}$ such that $zTz^{-1}$ is invariant under the Cartan
involution associated with the decomposition $\mathfrak{g} =
\mathfrak{k} \oplus \mathfrak{p},$ and consequently
\begin{equation}\label{E:G2007}
z\mathcal{C}z^{-1} \subset \mathcal{K} \ \ \text{and} \ \
zT_s(\R)^{\circ}z^{-1} \subset \mathrm{exp}\: \mathfrak{p}.
\end{equation}
Thus, here we do not require the inclusion $zT_s(\R)^{\circ}z^{-1}
\subset \exp \mathfrak{a}^+,$ however, all the $z$'s satisfying
(\ref{E:G2007}) lie in the same coset modulo $\cK,$ and therefore
define the same point $\cK z \in \fX.$ So, if we write $\gamma =
\gamma_e \cdot \gamma_h$ with $\gamma_e \in \mathcal{C}$ and
$\gamma_h \in T_s(\R)^{\circ},$ and then $\gamma_h = z^{-1}
\mathrm{exp}(X) z$ for some $X \in \mathfrak{p}$  that commutes with
$z\gamma_ez^{-1};$ moreover, it follows from the above discussion
that $X$ is  a conjugate of $\ell (\gamma)$ under an element of
$\text{Ad}\,\mathcal{K}$. With these notations, the curve
$\tilde{c}_{\gamma}$ parametrized by $\widetilde{\varphi} \colon t
\mapsto \mathcal{K}\mathrm{exp}(tX)z$ for $t \in \R,$ is a geodesic
on $\mathfrak{X}$ which passes through the point $\mathcal{K}z$.
Furthermore,
$$
\widetilde{\varphi}(t) \cdot \gamma = \mathcal{K}\mathrm{exp}(tX) \cdot
z \gamma_e z^{-1} \cdot z\gamma_h z^{-1} \cdot z = \mathcal{K}
\mathrm{exp}(tX) \cdot \mathrm{exp}(X) \cdot z = \widetilde{\varphi}(t +
1),
$$
implying that the map $\varphi \colon \R \to \mathfrak{X}_{\Gamma},$
obtained by composing $\widetilde{\varphi}$ with the natural
map $\pi \colon \mathfrak{X} \to \mathfrak{X}_{\Gamma},$ is periodic
with period 1, and hence its smallest period is of the form
$1/n_{\gamma}$ for some integer $n_{\gamma} \geqslant 1.$ It follows
that the image $c_{\gamma}$ of $\tilde{c}_{\gamma}$ in
$\mathfrak{X}_{\Gamma}$ is a closed geodesic, and since
$$
\langle \varphi'(t) , \varphi'(t) \rangle = \langle
{\widetilde{\varphi}}'(t) , {\widetilde{\varphi}}'(t)\rangle = \langle X , X
\rangle,
$$
for all $t \in \R,$ we see that the length of $c_{\gamma}$ is
$(1/n_{\gamma}) \langle X , X \rangle.$
\begin{prop}\label{P:G-1}
{\rm (i)} Every closed geodesic in $\mathfrak{X}_{\Gamma}$ is of the
form $c_{\gamma}$ for some semi-simple $\gamma \in \Gamma.$

\vskip3mm

\noindent {\rm (ii)} \parbox[t]{11.8cm}{\baselineskip=5mm The length
of $c_{\gamma}$ is $(1/n_{\gamma})\lambda_{\Gamma}(\gamma)$ where
$n_{\gamma}$ is an integer $\geqslant 1$ and
$\lambda_{\Gamma}(\gamma)$ is given by the following formula:
\begin{equation}\label{E:G-1}
\lambda_{\Gamma}(\gamma)^2 =\langle \ell(\gamma),\ell (\gamma)\rangle = \left(\sum (\log \vert \alpha(\gamma)
\vert)^2\right),
\end{equation}
where the summation is over all roots of $G$ with respect to $T$ and
$\log$ denotes the natural logarithm.}

\vskip3mm

\noindent Thus, $$\Q \cdot L(\mathfrak{X}_{\Gamma}) = \Q \cdot \{
\lambda_{\Gamma}(\gamma) \ \vert \ \gamma \in \Gamma \ \
\text{semi-simple} \},$$ where $\lambda_{\Gamma}(\gamma)$ is given
by~{\rm (\ref{E:G-1})}.
\end{prop}
\begin{proof} (i) Any closed geodesic $c$ in $\mathfrak{X}_{\Gamma}$ is obtained as
the image under $\pi$ of a geodesic $\tilde{c}$ in $\mathfrak{X}.$
Fix a point $\mathcal{K}z \in \tilde{c}.$ It is known that
$\tilde{c}$ admits a parametrization of the form
$$\widetilde{\varphi}(t) = \mathcal{K} \exp(tX)z$$ for some $X \in
\mathfrak{p}$ (cf.\:\cite{H}, Theorem 3.3(iii) in Ch.\:IV). After
replacing $X$ by a suitable positive-real multiple, we can assume
that $\pi(\widetilde{\varphi}(0)) = \pi(\widetilde{\varphi}(1)),$ and
$d_{\widetilde{\varphi}(0)}\pi({\widetilde{\varphi}}'(0)) =
d_{\widetilde{\varphi}(1)}\pi({\widetilde{\varphi}}'(1))$. Then, in
particular, $\widetilde{\varphi}(1) = \widetilde{\varphi}(0) \gamma$ for
some $\gamma \in \Gamma.$ Since the map
$$
\mathcal{K} \times \mathfrak{p} \to \mathcal{G}, \ \ (\kappa , Y)
\mapsto \kappa \exp(Y),
$$
is a diffeomorphism, the element $z\gamma z^{-1}$ can be uniquely
written in the form $z\gamma z^{-1} = \kappa \exp(Y).$ Then
$\tilde{c}(1) = \tilde{c}(0)\gamma$ yields $X = Y,$ i.e.,
\begin{equation}\label{E:G100}
z\gamma z^{-1} = \kappa \exp(X).
\end{equation}
Furthermore, the curves in $\mathfrak{X}$ with the parametrizations
$$
\widetilde{\varphi}_1(t) = \widetilde{\varphi}(t) \cdot \gamma \ \
\text{and} \ \ \widetilde{\varphi}_2(t) = \tilde{\varphi}(t + 1)
$$
are both geodesics in $\mathfrak{X}$ such that
$$
\widetilde{\varphi}_1(0) = \widetilde{\varphi}(0) \cdot \gamma =
\widetilde{\varphi}(1) = \widetilde{\varphi}_2(0) =: p.
$$
Since $\pi(\widetilde{\varphi}_1(t)) = \pi(\widetilde{\varphi}(t)),$
we have
$$
d_p\pi({\widetilde{\varphi}}'_1(0)) =
d_{\widetilde{\varphi}(0)}\pi({\widetilde{\varphi}}'(0)) =
d_{\widetilde{\varphi}(1)}\pi({\widetilde{\varphi}}'(1)) =
d_p\pi({\widetilde{\varphi}}'_2(0)).
$$
Thus, ${\widetilde{\varphi}}'_1(0) = {\widetilde{\varphi}}'_2(0),$ hence by
the uniqueness of a geodesic through a given point in a given
direction, we get $\widetilde{\varphi}_1(t) = \widetilde{\varphi}_2(t)$ for
all $t.$ Combining the definitions of $\widetilde{\varphi},$
$\widetilde{\varphi}_1$ and $\widetilde{\varphi}_2$ with  (\ref{E:G100}), we
now obtain that
$$
\mathcal{K}\exp(tX)\kappa = \mathcal{K}\exp(t(\mathrm{Ad}\:
\kappa^{-1}(X))) = \mathcal{K}\exp(tX),
$$
which implies that $\kappa$ commutes with $\exp(tX)$ for all $t.$
Since the elements $\kappa$ and $\exp(X)$ are semi-simple, we
conclude that $\gamma = z^{-1}(\kappa \exp(X))z$ is semi-simple.
Moreover, $\kappa$ and $\exp(X)$ are contained in a maximal
$\R$-torus $T_0$ of $G$ which is invariant under the Cartan involution. Let $T = z^{-1}T_0z$. 
Then $T(\R) = z^{-1}T_0(\R)z$ contains $\gamma,$ and $\gamma_e = z^{-1}\kappa
z$ and $\gamma_h = z^{-1}\exp(X)z$ in the notations introduced prior
to the statement of the proposition. It is now obvious that $c$
coincides with the geodesic $c_{\gamma}.$ As we already explained,
its length is $(1/n_{\gamma})\langle X , X \rangle^{1/2},$ where
$n_{\gamma}$ is the integer $\geqslant 1$ such that $1/n_{\gamma}$
is the smallest positive period of $\varphi(t) =
\pi(\tilde{\varphi}(t)).$

\vskip3mm

(ii) We need to show that $\lambda_{\Gamma}(\gamma) := \langle X , X
\rangle^{1/2}$ ($=\langle \ell(\gamma),\ell(\gamma)\rangle^{1/2}$)
is given by the equation (\ref{E:G-1}). Since the Killing form is
invariant under the adjoint action of ${\mathcal G}$ on $\mathfrak
g$, we have $\langle X , X \rangle = \langle X' , X' \rangle$, where
$X' = {\rm Ad} z^{-1} (X)$ so that $\gamma_h = \exp(X').$ In a
suitable basis of $\mathfrak{g},$ $\mathrm{Ad}\: \gamma_h$ is
represented by a diagonal matrix whose diagonal entries are $1$
(repeated $\dim T$ times) and $\alpha(\gamma_h)$ for all $\alpha \in
\Phi(G , T);$ notice that all these numbers are real and positive.
In the same basis, $\mathrm{ad}\: X'$ is represented by a diagonal
matrix with the diagonal entries $0$ (repeated $\dim T$ times) and
$d\alpha(X')$ for all $\alpha \in \Phi(G , T).$ For every $\alpha$
we clearly have
$$
\vert \alpha(\gamma) \vert = \vert \alpha(\gamma_h) \vert =
\exp(d\alpha(X')).
$$
So,
$$
\langle X, X\rangle = \langle X' , X' \rangle = \sum_{\alpha \in
\Phi(G , T)} (d\alpha(X'))^2 = \sum_{\alpha \in \Phi(G , T)} (\log
\vert \alpha(\gamma) \vert)^2,
$$
and (\ref{E:G-1}) follows.
\end{proof}

In order to relate the notion of length-commensurability  with that
of weak commensurability, we need to recast formula
(\ref{E:G-1}) in a slightly different form. As a root $\alpha$ of
$G$ with respect to $T$ is a character of $T$, $|\alpha(\gamma)|^2$
is the value assumed by a positive character of $T$, and therefore,
\begin{equation}\label{E:G-2}
\lambda_{\Gamma}(\gamma)^2 = \sum_{i = 1}^p s_i (\log
\chi_{i}(\gamma))^2 ,
\end{equation}
where $\chi_1, \ldots , \chi_p$ are certain positive characters of
$T$ and $s_1, \ldots , s_p$ are positive rational numbers (whose
denominators are divisors of 4).

\vskip1mm

We will now elaborate on (\ref{E:G-2}) in the rank one case.
\begin{lemma}\label{L:G2007}
Assume that $\mathrm{rk}_{\R}\: G = 1,$ and let $\Gamma$ be a
discrete torsion-free subgroup of $\cG = G(\R).$ Let $\gamma \in
\Gamma$ be a semi-simple element $\neq 1,$ and let $T$ be a maximal
$\R$-torus containing it. Then

\vskip2mm

{\rm (1)} \parbox[t]{11.5cm}{$\mathrm{rk}_{\R}\: T = 1,$ so the
group of positive characters of $T$ is cyclic with a generator, say,
$\chi.$}

\vskip1mm

{\rm (2)} \parbox[t]{11.5cm}{$\chi(\gamma) \neq 1.$}

\vskip1mm

{\rm (3)} \parbox[t]{11.5cm}{There exists $t > 0,$ depending only on
$G,$ but not on $\gamma,$ $\Gamma$ or $T$ such that
$$\lambda_{\Gamma}(\gamma) = t \vert \log \chi(\gamma) \vert.$$}
\end{lemma}
\begin{proof}
(1): $\mathrm{rk}_{\R}\: T = 0$ would imply that $T(\R)$ is compact,
so the discreteness of $\langle \gamma \rangle$ would imply its
finiteness. Since $\Gamma$ is torsion-free, we would get $\gamma =
1,$ a contradiction.

\vskip1mm

(2): Proved similarly using the fact that $(\ker \chi)(\R)$ is
compact.

\vskip1mm

(3): This follows from (\ref{E:G-1}) and (\ref{E:G-2}) combined with
the fact that any two maximal $\R$-tori of $G$ having real rank one
are conjugate under an element of $\cG.$
\end{proof}

\begin{cor}\label{C:G2007}
Assume that $\mathrm{rk}_{\R}\: G = 1.$ Let $K$ be a number field
contained in $\R,$ and assume that $G_1$ and $G_2$ are two $K$-forms
of $G$ having the same set of $K$-isomorphism classes of maximal
$K$-tori. Furthermore, for $i = 1 , 2,$ let $\Gamma_i$ be a discrete
torsion-free $(G_i , K)$-arithmetic subgroup of $\cG.$ Then
\begin{equation}\label{E:G2008}
\Q \cdot \lambda_{\Gamma_1}(\Gamma_1^{ss}) = \Q \cdot
\lambda_{\Gamma_2}(\Gamma_2^{ss}),
\end{equation}
and consequently, $\fX_{\Gamma_1}$ and $\fX_{\Gamma_2}$ are
length-commensurable.
\end{cor}
Indeed, let $\gamma_1 \in \Gamma_1^{ss} \setminus \{ 1 \},$ and let
$T_1$ be a maximal $K$-torus of $G_1$ containing $\gamma_1.$ By our
assumption, for a suitable maximal $K$-torus $T_2$ of $G_2,$ there
exists a $K$-isomorphism $\varphi \colon T_1 \to T_2.$ Since
$\varphi(T_1(K) \cap \Gamma_1)$ is an arithmetic subgroup of $T_2(K),$
there exists $n
> 0$ such that $\gamma_2 := \varphi(\gamma_1)^n \in T_2(K) \cap
\Gamma_2.$ Let $\chi^{(1)}$ be a generator of the group of positive
characters of $T_1$ (cf. Lemma \ref{L:G2007})(1)). Then $\chi^{(2)}
:= (\varphi^*)^{-1}(\chi_1)$ is a generator of the group of positive
characters of $T_2,$ and $\chi^{(2)}(\gamma_2) =
\chi^{(1)}(\gamma_1)^n.$ It follows from Lemma \ref{L:G2007}(3) that
$$
\vert \lambda_{\Gamma_2}(\gamma_2)/\lambda_{\Gamma_1}(\gamma_1)
\vert = n,
$$
yielding the inclusion $$\Q \cdot \lambda_{\Gamma_1}(\Gamma_1^{ss})
\subset \Q \cdot \lambda_{\Gamma_2}(\Gamma_2^{ss}).$$ By symmetry,
we get (\ref{E:G2008}). The last assertion follows from
(\ref{E:G2008}) and Proposition \ref{P:G-1}.

\vskip2mm

To deal with the higher rank case, we need the following.
\begin{lemma}\label{L:G-1}
Let $\gamma_1 ,
\gamma_2 \in G(\R)$ be two semi-simple elements contained in the
maximal \,$\R$-tori $T_1$ and $T_2$ of $G,$ respectively. Given two
collections of characters $\chi_1^{(1)}, \ldots , \chi_{d_1}^{(1)}
\in X(T_1)$ and $\chi_1^{(2)}, \ldots , \chi_{d_2}^{(2)} \in
X(T_2),$ we set
$$S_i = \{\log \vert \chi_1^{(i)}(\gamma_i)\vert, \ldots , \log
\vert \chi_{d_i}^{(i)}(\gamma_i)\vert \}.$$ If  $\gamma_1 ,
\gamma_2$ are not weakly commensurable and each of the sets (of real
numbers) $S_1$ and $S_2$ is linearly independent over $\Q,$ then so
is their union $S_1 \cup S_2.$
\end{lemma}
\begin{proof}
According to the above discussion, there exist positive characters
$\theta_1^{(1)}, \ldots ,  \theta_{d_1}^{(1)} \in X(T_1)$ and
$\theta_1^{(2)}, \ldots , \theta_{d_2}^{(2)} \in X(T_2)$ such that
$$
\theta_j^{(i)}(x) = \vert \chi_j^{(i)}(x) \vert^2 \ \ \ \text{for
all} \ \ x \in T_i(\R).
$$
If the set $S_1 \cup S_2$ is linearly dependent over $\Q,$ there
exist integers $s_1, \ldots , s_{d_1}$, $ t_1, \ldots , t_{d_2},$
not all zero,  such that
$$
s_1\log \theta_1^{(1)}(\gamma_1) + \cdots + s_{d_1}\log
\theta_{d_1}^{(1)}(\gamma_1) + t_1\log \theta_1^{(2)}(\gamma_2) +
\cdots + t_{d_2} \log \theta_{d_2}^{(2)}(\gamma_2) = 0.
$$
Consider the characters
$$\psi_1 = s_1\theta_1^{(1)} + \cdots + s_{d_1}\theta_{d_1}^{(1)} \  \
\text{of} \  \  T_1\  \  \text{and} \ \   \psi_2 =
-(t_1\theta_1^{(2)} +\cdots + t_{d_2}\theta_{d_2}^{(2)}) \ \
\text{of}  \  \  T_2.$$  Then $\psi_1(\gamma_1) = \psi_2(\gamma_2),$
and hence,
$$
\psi_1(\gamma_1) = 1 = \psi_2(\gamma_2)
$$
because $\gamma_1$ and $\gamma_2$ are not commensurable. This means
that
$$
s_1\log \theta_1^{(1)}(\gamma_1) + \cdots + s_{d_1}\log
\theta_{d_1}^{(1)}(\gamma_1) = 0 =  t_1\log \theta_1^{(2)}(\gamma_2)
+ \cdots + t_{d_2} \log \theta_{d_2}^{(2)}(\gamma_2),
$$
and therefore all the coefficients are zero because the sets $S_1$
and $S_2$ are linearly independent.
\end{proof}

Some of our results depend on the validity of Schanuel's conjecture
in transcendental number theory (cf.\,\,\cite{A}), and we recall
here its statement.

\vskip2mm

\noindent {\bf Schanuel's conjecture.} {\it If $z_1, \ldots , z_n
\in \C$ are linearly independent over $\Q,$ then the transcendence
degree (over $\Q$) of the field generated by
$$
z_1, \ldots , z_n; \  e^{z_1}, \ldots , e^{z_n}
$$
is $\geqslant n.$}

\vskip2mm

\noindent We will only use the fact that the truth of this
conjecture implies that for algebraic numbers $z_1, \ldots , z_n,$
(any values of) their logarithms
$$
\log z_1, \ldots , \log z_n
$$
are algebraically independent once they are linearly independent
(over $\Q$).
\begin{prop}\label{P:G-2}
Let $G$ be a connected semi-simple real algebraic subgroup of $\mathrm{SL}_n$
and $\mathcal{G}=G(\R)$. Let $\Gamma_1$, $\Gamma_2$ be two discrete
torsion-free subgroups of $\mathcal{G}$.  Suppose that nontrivial
semi-simple elements $\gamma_1 \in \Gamma_1$ and $\gamma_2 \in
\Gamma_2$ are not weakly commensurable. Then

\vskip2mm

{\rm \ (i)} \parbox[t]{11.5cm}{If $\mathrm{rk}_{\R} \: G = 1$, then
$\theta = \lambda_{\Gamma_1}(\gamma_1)/\lambda_{\Gamma_2}(\gamma_2)$
is irrational. Moreover, if there exists a number field $K$ such
that $\Gamma_1$ and $\Gamma_2$ can be conjugated into
$\mathrm{SL}_n(K)$, then $\theta$ is transcendental over $\Q$.}

\vskip2mm

{\rm (ii)} \parbox[t]{11.5cm}{If there exists a number field $K$
such that $\Gamma_1$ and $\Gamma_2$ can be conjugated into
$\mathrm{SL}_n(K)$, and Schanuel's conjecture holds, then
$\lambda_{\Gamma_1}(\gamma_1)$ and $\lambda_{\Gamma_2}(\gamma_2)$
are algebraically independent over $\Q.$}
\end{prop}
\begin{proof}
We fix maximal $\R$-tori $T_1$ and $T_2$ of $G$ which contain
$\gamma_1$ and $\gamma_2$ respectively.

(i) Using Lemma \ref{L:G2007}, (1) and (2), for $i = 1 , 2,$ we can
pick a generator $\chi^{(i)}$ of the group of positive characters of
$T_i$ so that $\chi^{(i)}(\gamma_i) > 1$ for $i = 1 , 2.$ Then by
Lemma \ref{L:G2007}(3) we have
$$
\lambda_{\Gamma_i}(\gamma_i) = t \log \chi^{(i)}(\gamma_i).
$$
Since the elements $\gamma_1$ and $\gamma_2$ are not weakly
commensurable, for every nonzero integers $m$, $n$, we have
$$
\chi^{(1)}(\gamma_1)^m \neq \chi^{(2)}(\gamma_2)^n,
$$
i.e., the ratio
$$
\theta=\frac{\lambda_{\Gamma_1}(\gamma_1)}{\lambda_{\Gamma_2}(\gamma_2)}=
\frac{\log \chi^{(1)}(\gamma_1)}{\log
\chi^{(2)}(\gamma_2)}
$$
is irrational. If there exists a number field $K$ such that $\Gamma_1$ and $\Gamma_2$ can be conjugated into
$\mathrm{SL}_n(K),$  then the numbers $\chi^{(i)}(\gamma_i)$ are
algebraic, and therefore by a theorem proved independently by
Gel'fond and Schneider in 1934 (cf.\,\cite{Ba}), $\theta$ is
transcendental over $\Q.$

\vskip2mm

(ii) According to (\ref{E:G-2}), we have the following expressions
$$
\lambda_{\Gamma_1}(\gamma_1)^2 = \sum_{i = 1}^{p} s_i^{(1)} (\log
\chi_i^{(1)}(\gamma_1))^2 \ \ \ \text{and} \ \ \
\lambda_{\Gamma_2}(\gamma_2)^2 = \sum_{i = 1}^{p} s_i^{(2)} (\log
\chi_i^{(2)}(\gamma_2))^2
$$
After renumbering the characters, we can assume that
$$
a_1 := \log \chi_1^{(1)}(\gamma_1), \ldots , a_{m_1} := \log
\chi_{m_1}^{(1)}(\gamma_1)
$$
$$
(\text{resp.,} \ b_1 := \log \chi_1^{(2)}(\gamma_2), \ldots ,
b_{m_2} = \log \chi^{(2)}_{m_2}(\gamma_2))
$$
for some $m_1 , m_2 \leqslant p,$ form a basis of the $\Q$-subspace
of $\R$ spanned by $\log \chi_i^{(1)}(\gamma_1)$ (resp., $\log
\chi_i^{(2)}(\gamma_2)$) for $i \leqslant p$ (notice that $m_1 , m_2
\geqslant 1$ as otherwise the length of the corresponding geodesic
would be zero, which is impossible). It follows from Lemma
\ref{L:G-1} that the numbers
$$
a_1, \ldots , a_{m_1};\  b_1, \ldots , b_{m_2}
$$
are linearly independent over $\Q.$ Since by our assumption the
subgroups $\Gamma_1$ and $\Gamma_2$ can be conjugated into
$\mathrm{SL}_n(K),$ the values $\chi_i^{(j)}(\gamma_j)$ are algebraic
numbers, so it follows from Schanuel's conjecture  that $a_1, \ldots
, a_{m_1};\,\, b_1, \ldots , b_{m_2}$ are algebraically independent
over $\Q.$ It remains to observe that
$\lambda_{\Gamma_1}(\gamma_1)^2$ and
$\lambda_{\Gamma_2}(\gamma_2)^2$ are given by nonzero homogeneous
polynomials of degree two, with rational coefficients, in $a_1,
\ldots , a_{m_1}$ and $b_1, \ldots , b_{m_2},$ respectively, and
therefore they are algebraically independent.
\end{proof}

By combining Propositions \ref{P:G-1} and \ref{P:G-2} we obtain the
following:
\begin{thm}\label{T:G-3} Let $\Gamma_1$, $\Gamma_2$ be discrete torsion-free
subgroups of
$\mathcal{G}$. If $\Gamma_1$ and $\Gamma_2$ are not weakly
commensurable, then, possibly after interchanging them, the
following assertions hold.

\vskip2mm

{\rm \ (i)} \parbox[t]{11.5cm}{If $\mathrm{rk}_{\R}\: G = 1$,  then there
exists $\lambda_1 \in L(\mathfrak{X}_{\Gamma_1})$ such that for {\rm
any} $\lambda_2 \in L(\mathfrak{X}_{\Gamma_2}),$ the ratio
$\lambda_1/\lambda_2$ is irrational.}

\vskip2mm

{\rm (ii)} \parbox[t]{11.5cm}{If there exists a number field $K$
such that both $\Gamma_1$ and $\Gamma_2$ can be conjugated into
$\mathrm{SL}_n(K)$, and Schanuel's conjecture holds, then there
exists $\lambda_1 \in L(\mathfrak{X}_{\Gamma_1})$ which is
algebraically independent from {\rm any} $\lambda_2 \in
L(\mathfrak{X}_{\Gamma_2}).$}

\vskip2mm

\noindent In either case, (under the above assumptions)
$\fX_{\Gamma_1}$ and $\fX_{\Gamma_2}$ are not length-commensurable.
\end{thm}

\vskip3mm

If $\cG$ does not contain any nontrivial connected compact normal subgroups, and it is not locally isomorphic to either ${\rm SL}_2(\R)$ or ${\rm SL}_2(\C)$, and $\Gamma$ is an irreducible lattice in $\cG$, then there exists a real number field $K$ such that $\Gamma$ can be conjugated into ${\rm SL}_n(K)$, see \cite{Ra}, Proposition 6.6.

\vskip3mm

\noindent {\it The results in the rest of this section for locally
symmetric spaces of rank $> 1$ assume the truth of Schanuel's
conjecture.} \vskip2mm

\vskip3mm

Henceforth, we will  study locally symmetric spaces of $\cG =
G(\R)$, where $G$ is an {\it absolutely simple} real algebraic
group. It follows from Theorem \ref{T:G-3} that
length-commensurability of the locally symmetric spaces
$\mathfrak{X}_{\Gamma_1}$ and $\mathfrak{X}_{\Gamma_2}$ implies weak
commensurability of the subgroups $\Gamma_1$ and $\Gamma_2.$ On the
other hand, commensurability of $\Gamma_1$ and $\Gamma_2$ up to an
$\R$-automorphism of $G$ is equivalent to commensurability of
$\mathfrak{X}_{\Gamma_1}$ and $\mathfrak{X}_{\Gamma_2}.$ Now Theorem F 
immediately implies the following.
\begin{thm}\label{T:G25}
If $\mathfrak{X}_{\Gamma_1}$ and $\mathfrak{X}_{\Gamma_2}$ are of
finite volume, length-commensurable, and $\Gamma_1$ is arithmetic,
then so is $\Gamma_2$.
\end{thm}

\vskip1mm

We will now focus on {\it arithmetically defined} locally symmetric
spaces. Using the above observation and applying Theorems C and D,
we obtain the following.
\begin{thm}\label{T:LS1}
Each class of length-commensurable arithmetically defined locally
symmetric spaces of $\cG = G(\R)$ is a union of finitely many
commensurability classes. It in fact consists of a single
commensurability class if $G$ is not of type $A_n$ $(n > 1),$ $D_{2n+1}$ $(n\geqslant 1)$,
$D_4$ or $E_6.$
\end{thm}

To see what this theorem means for hyperbolic spaces, we recall that
the  even-dimensional real hyperbolic space $\mathbf{H}^{2n}$ is the
symmetric space of a group of type $B_n,$ the odd-dimensional real
hyperbolic space $\mathbf{H}^{2n-1}$ - of a group of type $D_n,$
the complex hyperbolic space $\mathbf{H}^n_{\C}$ - of a group of
type $A_n,$ and the quaternionic hyperbolic space
$\mathbf{H}_{\mathbb{H}}^n$ - of a group of type $C_{n+1}.$ All these
spaces are of rank one.  Using Theorem C and Proposition  \ref{P:G-2}(i), 
we obtain the following result.

\begin{cor}\label{C:LS1}
Let $M$ be either the real hyperbolic space
$\mathbf{H}^{2n}$, or  $\mathbf{H}^{4n+3}$, or the quaternionic hyperbolic space
$\mathbf{H}_{\mathbb{H}}^n$, for any $n\geqslant 1$, and let $M_1$ and
$M_2$ be two arithmetic quotients of $M.$ If $M_1$ and $M_2$ are not
commensurable, then after a possible interchange of $M_1$ and $M_2,$
there exists $\lambda_1 \in L(M_1)$ such that for any $\lambda_2 \in
L(M_2),$ the ratio $\lambda_1/\lambda_2$ is transcendental over
$\Q.$
\end{cor}

\noindent  {\bf Remark 8.11.} In Example 6.6, we indicated that for
the $\R$-group $G = \mathrm{SL}_{2 , \mathbb{H}},$ one can construct
two anisotropic $\Q$-forms $G_1$ and $G_2$ that have the same set of
$\Q$-isomorphism classes of maximal $\Q$-tori. For $i = 1 , 2,$ fix
a torsion-free $(G_i , \Q)$-arithmetic subgroup $\Gamma_i$ of $\cG.$
Since $G \simeq \mathrm{Spin}(q),$ where $q$ is a real quadratic
form of signature $(5 , 1),$ the corresponding symmetric space $\fX$
is $\mathbf{H}^5.$ Using Corollary \ref{C:G2007}, we now conclude
that $\fX_{\Gamma_1}$ and $\fX_{\Gamma_2}$ are length-commensurable,
but noncommensurable, compact hyperbolic 5-manifolds. A similar
argument applied to a suitable modification of Example 6.6 enables
one to construct examples of noncommensurable length-commensurable
complex hyperbolic manifolds of any even dimension. These examples
will be subsumed by general constructions in \S 9, which in
particular, allow one to construct examples of this nature for real
hyperbolic manifolds of any dimension of the form $4n + 1,$ and for
complex hyperbolic manifolds of any dimension, cf.\:9.14.


\vskip2mm

\setcounter{thm}{11}

We now recall that given a discrete $(G_i , K_i)$-arithmetic
subgroup $\Gamma_i \subset \cG,$ the compactness of the quotient
$\cG/\Gamma_i,$ and hence of the locally symmetric subspace
$\mathfrak{X}_{\Gamma_i},$ is equivalent to $G_i$ being
$K_i$-anisotropic (cf.\:\cite{PlR}, Theorem 4.17). Combining this
with Theorem E, we obtain the following.
\begin{thm}\label{T:G2010}
Let $\fX_{\Gamma_1}$ and $\fX_{\Gamma_2}$ be two arithmetically
defined locally symmetric spaces of the same absolutely simple real
Lie group $\cG$. If they are length-commensurable, then the
compactness of one of them implies the compactness of the other.
\end{thm}

We close this section with a result which applies also to
nonarithmetic subgroups.
\begin{thm}\label{T:G2011}
Let $\fX_{\Gamma_1}$ and $\fX_{\Gamma_2}$ be two locally symmetric
spaces of the same absolutely simple real Lie group $\cG$, modulo torsion-free lattices
$\Gamma_1$ and $\Gamma_2.$ Denote by $K_{\Gamma_i}$ the field
generated by the traces ${\tr}\mathrm{Ad}\: \gamma$ for $\gamma \in
\Gamma_i.$ If $\fX_{\Gamma_1}$ and $\fX_{\Gamma_2}$ are
length-commensurable, then $K_{\Gamma_1} = K_{\Gamma_2}.$
\end{thm}
Indeed, by Theorem \ref{T:G-3}, $\Gamma_1$ and $\Gamma_2$ are weakly
commensurable, so the assertion follows from Theorem A.

\vskip1mm

\section{Construction of nonisomorphic groups with the same tori and
noncommensurable length-commensurable locally  symmetric spaces of
type $A_n$, $D_n$ and $E_6$ .}\label{S:NI}

\vskip4mm

According to Theorem 7.3, if $K$ is a number field and $G_1$ and
$G_2$ are two $K$-forms of a connected absolutely simple group of
type different from $A_n$ $(n > 1),$ $D_{2n+1}$ and
$E_6$, then the fact that every maximal $K$-torus $T_1$ of $G_1$ is
$K$-isomorphic to some maximal $K$-torus $T_2$ of $G_2,$ and vice
versa, implies that $G_1$ and $G_2$ are $K$-isomorphic. The goal of
this section is to describe a general construction of nonisomorphic
$K$-forms of each of the types $A_n$, $D_{2n+1}$, $n>1$, and $E_6$,
which have the ``same" systems of maximal $K$-tori in a very strong
sense (see below for the definition of groups with coherently equivalent
systems of maximal $K$-tori). Furthermore, we show that
arithmetic subgroups of the forms we construct lead to
noncommensurable length-commensurable locally symmetric spaces, cf.\:Proposition 9.13. 
\vskip2mm

We begin by recalling the well-known Galois-cohomological parametrization
of the conjugacy classes of maximal $K$-tori of a given group. Let
$G$ be a connected semi-simple simply connected algebraic group over
a number field $K.$ Fix a maximal $K$-torus $T^0$ of $G,$ and let $N
= N_G(T^0)$ and $W = N/T^0$ denote its normalizer and the
corresponding Weyl group. For any field extension $\mathscr{K}/K,$
we let $\theta_{\mathscr{K}} \colon H^1({\mathscr K} , N) \to
H^1({\mathscr K} , W)$ denote the map induced by the natural homomorphism $N\to W$,  and let
$$\mathscr{C}_{\mathscr K}:= \mathrm{Ker}(H^1({\mathscr K}, N) \longrightarrow
H^1({\mathscr K} , G)).$$
The maximal $\mathscr{K}$-tori of $G$ bijectively correspond to the
$\mathscr{K}$-rational points of the variety $\mathscr{T} = G/N$ of
maximal tori of $G.$ Furthermore, $G$ acts on $\mathscr{T}$ by left
multiplication (which corresponds to the conjugation action of $G(\mathscr{K})$
on the set of maximal $\mathscr{K}$-tori), and the elements of the orbit set
$G({\mathscr K})
\backslash \mathscr{T}({\mathscr K})$ are in one-to-one
correspondence with the $G({\mathscr K})$-conjugacy classes of
maximal $\mathscr K$-tori of $G.$ The following is well-known.
\begin{lemma}\label{L:T1}
There is a natural bijection $\delta_{\mathscr K}$ from
$\mathscr{C}_{\mathscr K}$ onto $G({\mathscr K}) \backslash
\mathscr{T}({\mathscr K}).$
\end{lemma}
We just recall the construction of $\delta_{\mathscr K}.$ If $n:\,\sigma\mapsto n_{\sigma}$,
$\sigma \in \text{Gal} (\overline{{\mathscr K}}/{\mathscr{K}})$, is a
$N (\overline{\mathscr{K}})$-valued Galois cocycle representing an element of
$\mathscr{C}_{\mathscr K},$ then there exists $g \in
G(\overline{{\mathscr K}})$ such that $n_{\sigma} = g^{-1}\sigma(g)$
for all $\sigma \in \Ga(\overline{{\mathscr K}}/{\mathscr K}).$ Then
the torus $T = gT^0g^{-1}$ is defined over $\mathscr K$, and
$\delta_{\mathscr K}$ carries the cohomology class of $n$ to the
$G({\mathscr K})$-conjugacy class of $T.$

\vskip2mm

We now establish a local-global principle pertaining to the
description of maximal $K$-tori of $G.$ To formulate it, we observe
that there is an obvious map $W \longrightarrow \mathrm{Aut}\: T^0,$
so for any $x \in H^1({\mathscr K} , W),$ one can consider the
corresponding twisted $\mathscr K$-torus $_{x}T^0.$
\begin{thm}\label{T:T1}
Fix $x \in H^1(K , W)$ and  suppose that

\vskip2mm

\ {\rm (i)} $x \in \theta_{K_v}(\mathscr{C}_{K_v})$ for all $v \in
V^K;$

\vskip1mm

{\rm (ii)} \parbox[t]{11.5cm}{{\brus SH}$^2(_{x}T^0) :=
\mathrm{Ker}(H^2(K , {_{x}T^0}) \longrightarrow \prod_{v \in V^K}
H^2(K_v ,{ _{x}T^0}))$ is trivial (which holds if, for example,  
there exists $v_0 \in V^K$ such that $_{x}T^0$ is
$K_{v_0}$-anisotropic, cf.\,\cite{PlR}, Proposition 6.12).}

\vskip2mm

\noindent Then $x \in \theta_K(\mathscr{C}_K).$
\end{thm}
\begin{proof}
Applying the constructions from \cite{S}, Ch.\:I, \S 5.6, to the
exact sequence
$$
1 \to T^0 \longrightarrow N \longrightarrow W \to 1,
$$
we see that to any field extension ${\mathscr K}/K,$ one can
associate a natural cohomology class $\Delta_{\mathscr K}(x) \in
H^2({\mathscr K} , {_{x}T^0})$ such that $x \in \theta_{\mathscr
K}(H^1({\mathscr K} , N))$ if and only if $\Delta_{\mathscr K}(x)$
is trivial. It follows from (i) that $\Delta_K(x) \in ${\brus
SH}$^2(_{x}T^0),$ which is trivial by (ii). Thus, $x = \theta_K(y)$
for some $y \in H^1(K , N).$ Furthermore, according to {\it
loc.cit.,} \S 5.5, for any ${\mathscr K}/K$ there is a natural
surjective map $\nu_{\mathscr K} \colon H^1({\mathscr K} ,
{_{x}T^0}) \to \theta_{\mathscr K}^{-1}(x).$ For each $v \in
V^K_{\infty},$ by (i), we can find $z_v \in \mathscr{C}_{K_v}$ such
that $\theta_{K_v}(z_v) = x,$ and then pick $t_v \in H^1(K_v ,
{_{x}T^0})$ for which $\nu_{K_v}(t_v) = z_v.$ By \cite{PlR},
Proposition 6.17, the diagonal map $H^1(K , {_{x}T^0})
\longrightarrow \prod_{v \in V^K_{\infty}} H^1(K_v ,{ _{x}T^0})$ is
surjective, so there is $t \in H^1(K , {_{x}T^0})$ that maps to
$(t_v)_{v \in V^K_{\infty}}.$ Set $z = \nu_K(t).$ Then $z$  maps onto
$(z_v)_{v \in V^K_{\infty}}$ under the diagonal map $H^1(K , N)
\longrightarrow \prod_{v \in V^K_{\infty}} H^1(K_v , N).$ Combining
the fact that $z_v \in \mathscr{C}_{K_v}$ with the injectivity of
the map $H^1(K , G) \longrightarrow \prod_{v \in V^K_{\infty}}
H^1(K_v , G)$ (\cite{PlR}, Theorem 6.6), we obtain that $z \in
\mathscr{C}_K.$ Thus, $x = \theta_K(z) \in \theta_K(\mathscr{C}_K),$
as required.
\end{proof}

\vskip2mm

We now turn to the comparison of the sets of maximal $K$-tori of two
absolutely simple simply connected $K$-groups $G_1$ and $G_2.$ We
assume that there exist maximal $K$-tori $T_1^0$ of $G_1$ and
$T_2^0$ of $G_2$, and a $\overline{K}$-isomorphism $\varphi_0 \colon
G_1 \to G_2$ whose restriction to $T_1^0$ is an isomorphism onto
$T_2^0$ defined over $K,$ and we fix these $T_1^0,$ $T_2^0$ and
$\varphi_0$ for the rest of the section. Clearly, $\varphi_0$
induces an isomorphism between $N_1 = N_{G_1}(T_1^0)$ and $N_2 =
N_{G_2}(T_2^0),$ and hence an isomorphism $\varphi_0^W$ between the
Weyl groups $W_1 = N_1/T_1^0$ and $W_2 = N_2/T_2^0.$
\begin{lemma}\label{L:T2}
The map $\varphi_0^W:\, W_1\rightarrow W_2$ is defined over $K$.
\end{lemma}
\begin{proof}
Since $\varphi_0 \vert T_1^0$ is defined over $K$, for any $n \in
N_1(\overline{K}),$ $t \in T_1^0(\overline{K})$ and any $\sigma \in
\Ga(\overline{K}/K),$ we have
$$
\varphi_0(\sigma(ntn^{-1})) = \sigma(\varphi_0(ntn^{-1})),
$$
which implies that
$$
\varphi_0(\sigma(n)) \varphi_0(\sigma(t)) \varphi_0(\sigma(n))^{-1}
= \sigma(\varphi_0(n))
\sigma(\varphi_0(t))\sigma(\varphi_0(n))^{-1}.
$$
Since $\varphi_0(\sigma(t)) = \sigma(\varphi_0(t)),$ we conclude
that $\sigma(\varphi_0(n))\equiv \varphi_0(\sigma(n))$ modulo
$T^0_2({\overline{K}})$. This means that $\varphi_0^W$ commutes with
every $\sigma \in \Ga(\overline{K}/K),$ hence it is defined over
$K$.
\end{proof}

\vskip2mm

Lemma \ref{L:T2} enables us to define, for any field extension
${\mathscr K}/K,$ the induced isomorphism $H^1({\mathscr K},
W_1)\rightarrow H^1({\mathscr K}, W_2),$ which will also be denoted
by $\varphi_0^W.$ This isomorphism will play a critical role in
comparing the maximal $K$-tori of $G_1$ and $G_2.$ More precisely,
for $i = 1 , 2,$ we let
$\theta_{\mathscr K}^{(i)} \colon H^1({\mathscr K} , N_i) \to
H^1({\mathscr K} , W_i)$ be the map induced by the canonical
homomorphism $N_i \to W_i.$ Furthermore, let
$\mathscr{C}^{(i)}_{\mathscr K} = \mathrm{Ker}(H^1({\mathscr K} ,
N_i) \to H^1({\mathscr K} , G_i)),$ and let $\delta^{(i)}_{\mathscr
K} \colon \mathscr{C}^{(i)}_{\mathscr{K}} \to G_i({\mathscr K})
\backslash \mathscr{T}_i({\mathscr K})$ (where $\mathscr{T}_i$ is
the variety of maximal tori of $G_i$) be the bijection provided by
Lemma \ref{L:T1}. Then the condition that $G_1$ and $G_2$ have the
``same"  maximal $K$-tori is basically equivalent to the following
\begin{equation}\label{T:E125}
\varphi_0^W(\theta^{(1)}_K(\mathscr{C}^{(1)}_K)) =
\theta^{(2)}_K(\mathscr{C}^{(2)}_K).
\end{equation}
To give a precise interpretation of (\ref{T:E125}), we need to introduce
the following definition.

\vskip2mm

\vskip2mm \noindent {\bf Definition.} Let $\mathscr K$ be a field
extension of $K$ and let $T_1$ be a maximal $\mathscr K$-torus of
$G_1.$ A $\mathscr K$-embedding $\iota \colon T_1 \to G_2$ will be
called {\it coherent} (relative to $\varphi_0$) if there exists a
$\overline{\mathscr K}$-isomorphism $\varphi \colon G_1 \to G_2$ of
the form $\varphi = \mathrm{Int}\: h \circ \varphi_0,$ where $h \in
G_2(\overline{\mathscr K}),$ such that $\iota = \varphi \vert T_1.$
Furthermore, we say that $G_1$ and $G_2$ {\it have coherently
equivalent systems of maximal $K$-tori} if every maximal $K$-torus
$T_1$ of $G_1$ admits a coherent $K$-embedding into $G_2,$ and every
maximal $K$-torus $T_2$ of $G_2$ admits a coherent $K$-embedding
into $G_1.$

\vskip2mm

\begin{lemma}\label{L:T4}
Let $T_1$ be a maximal $\mathscr K$-torus of $G_1,$ and let $x_1 \in
\mathscr{C}^{(1)}_{\mathscr K}$ be the cohomology class that
corresponds to $T_1$ under $\delta^{(1)}_{\mathscr K}.$ Then $T_1$
admits a coherent (relative to $\varphi_0$) $\mathscr K$-embedding
into $G_2$ if and only if $\varphi_0^W(\theta^{(1)}_{\mathscr
K}(x_1)) \in \theta^{(2)}_{\mathscr K}(\mathscr{C}^{(2)}_{\mathscr
K}).$ Thus, {\rm (\ref{T:E125})} is equivalent to the condition that
$G_1$ and $G_2$ have coherently equivalent systems of maximal
$K$-tori.
\end{lemma}
\begin{proof}
Pick $g_1 \in G_1(\overline{\mathscr K})$ so that $T_1 =
g_1T_1^0g_1^{-1}.$ Then $x_1$ is represented by the $N_1(\overline{\mathscr K})$-valued
Galois cocycle $\sigma \mapsto \alpha_{\sigma} := g_1^{-1}\sigma(g_1)$, $\sigma \in
\Ga(\overline{\mathscr K}/{\mathscr K}),$ and therefore,
$\varphi_0^W(\theta^{(1)}_{\mathscr K}(x_1))$ is represented by the cocycle 
\begin{equation}\label{E:T5}
\sigma\mapsto \beta_{\sigma} := \varphi_0(g_1^{-1}\sigma(g_1))T_2^0 \in W_2.
\end{equation}
Let $\varphi \colon G_1 \to G_2$ be an isomorphism of the form
$\varphi = \mathrm{Int}\: h \circ \varphi_0,$ where $h \in
G_2(\overline{\mathscr K}).$ Then $T_2 := \varphi(T_1)$ can be
written in the form $T_2 = g_2T_2^0g_2^{-1},$ where $g_2 =
h\varphi_0(g_1).$ So, $T_2$ is defined over $\mathscr K$ if and only
if $g_2^{-1}\sigma(g_2) \in N_2 ({\overline{\mathscr K}})$ for all $\sigma \in
\Ga(\overline{\mathscr K}/{\mathscr K}),$ in which case the class
$x_2$ corresponding to $T_2$ is represented by the $N_2({\overline{\mathscr K}})$-valued Galois cocycle $\sigma\mapsto g_2^{-1}\sigma(g_2).$ Then
$\theta^{(2)}_{\mathscr K}(x_2)$ is represented by the cocycle
\begin{equation}\label{E:T6}
\sigma\mapsto \gamma_{\sigma} := g_2^{-1}\sigma(g_2)T_2^0 =
\varphi_0(g_1)^{-1}h^{-1}\sigma(h)\sigma(\varphi_0(g_1)) T_2^0 \in
W_2.
\end{equation}
Finally, notice that the condition that $\varphi \vert T_1$ is defined
over $\mathscr K$ is equivalent to
\begin{equation}\label{E:T7}
\varphi(\sigma(g_1tg_1^{-1})) = \sigma(\varphi(g_1tg_1^{-1})) \ \
\text{for all} \ \ t \in T^0(\overline{\mathscr K}) \  \ \text{and}
\ \ \sigma \in \Ga(\overline{\mathscr K}/{\mathscr K}).
\end{equation}
The left- and right-hand sides of (\ref{E:T7}) can be expanded as
follows:
$$
\varphi(\sigma(g_1tg_1^{-1})) = h \varphi_0(\sigma(g_1tg_1^{-1}))
h^{-1} = h\varphi_0(\sigma(g_1))
\varphi_0(\sigma(t))\varphi_0(\sigma(g_1))^{-1}h^{-1}
$$
and
$$
\sigma(\varphi(g_1tg_1^{-1})) =
\sigma(h\varphi_0(g_1tg_1^{-1})h^{-1}) =
\sigma(h)\sigma(\varphi_0(g_1)) \sigma(\varphi_0(t))
\sigma(\varphi_0(g_1))^{-1} \sigma(h)^{-1}.
$$
So, since $\varphi_0(\sigma(t)) = \sigma(\varphi_0(t)),$ we see that
(\ref{E:T7}) is equivalent to
\begin{equation}\label{E:T11}
\varphi_0(\sigma(g_1))^{-1}h^{-1}\sigma(h)\sigma(\varphi_0(g_1)) \in
T_2^0 \ \ \text{for all} \ \ \sigma \in \Ga(\overline{\mathscr
K}/{\mathscr K}).
\end{equation}

\vskip1mm

Now, suppose $\varphi \vert T_1$ is defined over $\mathscr K$,
i.e.,\,(\ref{E:T11}) holds. We claim that
$\varphi_0^W(\theta^{(1)}_{\mathscr K}(x_1)) =
\theta^{(2)}_{\mathscr K}(x_2) \in \theta^{(2)}_{\mathscr
K}(\mathscr{C}^{(2)}_{\mathscr K}).$ Indeed, combining (\ref{E:T11})
with (\ref{E:T6}) and (\ref{E:T5}), we see that
$$
\gamma_{\sigma} =
\varphi_0(g_1)^{-1}h^{-1}\sigma(h)\sigma(\varphi_0(g_1))T_2^0 =
\varphi_0(g_1^{-1}\sigma(g_1))T_2^0 = \beta_{\sigma},
$$
as required.

\vskip3mm

Conversely, suppose $\varphi_0^W(\theta^{(1)}_{\mathscr K}(x_1)) \in
\theta^{(2)}_{\mathscr K}(\mathscr{C}^{(2)}_{\mathscr K}).$ This
means that there exists $g_2 \in G_2(\overline{\mathscr K})$ such
that
\begin{equation}\label{E:T12}
\beta_{\sigma} = g_2^{-1}\sigma(g_2)T_2^0 \ \ \text{for all} \ \
\sigma \in \Ga (\overline{\mathscr K}/{\mathscr K}).
\end{equation}
Set $h = g_2\varphi_0(g_1)^{-1}$ and $\varphi = \mathrm{Int}\: h
\circ \varphi_0.$ We need to show that $\varphi \vert T_1$ is
defined over $\mathscr K$, in other words, (\ref{E:T11}) holds. But
this is obtained directly by combining (\ref{E:T5}) with
(\ref{E:T12}).
\end{proof}

\vskip1mm

Combining  Theorem \ref{T:T1} with Lemma \ref{L:T4}, we obtain the
following local-global principle for the existence of a coherent
$K$-embedding of a $K$-torus as a maximal torus in a semi-simple group.
\begin{thm}\label{T:T2}
Let $G_1$ and $G_2$ be two connected semi-simple simply connected
algebraic groups over a number field $K.$ Assume that

\vskip1mm

\noindent {$(*)$} \parbox[t]{12cm}{there exist maximal $K$-tori
$T_1^0$ of $G_1$ and $T_2^0$ of $G_2,$ and a
$\overline{K}$-isomorphism $\varphi_0 \colon G_1 \to G_2$ whose
restriction to $T_1^0$ is an isomorphism onto $T_2^0$ defined over
$K$.}

\vskip1mm

\noindent Let $T_1$ be a maximal $K$-torus of $G_1$ such that {\brus
SH}$^2(T_1)$ is trivial (which automatically holds if there exists
$v_0 \in V^K$ such that $T_1$ is $K_{v_0}$-anisotropic). If $T_1$
admits a coherent (relative to $\varphi_0$) $K_v$-embedding into
$G_2$ for every $v \in V^K$, then it admits a coherent $K$-embedding
into $G_2.$
\end{thm}
\vskip2mm

The following lemma explains why coherent embeddings of tori are easier to
analyze if the ambient group is not of type  $D_{2n}.$
\begin{lemma}\label{L:T3}
Assume that $G_1$ and $G_2$ are of type different from $D_{2n},$ and
let ${\mathscr K}/K$ be a field extension. If \,$T_1$ is a maximal
$\mathscr K$-torus of $G_1$ and $\varphi \colon G_1 \to G_2$ is a
$\overline{{\mathscr K}}$-isomorphism such that $\iota := \varphi
\vert T_1$ is defined over $\mathscr K$, then either $\iota$, or
$\iota',$ defined by $\iota'(t) = \iota(t)^{-1},$ is a {\rm
coherent} $\mathscr K$-embedding of $T_1$ into $G_2$ (in particular,
$T_1$ admits such an embedding). Thus, if $G_1$ and $G_2$ are
$\mathscr K$-isomorphic, then they have coherently equivalent
systems of maximal $\mathscr K$-tori.
\end{lemma}
\begin{proof}
Obviously, $T_2 := \varphi(T_1)$ is defined over $\mathscr K$. Let
$\Phi_2$ be the root system of $G_2$ with respect to $T_2.$ Since
$G_2$ is not of type $D_{2n},$ the quotient
$\mathrm{Aut}(\Phi_2)/W(\Phi_2)$ is of order $\leqslant 2,$ and in
case it is of  order 2, the automorphism $\alpha \mapsto -\alpha$
represents the nontrivial coset. Equivalently,
$\mathrm{Aut}\,G_2/\mathrm{Int}\,G_2$ has order $\leqslant 2,$ and
in case it has order 2, there is an outer automorphism $\tau$ of
$G_2$ defined over $\overline{\mathscr K}$ such that $\tau(t) =
t^{-1}$ for all $t \in T_2.$ Set $\varphi' = \tau \circ \varphi$,
then $\varphi' \vert T_1 = \iota'$. Since one of $\varphi$ and
$\varphi'$ is of the form $\mathrm{Int}\: h \circ \varphi_0,$ the
lemma follows.
\end{proof}

\vskip0.5mm

Combined with Theorem \ref{T:T2}, this lemma yields the following.
\begin{cor}\label{C:T10}
Let $G_1$ and $G_2$ be two connected absolutely simple simply connected
algebraic groups of type different from $D_{2n}$, and suppose that
the condition {$(*)$} of Theorem \ref{T:T2} holds. Assume in
addition that {\brus SH}$^2$ is trivial for all maximal $K$-tori of
$G_1$ and $G_2$ (which automatically holds if there exists a place
$v_0$ of $K$ such that $G_i$ is $K_{v_0}$-anisotropic
for $i = 1 , 2$). If $G_1 \simeq G_2$ over $K_v$, for all $v \in
V^K$, then $G_1$ and $G_2$ have coherently equivalent systems of
maximal $K$-tori.
\end{cor}

Of course, if $G_1$ and $G_2$ are not of type $A$, $D$
or $E_6,$ then the assumption that $G_1 \simeq G_2$ over $K_v$ for
all $v \in V^K$ implies that $G_1 \simeq G_2$ over $K,$ and our
assertion becomes obvious (cf.\,Lemma \ref{L:T3}). We will use
Corollary \ref{C:T10} to show that for each of the types $A_n,$
$D_{2n+1},$ or $E_6$, one can construct an arbitrarily large number
of pairwise nonisomorphic absolutely simple simply connected
$K$-groups of this type with coherently equivalent systems of
maximal $K$-tori (cf.\,Theorem \ref{T:T5}).
\vskip1mm

Let $G_0$ be a connected absolutely simple simply connected quasi-split
$K$-group of one of the following types: $A_n$ $(n > 1),$ $D_{2n+1}$
and $E_6.$ We first describe a general construction of nonisomorphic
inner twists $G_1$ and $G_2$ of $G_0$ which are isomorphic over
$K_v$ for all $v \in V^K.$ Let $L$ be the minimal Galois extension
of $K$ over which $G_0$ splits, and let $V_0$ be the set of $v \in
V^K_f$ that split in $L.$ We let $C$ denote the center of $G_0;$
clearly, $C$ is $L$-isomorphic to $\mu_{\ell},$ the group of
$\ell$-th roots of unity, where $\ell = n + 1$ for $G_0$ of type
$A_n,$ $\ell = 4$ for type $D_{2n+1},$ and $\ell = 3$ for type
$E_6.$ Each $x \in G_0$ gives the inner automorphism $z \mapsto
xzx^{-1}$ of $G_0$.  This leads to the natural isomorphism $i$ from
the adjoint group $\overline{G}_0$ of $G_0$  onto the group of inner
automorphisms $\mathrm{Int}\: G_0$ ($\subset\mathrm{Aut}\:G_0$). Any automorphism $g$ of $G_0$
can be regarded as an automorphism of $\overline{G}_0,$ and then for
every $x \in \overline{G}_0,$ we have $g \circ i(x) \circ g^{-1} =
i(g(x))$ in $\mathrm{Aut}\: G_0.$ \vskip1mm

For a class $\mathfrak{c}\in H^1(K,{\overline G}_0)$, in the sequel we
will let $\sigma\mapsto \mathfrak{c}_{\sigma}$, $\sigma\in \text{Gal}({\overline K}/K)$,
denote a Galois cocycle representing $\mathfrak{c}$.
\vskip1mm

For any $v \in V^K,$ we have the following commutative diagram
$$
\begin{array}{ccc}
H^1(K , \overline{G}_0) & \stackrel{\alpha}{\longrightarrow} & H^1(K
, \mathrm{Aut}\: G_0) \\
\gamma_v \downarrow & & \downarrow \beta_v \\
H^1(K_v , \overline{G}_0) & \stackrel{\alpha_v}{\longrightarrow} &
H^1(K_v , \mathrm{Aut}\: G_0),
\end{array}
$$
in which $\alpha$ and $\alpha_v$ are induced by $i.$ Furthermore,
for any extension ${\mathscr K}/K$ there is a natural map
$\rho_{\mathscr K} \colon H^1({\mathscr K} , \overline{G}_0)
\rightarrow H^2({\mathscr K} , C).$ We will also need the map $\mu
\colon H^2(K , C) \to \bigoplus_v H^2(K_v , C).$


\vskip1mm

\begin{lemma}\label{L:T5}
Let $\xi_1,\,\, \xi_2 \in H^1(K , \overline{G}_0).$

\vskip2mm

\ {\rm (i)} If $\rho_K(\xi_1) \neq \pm \rho_K(\xi_2)$, then
$\alpha(\xi_1) \neq \alpha(\xi_2).$

\vskip1mm

{\rm (ii)} \parbox[t]{11.5cm}{If $v \in V^K_f$ and
$\rho_{K_v}(\gamma_v(\xi_1)) = \pm \rho_{K_v}(\gamma_v(\xi_2))$, then
$\beta_v(\alpha(\xi_1)) = \beta_v(\alpha(\xi_2)).$}
\end{lemma}
\begin{proof}
Notice that $\mathrm{Aut}\: G_0$ has the following semi-direct product
decomposition
$$
\mathrm{Aut}\: G_0 = \mathrm{Int}\: G_0 \rtimes \Sigma,
$$
where $\Sigma$ is a $K$-subgroup of order two, whose nontrivial element $s$
is defined over $K$ and acts on $C$ as $c \mapsto c^{-1}.$

(i): Suppose $\alpha(\xi_1) = \alpha(\xi_2).$ Then there exists $g
\in \mathrm{Aut}\: G_0$ such that
$$i({\xi_{2}}_{\sigma}) = g\circ i({\xi_{1}}_{\sigma}) \circ
\sigma(g)^{-1} \ \ \text{for all} \ \ \sigma \in \Ga(\overline{K}/K).
$$
If $g \in {\rm {Int}}\,G_0,$ then $\xi_1 = \xi_2$, and therefore,
$\rho_K(\xi_1) = \rho_K(\xi_2).$ Now, suppose $g \notin {\rm Int
}\,G_0.$ Then $g = hs,$ $h \in {\rm Int}\, G_0.$ The cohomology
class $\xi'_2$ in $H^1(K , \overline{G}_0)$ corresponding to the cocycle
$$\sigma\mapsto {\xi'_{2}}_{\sigma} = s({\xi_{1}}_{\sigma}),\ \ \ \sigma\in \text{Gal}
({\overline K}/K),$$ clearly equals $\xi_2$.
As $s(c) = c^{-1}$ for $c\in C$, we conclude that
$$
\rho_K(\xi_2) = \rho_K(\xi'_2) = - \rho_K(\xi_1),
$$
a contradiction.

\vskip2mm

(ii): Recall that $\rho_{K_v}$ is a bijection for any $v \in V^K_f$
(cf.\,\cite{PlR}, Corollary of Theorem 6.20), so our claim is
obvious if
$\rho_{K_v}(\gamma_v(\xi_1)) = \rho_{K_v}(\gamma_v(\xi_2)).$ Suppose
now that $\rho_{K_v}(\gamma_v(\xi_1)) =
-\rho_{K_v}(\gamma_v(\xi_2)).$ Consider the ${\overline
G}({\overline K})$-valued Galois cocycle $\sigma \mapsto {\xi'_{2}}_{\sigma}:= s({\xi_{2}}_{\sigma})$, and let $\xi'_2$ be the associated
cohomology class. Then for $\sigma\in \text{Gal}({\overline K}/K)$
we have
$$i({\xi'_{2}}_ {\sigma}) = s\circ i({\xi_{2}}_{\sigma}) \circ
s^{-1} = s\circ i({\xi_{2}}_{\sigma}) \circ \sigma(s)^{-1},
$$
so $\alpha(\xi'_2) = \alpha(\xi_2).$ On the other hand,
$$\rho_{K_v}(\gamma_v(\xi'_2)) = - \rho_{K_v}(\gamma_v(\xi_2)) =
\rho_{K_v}(\gamma_v(\xi_1)).$$ Then $\gamma_v(\xi'_2) =
\gamma_v(\xi_1),$ and
$$
\beta_v (\alpha(\xi_1)) = \beta_v (\alpha(\xi'_2)) =
\beta_v (\alpha(\xi_2)).$$\end{proof}

\vskip2mm

Let $\widehat{C}$ be the character group of $C.$ Fix a generator
$\chi$ of  $\widehat{C}(K),$ and let $d$ denote its order. For each
$v \in V^K,$ $\chi$ induces a character $$\chi_v \colon H^2(K_v , C)
\to H^2(K_v, {\rm{GL}}_1)\subset \Q/\Z.$$ If $v \in V_0$, then
$H^2(K_v , C) \simeq \mathrm{Br}(K_v)_{\ell}$ is cyclic of order
$\ell,$ and one can choose a generator $b_v \in H^2(K_v , C)$ such
that $\chi_v(b_v) = 1/d.$ Now, let $V$ be a finite
subset of $V^K$ containing $V^K_{\infty}$,  and suppose that for
each $v \in V$ we are given $\xi^{(v)} \in H^1(K_v ,
\overline{G}_0).$ Fix an integer $t \geqslant 1,$ and pick $2(t+1)$
places
$$
v'_0 ,\, v''_0,\, v'_1,\, v''_1,\, \ldots ,\, v'_t,\, v''_t \in V_0
\setminus (V_0 \cap V).$$
Let $V_t =\{ v'_0 ,\, v''_0,\, v'_1,\, v''_1,\, \ldots ,\, v'_t,\, v''_t \}$. Now pick $x_{v''_0} \in H^2(K_{v''_0} , C)$ so that
$$
\sum_{v \in V} \chi_v(\rho_{K_v}(\xi^{(v)})) + \chi_{v'_0}(b_{v'_0})
+ \chi_{v''_0}(x_{v''_0}) = 0.
$$
Next, fix $\displaystyle \varepsilon = (\varepsilon_1, \ldots ,
\varepsilon_t) \in E_t := \prod_{i = 1}^t \{ \pm 1\},$ and  consider
$(x(\varepsilon)_v) \in \bigoplus_v H^2(K_v , C)$ with the following
components:
\begin{equation}\label{E:T25}
x(\varepsilon)_v = \left\{ \begin{array}{ccl} \rho_{K_v}(\xi^{(v)})
& , & v \in V \\ b_{v'_0} & , & v = v'_0 \\ x_{v''_0} & , & v = v''_0 \\
\varepsilon_j b_{v'_j} & , &
v = v'_j, \   j \geqslant 1 \\
-\varepsilon_j b_{v''_j}  & , & v = v''_j, \ j \geqslant 1 \\
0 &  , &  \text{for all other}\ v
\end{array} \right.
\end{equation}
We obviously have $\sum_v \chi_v(x(\varepsilon)_v) = 0,$ so it
follows from a theorem of Poitou-Tate (cf.\:\cite{S}, Ch.\,II, \S 6, Theorem C)
that there exists $x(\varepsilon) \in H^2(K , C)$ such that
$\mu(x(\varepsilon)) = (x(\varepsilon)_v).$ We now want to construct
a maximal $K$-torus $\overline{T}_0$ of $\overline{G}_0$ (depending
on $V,$ $\xi^{(v)}$ for $v \in V$, and $V_t$) such that for each $\varepsilon \in E_t$,
$x(\varepsilon)$ lifts to a class $\zeta(\varepsilon) \in H^1(K ,
\overline{T}_0)$ whose image in $H^1(K_v, {\overline G}_0)$ is $\xi^{(v)}$ for all $v\in V$.
\vskip1mm

For every real $v,$ $\xi^{(v)}$ is given by an element $g_v \in
\overline{G}_0({\overline K}_v)$ such that $g_v \overline{g}_v = 1,$
where $\overline{g}_v$ denotes the conjugate of $g_v$ under the
nontrivial automorphism of $\overline{K}_v/K_v = \C/\R.$ It follows
from the uniqueness of the Jordan decomposition that the semi-simple
and the unipotent components $g_v^s$, $g_v^u$ of $g_v$ also define
cocycles. If $g_v^u \neq 1,$ then the 1-dimensional connected unipotent
subgroup $U$ generated by $g_v^u$ is defined over $K_v = \R.$ Using
the fact that $H^1(K_v , U)$ is trivial, one sees that $\xi^{(v)}$
is the cohomology class given by $g_v^s$. So we can assume that
$g_v$ is semi-simple. Then $g_v$ is contained in the connected
centralizer $H := Z_{\overline{G}_0}(g_v)^{\circ}$ (cf.\,\cite{Bo},
Corollary 11.12), and $H$ is defined over $K_v.$ Hence, $g_v$ is
contained in a maximal $K_v$-torus $\overline{T}^{(v)}$ of $H$ which
is also a maximal torus of $\overline{G}_0.$ For each $v \in (V
\setminus V^K_{\infty}) \cup V_t,$ we pick a maximal $K_v$-torus
$\overline{T}^{(v)}$ of $\overline{G}_0$ which is anisotropic over
$K_v$ (see \cite{PlR}, Theorem 6.21, or \cite{B}, \S 2.4). Using the
weak approximation property for the variety of maximal tori of
$\overline{G}_0$ (cf.\,\cite{PlR}, Corollary 3 in \S 7.1), we can
find a maximal $K$-torus $\overline{T}_0$ of $\overline{G}_0$ which
is conjugate to $\overline{T}^{(v)}$ under an element of
$\overline{G}_0(K_v)$ for all $v \in V \cup V_t.$ Let $\pi \colon
G_0 \to \overline{G}_0$ be the natural $K$-isogeny, and $T_0 =
\pi^{-1}(\overline{T}_0).$

\begin{lemma}\label{L:T50}
For every $\varepsilon \in E_t$, there exists $\zeta(\varepsilon) \in
H^1(K , \overline{T}_0)$ which maps onto $x(\varepsilon)$ under the
coboundary map $H^1(K , \overline{T}_0) \to H^2(K , C)$, and whose
image in $H^1(K_v , \overline{G}_0)$ equals
$\xi^{(v)}$ for all $v \in V.$ \end{lemma}
\begin{proof}
For any real $v$, as ${\overline T}_0$ is conjugate to ${\overline T}^{(v)}$ under an element of
${\overline G}_0(K_v)$, and $\xi^{(v)}$ is given
by $g_v\in {\overline T}^{(v)}({\overline K}_v)$,
there exists a cohomology class ${\xi'}^{(v)}$ in
$H^1(K_v, {\overline T}_0)$ which maps onto $\xi^{(v)}$
under the natural map $H^1(K_v, {\overline T}_0)\rightarrow H^1(K_v,{\overline G}_0)$.
On the other hand,  for every nonarchimedean $v\in V$,
as ${\overline T}_0$ is anisotropic over $K_v$, the natural
map $H^1(K_v, {\overline T}_0)\rightarrow H^1(K_v, {\overline G}_0)$ is
onto (see the proof of Theorem 6.20 on p.\,326 of \cite{PlR}),
there is a ${\xi'}^{(v)}\in H^1(K_v, {\overline T}_0)$ which maps onto $\xi^{(v)}$.
\vskip1mm

We have the following commutative diagram with exact rows:
$$
\begin{array}{ccccc}
H^1(K , \overline{T}_0) & \stackrel{\delta_1}{\longrightarrow} &
H^2(K , C) & \stackrel{\delta_2}{\longrightarrow} & H^2(K , T_0) \\
\eta_1 \downarrow & & \eta_2 \downarrow & & \eta_3\downarrow \\
\bigoplus_{v} H^1(K_v , \overline{T_0}) &
\stackrel{\Delta_1}{\longrightarrow} & \bigoplus_v H^2(K_v , C) &
\stackrel{\Delta_2}{\longrightarrow} & \bigoplus_v H^2(K , T_0)
\end{array}
$$
(notice that $\eta_2$ actually coincides with $\mu$). First, we will
show that $x(\varepsilon) \in \mathrm{Im}\: \delta_1 =
\mathrm{Ker}\: \delta_2.$ Observe that
\begin{equation}\label{E:T10}
x(\varepsilon)_v \in \mathrm{Im}(H^1(K_v , \overline{T}_0) \to
H^2(K_v , C))
\end{equation}
for all $v.$ This is obvious if $v \not\in V \cup
V_t.$ For any real $v $, this follows from the fact that
$x(\varepsilon)_v = \rho_{K_v}(\xi^{(v)})$, and $\xi^{(v)}$ is
the image of ${\xi'}^{(v)}\in H^1(K_v
, \overline{T}_0)$. For a nonarchimedean $v \in V \cup
V_t,$ by our construction $T_0$ is $K_v$-anisotropic, and it follows
from the Nakayama-Tate Theorem (cf.\,\cite{PlR}, Theorem 6.2) that
$H^2(K_v , T_0)$ is trivial. So the map $H^1(K_v , \overline{T}_0)
\to H^2(K_v , C)$ is surjective, and (\ref{E:T10}) is automatic.
Thus, $\eta_2(x(\varepsilon)) = (x(\varepsilon)_v) \in \mathrm{Im}\
\Delta_1,$ so
$$
\Delta_2(\eta_2(x(\varepsilon))) = \eta_3(\delta_2(x(\varepsilon)))
= 0.
$$
Since $T_0$ is anisotropic at every $v \in V_t,$ we have that
{\brus SH}$^2(T_0) = \mathrm{Ker}\: \eta_3$ is trivial, and hence
$\delta_2(x(\varepsilon)) = 0,$ as required. Fix
$\zeta'(\varepsilon) \in H^1(K , \overline{T}_0)$ such that
$\delta_1(\zeta'(\varepsilon)) = x(\varepsilon).$

For an extension $\mathcal{K}/K,$ we consider the natural
homomorphism $$\lambda_{\mathcal{K}} \colon H^1(\mathcal{K} , T_0)
\to H^1(\mathcal{K} , \overline{T}_0),$$ and for $v \in V^K$, we let
$\zeta'(\varepsilon)^{(v)}$ denote the image of
$\zeta'(\varepsilon)$ under the restriction map $H^1(K ,
\overline{T}_0) \to H^1(K_v , \overline{T}_0).$ For each $v \in V,$
the cohomology classes $\zeta'(\varepsilon)^{(v)}$ and ${\xi'}^{(v)}$ have the
same image in $H^2(K_v , C),$ so there exists $\theta(\varepsilon)_v
\in H^1(K_v , T_0)$ such that
$$
{\xi'}^{(v)} = \lambda_{K_v}(\theta(\varepsilon)_v) \cdot
\zeta'(\varepsilon)^{(v)}.
$$
By (\cite{PlR}, Proposition 6.17), the map $H^1(K , T_0) \to
\prod_{v \in V^K_{\infty}} H^1(K_v , T_0)$ is surjective. Pick
$\theta(\varepsilon) \in H^1(K , T_0)$ which maps onto
$(\theta(\varepsilon)_v)_{v \in V^K_{\infty}},$ and set
$\zeta(\varepsilon) = \lambda_K(\theta(\varepsilon)) \cdot
\zeta'(\varepsilon).$ Let $\zeta(\varepsilon)^{(v)}$ be
the image of $\zeta(\varepsilon)$ under
the map $H^1(K,{\overline{T}}_0)\to H^1(K_v, {\overline{T}}_0)$. Then $\delta_1(\zeta(\varepsilon)) =
\delta_1(\zeta'(\varepsilon)) = x(\varepsilon)$ and
$\zeta(\varepsilon)^{(v)} = {\xi'}^{(v)}$ for all $v \in
V^K_{\infty}.$ Finally, to show that the image of
$\zeta(\varepsilon)^{(v)}$ in $H^1(K_v , \overline{G}_0)$ coincides
with $\xi^{(v)}$ for nonarchimedean $v \in V,$ we observe that these
elements have the same image under $\rho_{K_v},$ which is a
bijection for all $v\in V^K_f$ (Corollary in \S 6.4 of \cite{PlR}).
\end{proof}

 Let $\zeta(\varepsilon)$ be as in the preceding lemma, and $\xi(\varepsilon)$ be the image of $\zeta(\varepsilon)$ under the
 natural map $H^1(K , \overline{T}_0) \to H^1(K ,
 \overline{G}_0)$. Then $\rho_K(\xi(\varepsilon)) = x(\varepsilon)$ and
$\gamma_v(\xi(\varepsilon)) = \xi^{(v)}$ for all $v \in V.$ Fix two
distinct $\varepsilon_1 , \varepsilon_2 \in E_t,$ and let $\xi_j =
\xi(\varepsilon_j).$ Since each $b_v$ has order $\ell > 2,$ it
follows from (\ref{E:T25}) that $\mu(\rho_K(\xi_1)) \neq \pm
\mu(\rho_K(\xi_2)),$ hence $\rho_K(\xi_1) \neq \pm \rho_K(\xi_2),$
so according to Lemma \ref{L:T5}(i), $\alpha(\xi_1) \neq
\alpha(\xi_2).$ On the other hand, we have
$$
\rho_{K_v}(\gamma_v(\xi_1)) = 0 = \rho_K(\gamma_v(\xi_2))  \ \
\text{for any} \ \ v \in V^K  \setminus (V \cup V_0),
$$
$$
\rho_{K_v}(\gamma_v(\xi_1)) = \pm \rho_{K_v}(\gamma_v(\xi_2)) \ \
\text{for any} \ \ v \in V_0,
$$
and
$$
\gamma_v(\xi_1) = \xi^{(v)} = \gamma_v(\xi_2) \ \ \text{for any} \ \
v \in V.
$$
Using Lemma \ref{L:T5}(ii), we now see that $\beta_v(\alpha(\xi_1))
= \beta_v(\alpha(\xi_2))$ for all $v \in V^K.$ Thus, we obtain the
following proposition.
\begin{prop}\label{P:T2}
The $2^t$ elements $\xi(\varepsilon)\in H^1(K,{\overline G}_0)$, $\varepsilon \in E_t,$ have
the following properties: the elements $\alpha(\xi(\varepsilon)) \in
H^1(K , \mathrm{Aut}\: G_0)$ are pairwise distinct, while for any $v
\in V^K,$ the elements $\beta_v(\alpha(\xi(\varepsilon))) \in
H^1(K_v , \mathrm{Aut}\: G_0)$ are all equal, and, in addition,
$\gamma_v(\xi(\varepsilon)) = \xi^{(v)}$ for all $v \in V.$
\end{prop}

For $\xi(\varepsilon)$ as above,  we let $G_{\varepsilon}$ denote the
form of $G_0$ obtained by twisting it by a cocycle representing $\alpha(\xi(\varepsilon)).$
Since the cohomology classes  $\alpha(\xi(\varepsilon)),$ $\varepsilon \in
E_t,$ are pairwise distinct, the corresponding groups
$G_{\varepsilon}$ are pairwise nonisomorphic over $K.$ Now, fix
$\varepsilon_1 , \varepsilon_2 \in E_t,$ and set
$$
\zeta_j = \zeta(\varepsilon_j) \in H^1(K , \overline{T}_0), \ \
\xi_j = \xi(\varepsilon_j) \in H^1(K , \overline{G}_0) \ \
\text{and} \ \ G_j = G_{\varepsilon_{j}}
$$
for $j = 1 , 2.$  As $\xi_j$ is the image of $\zeta_j$ under the
natural map $H^1(K,{\overline T}_0)\rightarrow H^1(K,{\overline
G}_0)$, there is a ${\overline T}_0({\overline K})$-valued Galois
cocycle $\sigma\mapsto z_{j\: \sigma}$, $\sigma\in
\text{Gal}({\overline K}/K)$, representing $\xi_j$. Therefore, there
exists  a $\overline{K}$-isomorphism $\varphi_j \colon G_0 \to G_j$
such that $\varphi_j^{-1} \circ \sigma(\varphi_j) = i(z_{j\:
\sigma}),$ for all $\sigma \in \Ga(\overline{K}/K)$, where $i$ is
the natural isomorphism ${\overline G}_0\rightarrow
\text{Int}\,G_0$. Then $\varphi_j \vert T_0$ is defined over $K$,
and hence, $T_j^0 := \varphi_j(T_0)$ is a maximal $K$-torus of
$G_j$. Now $\varphi_0 := \varphi_2 \circ \varphi_1^{-1}$ is a
$\overline{K}$-isomorphism from $G_1$ onto $G_2$ whose restriction
to $T_1^0$ is an isomorphism onto $T_2^0$ defined over $K$. Since
$\beta_v(\alpha(\xi_1)) = \beta_v(\alpha(\xi_2)),$ the groups $G_1$
and $G_2$ are $K_v$-isomorphic, for all $v \in V^K.$ In addition,
for each $j = 1 , 2$, and any $v \in V,$ the group $G_j$ is
$K_v$-isomorphic to the group $_{\xi^{(v)}}G_0$ obtained from $G_0$
by twisting over $K_v$ by any cocycle representing
$\alpha_v(\xi^{(v)}).$ So, applying Corollary \ref{C:T10}, we obtain
the following.
\begin{thm}\label{T:T5}
Let $\xi(\varepsilon) \in H^1(K , \overline{G}_0),$ $\varepsilon \in
E_t,$ be the cohomology classes as in Proposition
\ref{P:T2}, and let $G_{\varepsilon}$ be the group obtained by
twisting $G_0$ by a cocycle representing $\xi(\varepsilon).$ Then
$G_{\varepsilon},$ $\varepsilon \in E_t,$ are pairwise nonisomorphic
$K$-forms of $G_0.$ Moreover, if for every $\varepsilon \in E_t$,
and every maximal $K$-torus $T$ of $G_{\varepsilon}$, we have {\brus
SH}$^2(T) = 0$ (which is automatically the case if for some $v \in
V$ the twist $_{\xi^{(v)}}G_0$ is $K_v$-anisotropic), then all the
groups $G_{\varepsilon}$ have coherently equivalent systems of
maximal $K$-tori.
\end{thm}

\noindent{\bf Remark 9.12.} If $G$ is an absolutely simple simply connected 
inner $K$-form of type $A_n$, then the
condition {\brus SH}$^2(T) = \{ 0 \}$ is automatically satisfied for
any maximal $K$-torus $T$ of $G.$ Indeed, $T$ is of the form $T =
R^{(1)}_{A/K}(\mathrm{GL}_1),$ where $A$ is a commutative \'etale
${(n+1)}$-dimensional $K$-algebra. Letting $S =
R_{A/K}(\mathrm{GL}_1),$ we have the exact sequence
$$
1 \to T \longrightarrow S \longrightarrow \mathrm{GL}_1 \to 1,
$$
which in conjunction with Hilbert's Theorem 90 induces the following
commutative diagram with exact rows:
$$
\begin{array}{ccccc}
0 & \longrightarrow & H^2(K , T) & \longrightarrow & H^2(K , S) \\
 &  & \downarrow & & \downarrow \\
0 & \longrightarrow & \bigoplus_v H^2(K_v , T) & \longrightarrow &
\bigoplus_v H^2(K_v , S).
\end{array}
$$
Since the map $H^2(K , S) \longrightarrow \bigoplus_v H^2(K_v , S)$
is injective by the Albert-Hasse-Brauer-Noether Theorem, our
assertion follows.

\vskip3mm

We observe that if $G_1$ and $G_2$ have coherently equivalent
systems of maximal $K$-tori, then for any finite set $S \subset
V^K_{\infty}$ containing $V^K_{\infty},$ any  $(G_i , K ,
S)$-arithmetic subgroups $\Gamma_i \subset G_i(K)$ are weakly
commensurable (see the argument in Example 6.5). It turns out that
in this situation arithmetic subgroups provide 
length-commensurable locally symmetric spaces.

\vskip3mm

\noindent{\bf Proposition 9.13.} {\it Let $G$ be a connected semi-simple
real algebraic group
and $\mathfrak{X}$ be
the symmetric space of $\mathcal{G} = G(\R).$ For $i = 1 , 2,$ let
$\Gamma_i$ be a torsion-free $(G_i , K)$-arithmetic
subgroup of $\mathcal{G}.$ If $G_1$ and $G_2$ have coherently
equivalent systems of maximal $K$-tori, then the locally symmetric
spaces $\mathfrak{X}_{\Gamma_1}$ and $\mathfrak{X}_{\Gamma_2}$ are
length-commensurable.}
\vskip2mm

\begin{proof} (Cf.\:the proof of Corollary \ref{C:G2007}.)
We can assume that $\Gamma_i \subset G_i(K)$ for $i = 1 , 2.$ Let
$\gamma_1 \in \Gamma_1$ be a nontrivial semi-simple element, and let
$T_1 \subset G_1$ be a maximal $K$-torus containing it. By our
assumption, there exists an isomorphism $\varphi \colon G_1 \to G_2$
such that the restriction $\varphi \vert T_1$ is defined over $K$,
hence $T_2 := \varphi(T_1)$ is a maximal $K$-torus of $G_2.$ Since
$\varphi(T_1(K) \cap \Gamma_1)$ is an arithmetic subgroup of $T_2(K),$
there exists $n > 0$ such that $\gamma_2 : = \varphi(\gamma_1)^n$
belongs to $\Gamma_2.$ The map $\alpha \to \alpha \circ \varphi$
defines a bijection between the root systems $\Phi(G_2 , T_2)$ and
$\Phi(G_1 , T_1).$ It follows that the sets of complex numbers
$$
\{ \alpha(\gamma_1^n) \ \vert\  \alpha \in \Phi(G_2 , T_2) \} \ \
\text{and} \ \ \{ \alpha(\gamma_2) \ \vert\  \alpha \in \Phi(G_1 ,
T_1) \}
$$
are identical. Using the formula (\ref{E:G-1}) from Proposition
\ref{P:G-1}(ii), we see that
\vskip1mm

\centerline{$\lambda_{\Gamma_2}(\gamma_2)/\lambda_{\Gamma_1}(\gamma_1) \in
\Q.$}\end{proof}

\noindent{\bf 9.14.} We finally indicate how Theorem \ref{T:T5} can
be used to construct examples of weakly commensurable cocompact arithmetic and
$S$-arithmetic subgroups, and length-commensurable compact locally symmetric
spaces, which are not commensurable. Let $G$ be a connected
absolutely simple simply connected isotropic real algebraic group of
one of the following types: $A_n$,  $D_{2n+1}$, $n>1$, or $E_6,$ and
let $\mathcal{L}$ be either $\R$ or $\C$ depending on whether or not
$G$ is an inner form over $\R.$ Fix a real quadratic extension
$K/\Q,$ and let $v'_{\infty} , v''_{\infty}$ denote its two real
places. Next, pick a quadratic extension $L/K$ so that $L\otimes_K
K_{v'_{\infty}} = \mathcal{L}^{2/[\mathcal{L}:\R]}$ and $L\otimes_K
K_{v''_{\infty}} = \C,$ and let $G_0$ denote the nonsplit
quasi-split $K$-group of the same type as $G$ which splits over $L.$
Since for the types under consideration, the $\R$-anisotropic form
is an inner twist of the corresponding nonsplit quasi-split $\R$-group, there
exist cohomology classes $\xi^{(v'_{\infty})} \in
H^1(K_{v'_{\infty}} , \overline{G}_0)$ and $\xi^{(v''_{\infty})} \in
H^1(K_{v''_{\infty}} , \overline{G}_0)$ such that the twist
$_{\xi^{(v'_{\infty})}}G_0$ is isomorphic to $G$ and the twist
$_{\xi^{(v''_{\infty})}}G_0$ is $\R$-anisotropic. Then applying the
construction described in Theorem \ref{T:T5} to $V = \{ v'_{\infty}
, v''_{\infty} \}$ and the specified cocycles, we obtain $2^t$
groups $G_{\varepsilon},$ $\varepsilon \in E_t,$ which are pairwise
nonisomorphic over $K$ but have coherently equivalent systems of
maximal $K$-tori as these groups are all anisotropic over
$K_{v''_{\infty}}.$ Besides, $G_{\varepsilon}$ is isomorphic to $G$
over $K_{v'_{\infty}} = \R,$ for every $\varepsilon \in E_t.$ Thus,
torsion-free arithmetic subgroups of $G_{\varepsilon}$ yield
discrete torsion-free subgroups of $\mathcal{G} = G(\R),$ and it
follows from Proposition 9.13 that the resulting locally symmetric
spaces are length-commensurable, but not commensurable. Finally, for
any finite subset $S$ of $V^K$ containing $V^K_{\infty},$ the
$S$-arithmetic subgroups of $G_{\varepsilon},$ $\varepsilon \in
E_t,$ are weakly commensurable, but not commensurable (cf.\:Example
6.5).

\vskip2mm

\noindent {\bf Remark 9.15.} Most of the results of this section
immediately extend to a global function field $K.$ This applies, in
particular, to Theorem \ref{T:T2}, yielding a local-global
principle for the existence of a coherent embedding, and Theorem
\ref{T:T5}, containing a construction of forms of a quasi-split
group $G_0$ belonging to one of the types $A_n,$ $D_{2n+1}$ $(n >
1)$ or $E_6,$ which are not $K$-isomorphic, but are isomorphic over
$K_v$ for all $v \in V^K.$ It should be noted, however, that the
construction of nonisomorphic $K$-groups with coherently equivalent
systems of maximal $K$-tori, described in 9.14, extends to global
function fields only for groups of type $A_n.$ The reason is that we
ensured the triviality of {\brus SH}$^2(T)$ for all maximal tori of
a group under consideration by arranging that the group is
anisotropic at a certain archimedean place. Over global function
fields, however, any group of type different from $A_n,$ is
isotropic.

\vskip5mm

\section{Isospectral locally symmetric spaces}\label{S:Sp}

The following theorem is known. For locally symmetric spaces of rank
$1$, a proof is given in \cite{Ga}. However, for locally symmetric
spaces of rank $>1$, we have not been able to find a
reference for it. For the convenience of the reader we will give
below its proof which was supplied to us by Alejandro Uribe and
Steve Zelditch.

\begin{thm}\label{T:S1}
Let $M_1$ and $M_2$ be two compact locally symmetric spaces with
nonpositive sectional curvatures. Assume $M_1$ and $M_2$ are
isospectral, in the sense that the spectra of their Laplace-Beltrami
operators on functions are the same (their eigenvalues and their
multiplicities). Then the sets
$$
L(M_j) = \{\,\lambda \in\R\;;\; \text{there exists a periodic
geodesic in } M_j\ \text{of length}\ \lambda\,\},
$$
for $j=1$, $2$, are equal.
\end{thm}

As we will explain, this theorem is a direct consequence of theorems
of Duistermaat and Guillemin, \cite{DG}, and of Duistermaat, Kolk
and Varadarajan, \cite{DKV}.  (In fact, the results of the latter
paper alone imply this theorem, but it is conceptually better to
use the main theorem of \cite{DG} in the proof.)

The results of \cite{DKV} (cf.\:Proposition 5.15) include that, for
$M$ a compact locally symmetric space of non-compact type,
\begin{enumerate}
\item [(i)] $L(M)$ is a discrete subset of $\R$, and
\item[(ii)] if $\lambda\in L(M)$, the set
$$
Z_{\lambda}:= \{ \overline{x}\in T^1M\;;\; \text{the geodesic
through }\overline{x}\ \text{is periodic of length}\ \lambda \}
$$
is a finite union of closed submanifolds (possibly of different
dimensions) of the unit tangent bundle $T^1M$ of $M$.
\end{enumerate}
Denote by $Z_{\lambda}^\circ$ the union of connected components of
$Z_{\lambda}$ of maximal dimension. It turns out that, in addition
to the previous theorem, for $M$ as above
\begin{equation}\label{1}
\text {for \ \ all\ \ } \lambda\in L(M)\quad \text{dim }
Z_{\lambda}^\circ\ \text{and Vol }Z_{\lambda}^\circ\ \text{are
spectrally determined.}
\end{equation}
Here the volume is with respect to a measure naturally induced by
the geodesic flow. (Equation (5.47) of \cite{DKV} is a formula for
this volume.)

\bigskip

Let us now see how one proves Theorem \ref{T:S1} and the additional
statement, (\ref{1}). Proposition 5.8 of \cite{DKV} establishes that
each $Z_{\lambda}$ is a clean fixed-point set of the time $\lambda$
map of the geodesic flow.  We can therefore apply the
Duistermaat-Guillemin trace formula, \cite{DG}, to the square root
of the Laplace-Beltrami operator on $M$. Specifically, pick a length
$\lambda$ and a Schwartz function on the real line, $\varphi$, such
that its Fourier transform $\widehat\varphi$ is compactly supported
and satisfies:
\[
\widehat\varphi(\lambda )=1\quad\text{and}\quad L(M)\cap\text{supp
}\widehat\varphi = \{\lambda\}.
\]
(Such a $\varphi$ exists by item (i) above.) Let
$0=\mu_0<\mu_1\leqslant\mu_2\leqslant\cdots$ be the square roots of
the eigenvalues of the Laplacian on $M$, listed with their
multiplicities.  Then, by Theorem 4.5 of \cite{DG} one has an
asymptotic expansion as $\mu\to\infty$ of the form:
\begin{equation}\label{2}
\sum_j\varphi(\mu-\mu_j)\sim e^{i\mu \lambda}\sum_{j=0}^\infty
c_j\,\mu^{d_{\lambda}-j}.
\end{equation}
Here $d_{\lambda} =(\text{dim }Z_{\lambda}^\circ-1)/2$.  A key point
is that the  leading coefficient, $c_0$, is not zero because the
Maslov indices (the integers $\sigma_j$ in equation (4.7) in
\cite{DG}) of all periodic geodesics on $M$ are zero, by Proposition
5.15 of \cite{DKV}.  By equation (4.8) in \cite{DG}, $c_0$ is equal
to the volume of $Z_{\lambda}^\circ$ times a factor that depends
only on $d_{\lambda}$. The expansion (\ref{2}) in the present
context is explicitly discussed in \S 5.6 of \cite{DKV}  (see the
last formula in that section which, incidentally, contains a typo: a
$\tau$ is missing in the left-hand side exponent).  The dimension of
$Z_{\lambda}^\circ$ is determined spectrally by the size in $\mu$ of
the left-hand side of (\ref{2}), and therefore $c_0$ determines the
volume of $Z_{\lambda}^\circ$.

\smallskip
Theorem \ref{T:S1} and statement (\ref{1}) follow from
(\ref{2}), the information on $c_0$, and the basic fact that if
$L(M)\cap\text{supp }\widehat\varphi = \emptyset$, then the
left-hand side of (\ref{2}) is $O(\mu^{-\infty})$.  By considering
all possible test functions $\varphi$ as above, one can detect the
set $L(M)$ from the eigenvalues of the Laplacian.\hfill$\Box$

\vskip2mm

Let $\cG$ be a connected semi-simple real Lie group of adjoint type
without compact factors, and  $\fX$ be the symmetric space of $\cG$. 
Let $\Gamma_1$ and $\Gamma_2$ be two
torsion-free irreducible {\it cocompact} discrete subgroups of $\cG$, and for 
$i=1,\,2$, $ \fX_{\Gamma_i}= {\fX}/\Gamma_i$ be the corresponding locally symmetric
spaces. From Theorems \ref{T:S1} and \ref{T:G-3} we obtain the
following.

\begin{thm}\label{T:S2} If $\fX_{\Gamma_1}$ and $\fX_{\Gamma_2}$ are isospectral, then
$\Gamma_1$ and $\Gamma_2$ are weakly commensurable.

\end{thm}

We now assume that $\cG$ is {\it absolutely} simple. Then using
Theorem \ref{T:S2} in conjunction with Theorem F, we obtain the
following result.
\begin{thm}\label{T:S4}
If $\fX_{\Gamma_1}$ and $\fX_{\Gamma_2}$ are isospectral, and $\Gamma_1$ is arithmetic,
then so is $\Gamma_2.$
\end{thm}
Theorem \ref{T:S2} combined with Theorem C yields the following.
\begin{thm}\label{T:S3}
Any two arithmetically defined compact isospectral locally symmetric
spaces of an absolutely simple real Lie group of type other than
$A_n$ $(n > 1),$ $D_{2n+1}$ $(n\geqslant 1)$, $D_4$ and $E_6$, are commensurable
to each other.
\end{thm}
\vskip2mm

The following remark is due to Peter Sarnak.
\vskip1mm

\noindent{\bf Remark 10.5.} It was proved by Hermann Weyl that any
two isospectral Riemannian manifolds are of same volume (and of same
dimension), see, for example, \cite{Gi}, Theorem 4.2.1. Now, as
before, let $\cG$ be a connected semi-simple real Lie group of
adjoint type without compact factors, and $\mathfrak{X}$ be its
symmetric space. If $\Gamma$ is a torsion-free irreducible cocompact
discrete subgroup of $\cG$, then the set of conjugacy classes of
torsion-free irreducible cocompact discrete subgroups $\Gamma'$ of
$\cG$ such that $\mathfrak{X}/\Gamma'$ is isospectral to
$\mathfrak{X}/\Gamma$ is finite. This follows from H.C.\,Wang's
finiteness theorem (\cite{Ra}, Ch.\:IX) if $\cG$ is not isomorphic
to $\mathrm{PSL}_2(\mathbb{R})$, since according to a thereom of
Andr\'e Weil (\cite{Ra}, Theorem 7.63) cocompact irreducible
discrete subgroups in such a $\cG$ are locally rigid, and
$\mathfrak{X}/\Gamma$ and $\mathfrak{X}/\Gamma'$, and therefore,
$\cG /\Gamma$ and $\cG /\Gamma'$ have equal volume. On the other
hand, if $\cG$ is isomorphic to $\mathrm{PSL}_2(\mathbb{R})$, then
the finiteness of the conjugacy classes of $\Gamma'$s is proved in
\S 5.3 of \cite{Mc}.
\vskip1mm

\vskip0.5cm

\bibliographystyle{amsplain}

\end{document}